\documentclass[11pt]{article}
\usepackage[T1]{fontenc}
\usepackage[utf8]{inputenc}
\usepackage{lmodern}
\usepackage[margin=1in]{geometry}
\usepackage{amsfonts,amsmath,amssymb,amsthm,mathtools,mathrsfs,times,color}
\usepackage{microtype}
\usepackage{enumitem,booktabs,array}
\usepackage[colorlinks,linkcolor=blue,anchorcolor=blue,citecolor=blue]{hyperref}
\numberwithin{equation}{section}
\newtheorem{theorem}{Theorem}[section]
\newtheorem{lemma}[theorem]{Lemma}
\newtheorem{proposition}[theorem]{Proposition}
\newtheorem{corollary}[theorem]{Corollary}

\theoremstyle{definition}
\newtheorem{definition}[theorem]{Definition}
\newtheorem{remark}[theorem]{Remark}
\newcommand{\ddc}{\sqrt{-1}\partial\bar\partial}
\newcommand{\Ca}{\mathcal C}
\newcommand{\DF}{\operatorname{DF}}
\newcommand{\Aut}{\operatorname{Aut}}
\newcommand{\Vol}{\operatorname{Vol}}

\newcommand{\reg}{\mathrm{reg}}
\newcommand{\NA}{\mathrm{NA}}
\newcommand{\C}{\mathbb C}
\newcommand{\R}{\mathbb R}
\newcommand{\Q}{\mathbb Q}
\newcommand{\Z}{\mathbb Z}
\DeclareMathOperator{\Lie}{Lie}
\setlist[enumerate]{itemsep=4pt,topsep=5pt}

\hypersetup{pdftitle={Calabi flow in the controlled K\"ahler cone},
 pdfauthor={Haozhao Li, Bing Wang, Kai Zheng},
 pdfsubject={}}
\title{Calabi flow in the controlled K\"ahler cone}
\author{Haozhao Li\footnote{Supported by NSFC grant No. 12471058 and a
research fund from University of Science and Technology of China.},
\quad Bing Wang\footnote{Supported by NSFC grant No. 12431003 and a
research fund from University of Science and Technology of China.}
\quad and \quad Kai Zheng\footnote{Supported by NSFC grant No. 12571091.}}
\date{}
\begin{document}
\maketitle
\begin{abstract}
We prove that the geometric limit of a global Calabi flow in a real
class of the controlled K\"ahler cone is independent of the initial
metric. Under cohomological
K-polystability, the limit is a smooth cscK metric in the original
class and complex structure. We also establish varying-class
closedness and an equivariant continuity method. We apply this method
to construct a smoothly varying family of cscK metrics throughout an
explicit interval of symmetric real K\"ahler classes on the blow-up
of the projective plane at four points in general position.
\end{abstract}
\begingroup
\small
\setcounter{tocdepth}{2}
\tableofcontents
\endgroup

\section{Introduction}\label{app:sec:intro}
This is the fourth paper in our series on
Calabi flow on K\"ahler surfaces. The first paper \cite{LongTime} proves a
Sobolev extension criterion. It also gives
long-time existence for initial metrics with Calabi energy below the prescribed threshold
in the controlled K\"ahler cone. The second paper \cite{Compactness} proves
subsequential convergence to extremal orbifold flows under a
uniform Sobolev bound. The third paper \cite{Polarized} proves
Bergman estimates and uniqueness of the flow limit on polarized
surfaces. It also constructs test configurations.
Under K-polystability, the limit is a constant scalar curvature
K\"ahler (cscK) metric on the original surface.

The purpose of this paper is twofold. First, we prove uniqueness
of the extremal orbifold limit for a fixed complex surface and
real K\"ahler class. This limit is independent of the admissible
initial metric. Cohomological K-polystability identifies it with
a smooth cscK metric on the original surface. Second, we prove
closedness as the class varies and use a continuity method to
construct cscK metrics in real families. This step uses a separate
compactness argument for sequences of cscK metrics.
The precise hypotheses are as follows.

Let $M$ be a compact connected K\"ahler surface, put
$c\coloneqq c_1(M)$, and let $\mathcal K_M$ be its real K\"ahler
cone. Our conventions are
\begin{equation}\label{app:eq:conventions}
 [\omega(t)]=2\pi\Omega,\qquad
 \partial_t\omega=\ddc S,\qquad dV_t\coloneqq\omega(t)^2/2,
\end{equation}
where $S$ is complex scalar curvature, equal to one half of the
Riemannian scalar curvature. In the polarized case $\Omega=c_1(L)$,
so these conventions agree with those of the polarized paper.
Write
\[
 V\coloneqq\int_MdV_t=2\pi^2\Omega^2,\qquad
 \bar S\coloneqq V^{-1}\int_M S\,dV_t
          =\frac{2c\cdot\Omega}{\Omega^2}.
\]
The volume and average scalar curvature are constant.
The Calabi energy and its volume-normalized form are
\[
 \Ca(\omega)\coloneqq\int_M(S_\omega-\bar S)^2dV_\omega,
 \qquad \widehat{\Ca}(\omega)\coloneqq\frac{\Ca(\omega)}{V}.
\]
Along the flow we write $\Ca(t)\coloneqq\Ca(\omega(t))$.
The controlled cone and the initial-energy hypothesis are
\begin{equation}\label{vc:eq:cone}
 \mathcal U\coloneqq
 \left\{\Omega\in\mathcal K_M:\ c\cdot\Omega>0,\quad
 \delta(\Omega)\coloneqq
 c^2-\frac{2(c\cdot\Omega)^2}{3\Omega^2}>0\right\},
\end{equation}
\begin{equation}\label{ru:eq:barrier}
 \omega_0\in2\pi\Omega,\qquad
 \Ca(\omega_0)<8\pi^2\delta(\Omega).
\end{equation}
For $\Omega\in\mathcal U$, the long-time theorem
\cite[Corollary~1.3]{LongTime} gives a global smooth Calabi flow
from such an initial metric and a uniform Sobolev bound
\begin{equation}\label{tc:eq:sob}
 \sup_{t\geq0}C_S(\omega(t))\leq C_S<\infty.
\end{equation}
Here $C_S(\omega)$ is the least constant such that
\[
 \|u\|_{4}^2\leq C_S(\omega)
 \left(\|\nabla u\|_{2}^2+V^{-1/2}\|u\|_{2}^2\right)
\]
for all $u\in W^{1,2}(M,g)$, with norms taken using $dV_\omega$.

The cone has useful algebraic consequences. The inequality
$K_M\cdot\Omega<0$ makes all plurigenera vanish, so the
Enriques--Kodaira classification makes $M$ rational or ruled.
The inequality $c^2>0$ excludes ruled minimal models over curves
of genus $g\geq1$, whose first Chern number after $r$ blow-ups is
$8(1-g)-r\leq0$. Thus $M$ is rational and projective, with $q(M)=p_g(M)=0$.
In particular, the controlled-cone condition supplies the algebraic input
\begin{equation}\label{app:eq:algebraic-consequences}
 \chi(\mathcal O_M)=1,\qquad c_1(M)\cdot\Omega>0.
\end{equation}
These conditions prevent a jump in the limiting section spaces and
allow us to realize the geometric limits in fixed projective spaces.
The Sobolev bound \eqref{tc:eq:sob} supplies the analytic input for
geometric compactness.

For real classes, stability is understood cohomologically \cite{DRKahler,SDPoly}.
A normal relatively K\"ahler cohomological test configuration
carries an invariant real $(1,1)$-class restricting to $\Omega$
on the general fibre. Cohomological K-semistability requires its
Donaldson--Futaki invariant to be nonnegative for every such
configuration. K-polystability additionally requires equality
only for products, allowing a nontrivial $\C^*$-action on $M$.
Subsection~\ref{real:sec:stability-definitions} gives the precise
definitions and our normalization of the invariant.

\begin{theorem}[Uniqueness, stability and continuation in the controlled cone]
\label{ru:thm:unique}\label{rd:thm:convergence}
Let $\Omega\in\mathcal U$ be a real K\"ahler class, and let
$\omega_0\in2\pi\Omega$ be a smooth K\"ahler metric satisfying
\eqref{ru:eq:barrier}. Its Calabi flow exists smoothly for all
time and satisfies \eqref{tc:eq:sob}. It has the following
properties.
\begin{enumerate}[label=\textnormal{(\arabic*)}]
\item \textbf{\upshape Uniqueness of the limit.}
All subsequential limits supplied by
Theorem~\ref{app:thm:compactness} are holomorphically isometric
to one compact extremal K\"ahler orbifold
$(Z,J_Z,\omega_Z)$. More precisely, any two such limits admit
an orbifold biholomorphism $F:Z_1\to Z_2$ with
$F^*\omega_2=\omega_1$. The whole flow converges to this
holomorphic-isometry class in the geometric sense of
Definition~\ref{app:def:geometric-convergence}. This includes
convergence on every bounded interval of translated time, using
comparison maps independent of that time. This holomorphic-isometry
class depends only on
$(M,J,\Omega)$: every other smooth initial metric in $2\pi\Omega$
satisfying \eqref{ru:eq:barrier} has the same limit.
\item \textbf{\upshape Semistability.}
If $(M,\Omega)$ is cohomologically K-semistable, then
$\omega_Z$ is cscK. The limiting orbifold may still be singular.
\item \textbf{\upshape Polystability and the original class.}
If $(M,\Omega)$ is cohomologically K-polystable, then $Z$ is
smooth and there is a biholomorphism $F:M\to Z$ such that
\[
 F^*\Omega_Z=\Omega,\qquad
 \Omega_Z\coloneqq[\omega_Z]/(2\pi).
\]
Consequently $F^*\omega_Z$ is a smooth cscK metric in
$2\pi\Omega$ on the original complex surface, and the whole
flow converges to its holomorphic-isometry class as in part~(1).
\item \textbf{\upshape Continuation in the real cone.}
Suppose additionally that $\Aut^0(M)=\{1\}$. Let
$D\subset\mathcal U$ be a connected set containing $\Omega$,
and suppose $(M,B)$ is cohomologically K-polystable for every
$B\in D$. Then every class $2\pi B$, $B\in D$, contains a
unique smooth cscK metric, depending locally smoothly on $B$.
\end{enumerate}
For parts~(2)--(3), it suffices to impose the sign condition,
and respectively the product equality condition, only on normal
relatively K\"ahler cohomological test configurations whose
central fibre is normal, reduced and irreducible.
\end{theorem}

Part~(1) requires neither rationality nor stability of the class.
It also allows positive-dimensional automorphism groups.

Parts~(2)--(3) require no uniform positive stability bound.
The stated restricted stability conditions suffice because the proof
constructs a test configuration with the metric limit as its
central fibre and with the required Donaldson--Futaki sign.

Theorem~\ref{real:thm:poly-continuation} gives the precise closedness
and local smooth dependence statements underlying part~(4).
Here cohomological K-polystability is assumed for each real class,
so no quantitative stability gap is required.
An alternative continuation criterion uses a local restricted
stability gap in nearby rational auxiliary classes
(Theorem~\ref{real:thm:existence}). The neighborhood and gap constant
may depend on the real class. Its finite-group equivariant form is
used in the application below.

\medskip
\noindent\textbf{Outline of proof.}
The compactness theorem of \cite{Compactness} supplies subsequential
extremal orbifold limits, whose complex structures may a priori differ.
To compare them, we seek a common finite-dimensional description that
retains the metric class. When $\Omega$ is irrational, no line bundle
has first Chern class $\Omega$. We instead use finitely many auxiliary
bundles and their complete section spaces.

\medskip
\noindent\emph{Cohomological marking retains the metric class on each limit.}
Let $W_\Omega$ be the smallest rational real subspace of $H^2(M,\R)$
containing $\Omega$, and put $d\coloneqq\dim_\R W_\Omega$.
Harmonic representatives of an integral basis
have controlled $L^2$ norms and converge on regular regions.
The discreteness of the exceptional set is essential here.
The limiting forms extend across punctures in uniformizing four-balls,
and their connections extend after local gauge changes, giving
orbifold line bundles. The minimality of $W_\Omega$ makes the
coefficients of $\Omega$ in this basis linearly independent over $\Q$.
Together with rationality of orbifold intersection numbers, this turns
preservation of degrees against $\Omega$ into preservation of the
individual intersection pairings. The resulting marking retains the
numerical data needed for the Hilbert polynomials.

\medskip
\noindent\emph{Common auxiliary embeddings allow all subsequences to be compared.}
The bundles and embedding powers must be fixed in advance.
On a surface, adjunction and the Hodge index theorem bound the degrees
of the marked classes on every curve by its metric area. The Sobolev
bound controls the local orbifold group orders and Cartier denominators.
These facts give one neighborhood of $\Omega$ whose rational classes
are ample on every marked limit. Choose finitely many such classes
whose positive cone contains $\Omega$. Clearing denominators gives
line bundles with fixed Hilbert polynomials, and Matsusaka's theorem
\cite[Theorem~2.4]{Matsusaka} gives common very ample powers.
We obtain fixed bundles and positive coefficients satisfying
\[
 \Omega=\sum_\nu\alpha_\nu c_1(A_\nu),\qquad \alpha_\nu>0.
\]
Complete-section convergence identifies each geometric limit with its
scheme-theoretic Hilbert limit in the product of these embeddings.

\medskip
\noindent\emph{Harmonic Chern metrics give convergence of the Gram matrices.}
On each $A_\nu$, choose a Hermitian fibre metric $h_\nu(t)$ whose real
Chern curvature is the $\omega(t)$-harmonic representative of
$2\pi c_1(A_\nu)$. We call this a \emph{harmonic Chern metric}.
Harmonicity refers to its curvature $-\ddc\log h_\nu(t)$, computed
in a local holomorphic frame. This fibre metric is distinct from the
K\"ahler metric $\omega(t)$ on the surface. Its curvature has constant
trace by the K\"ahler identities, although it need not be positive.
Together with the Sobolev bound, this gives upper bounds for normalized
sections. Vanishing omitted volumes prevent loss of section mass.
On regular compact sets, nonvanishing of the complete limiting linear
systems and a Harnack argument control ratios of the fibre metrics.
These estimates give convergence of the Gram matrices of the varying
$L^2$ inner products in fixed section bases. On translated time
intervals, their limits record the extremal motion as matrix
exponentials, extending \cite[Lemma~4.1]{Polarized}.

\medskip
\noindent\emph{Countable spectra and orbit uniqueness identify the limit of one flow.}
The extremal spectra need not be rational in a real class.
Equivariant orbifold Riemann--Roch \cite{Kawasaki,Meinrenken}, applied
to tensor products of the auxiliary bundles, expresses the Hamiltonian
pairings and Futaki invariants through rational mixed coefficients.
Evaluation at the fixed $\alpha_\nu$ and the linear system for the
extremal field place every possible spectrum in the countable field
$\mathbb F=\Q(\alpha_1,\ldots,\alpha_d,\pi)$.
The ordered spectral data of relative Gram matrices vary continuously
with time and have a compact connected cluster set. Since a connected
subset of a countable subset of Euclidean space is a single point,
all limiting generators have the same spectra. The
Chen--Sun--Wang orbit argument \cite[Section~2]{CSW}, together with
fixed-orbifold extremal uniqueness in Theorem~\ref{un:thm:fixed},
then identifies the metric limits up to holomorphic isometry.

\medskip
\noindent\emph{Relative Gram matrices compare different initial metrics.}
For two flows in the same class, Mabuchi-distance contraction
\cite{CalabiChen,BDLFlow} and Darvas's endpoint estimates
\cite[Theorem~3]{DarvasFinite} bound the integral of the difference
of their normalized potentials against both evolving volume forms.
The positive expression for $\Omega$ relates this difference to the
logarithmic ratios of the auxiliary fibre metrics, up to a spatial
constant. On regular balls where selected limiting sections do not
vanish, section upper bounds control the relative Gram eigenvalues.
Reversing the flows gives a second inequality. Adding the inequalities
cancels the unknown constant, and positivity of every $\alpha_\nu$
bounds each relative condition number. After determinant normalization,
the changes of section bases and their inverses remain bounded.
Their limits identify the two Hilbert limits and all auxiliary bundles.
The same comparison on all positive and negative translated limiting
times forces the limiting change of basis to intertwine the extremal
generators. Fixed-orbifold uniqueness then gives the same metric limit.
Subsection~\ref{ru:sec:two-flow} contains this comparison.

\medskip
\noindent\emph{A degeneration preserving the real class makes stability applicable.}
The real extremal generator need not generate an algebraic
one-parameter subgroup. For degeneration towards zero, approximate
its negative by a rational generator in the same chamber of the
finitely many relevant weights. This preserves the intermediate
Hilbert limit and, for a non-cscK limit, the strict negative Futaki sign.
Donaldson's accessibility argument
\cite[Proposition~1(2), Section~5.2]{Donaldson} supplies a second
degeneration in the centralizer reaching the metric limit.
The two subgroups commute, so a sufficiently large positive power
of the first, multiplied by the second, reaches this same limit
\cite[Lemma~4.5]{Polarized}. The universal Hilbert family retains all
auxiliary bundles. Their positive combination
$\sum_\nu\alpha_\nu c_1(\mathcal L_\nu)$ restricts to exactly
$\Omega$ and its marked metric class on the general and central fibres.
The coefficients remain fixed when the generator is replaced.
Rational approximation of the coefficients is used only to verify
the Donaldson--Futaki formula on this same family
(Proposition~\ref{rd:prop:degeneration}). The invariant is negative
for a non-cscK limit and zero for a cscK limit. Semistability therefore
forces the limit to be cscK, and polystability makes the configuration
a product, identifying the limit with a smooth cscK metric in the
original class and complex structure.

\medskip
\noindent\emph{Static compactness and polystability give closedness for continuation.}
When the class varies, a separate argument is needed.
For cscK metrics in classes approaching $\Omega_0$ inside the
controlled cone, static compactness \cite{TV,AndersonOrbifold} gives
an orbifold limit. The marking and auxiliary embeddings realize it in
a degeneration of the exact pair $(M,\Omega_0)$. Its invariant is
zero because the central metric is cscK, so polystability forces a
product and yields closedness. For discrete automorphism groups,
the linearized scalar-curvature operator is invertible on normalized
potentials. The implicit-function theorem gives openness and local
smooth dependence, as in LeBrun--Simanca's deformation theory
\cite{LeBrunSimanca}. Part~(3) supplies the starting cscK metric in
part~(4). The cscK locus is consequently nonempty, open and closed
relative to $D$, so connectedness gives existence throughout $D$.

\medskip
\noindent\textbf{Application.}
We apply this continuity method to
$X=\operatorname{Bl}_{p_1,\ldots,p_4}\mathbb P^2$,
where the points are in general position, using the equivariant
criterion of Theorem~\ref{real:thm:existence}.
With $H$ the pullback of
the line class and $E\coloneqq E_1+\cdots+E_4$, set $A_a\coloneqq aH-E$.
Section~\ref{cone:sec:fourpoint} proves that every real class
\[
 2\pi A_a,\qquad 3\leq a<8+\frac{10}{\sqrt3},
\]
contains a unique smooth $S_4$-invariant cscK metric, varying smoothly
with $a$. The anticanonical K\"ahler--Einstein metric at $a=3$
starts the argument. Proposition~\ref{cone:prop:gap} verifies the
local restricted equivariant gap, independently of the cone
inequalities and including the two-dimensional neighborhood of an
irrational $A_a$. This gap excludes nonproduct zero-invariant
degenerations and provides the required closedness. In particular,
we obtain a continuation proof
for $A_{12}$. Admissible nonstationary initial metrics are then
obtained by small perturbations of the constructed cscK family.
Section~\ref{cone:sec:fourpoint} compares this continuation result
with the anticanonical existence theorem of Tian--Yau
\cite{TianYauKE} and the gluing constructions of Arezzo--Pacard
\cite{APBlowupII}.

\medskip
\noindent\textbf{Relations to earlier work.}
The polarized companion \cite{Polarized} supplies the complete-section
convergence, fixed-orbifold uniqueness and degeneration arguments used here.
Its polarization both determines the metric class and supplies the
section spaces, and its convergence theorem concerns a fixed flow.
For real classes these roles separate. An auxiliary embedding alone
does not specify the original class, and weak convergence gives no
canonical identification of the cohomology of the source and limit.
The marking records the additional data needed to preserve that class.
It is defined on $W_\Omega$ and may depend on the extraction.
The rational-spectrum argument is also replaced by a countability
argument. These extensions allow stability to be applied to the exact
real class in the cohomological framework of Dervan--Ross and
Sj\"ostr\"om Dyrefelt \cite{DRKahler,SDPoly}.

Our comparison of different initial metrics follows
Chen--Sun--Wang \cite[Section~4, concluding Remark~(1), p.~3170]{CSW}.
For two global flows in one polarized class with uniformly bounded
curvature and diameter, they combine Mabuchi-distance contraction
with Darvas's potential estimates and use bounded geometry to obtain
higher regularity of the potential difference. We keep this
metric-space starting point. Under the weak compactness of
\cite{Compactness}, the positive auxiliary decomposition instead
converts the endpoint integral estimate into relative Gram bounds.
This permits comparison when the metric class is real and the limits
have orbifold singularities. The argument is pairwise and requires
no connected family of admissible initial metrics.

The convergence statements concern holomorphic-isometry classes,
with subsequence-dependent geometric identifications. They do not
assert smooth convergence in a fixed coordinate gauge or uniqueness
of the test configuration. The combined degeneration retains the
metric limit and the real class, but need not retain the
flow-induced filtration or an optimality property.

For a fixed invariant flow, the convergence and degeneration
constructions admit finite-group equivariant versions. In that
setting it suffices to
test equivariant configurations with normal, reduced and irreducible
central fibre. The precise statement, including the condition for an
equivariant identification with the original surface, is
Theorem~\ref{eqv:thm:real}.

\medskip
\noindent\textbf{Organization of the paper.}
Section~\ref{app:sec:preliminaries}
recalls the geometric inputs and fixes cohomological stability.
Section~\ref{ru:sec:marking-construction} constructs the marked
auxiliary classes and common embeddings. Section~\ref{ru:app:proof}
proves limit uniqueness and the class-preserving degeneration.
Section~\ref{eqv:sec:extensions} establishes the equivariant version.
Section~\ref{real:sec:existence} first proves closedness and
continuation under cohomological K-polystability, then establishes
the local restricted gap criterion. Section~\ref{cone:sec:fourpoint} proves
the restricted stability gap and applies the criterion to the
symmetric four-point blow-up and the class $A_{12}$.

\medskip
\noindent\textbf{On the usage of AI.}
The authors used ChatGPT as an assisting tool in the development of
the paper. The authors take full responsibility for the contents of
the paper.

\par\addvspace{\baselineskip}
\noindent\textbf{Acknowledgement.}
The authors are deeply grateful to Professor Xiuxiong Chen
for his encouragement and constant support.
The authors also thank Professor Song Sun for his interest in this work.

\section{Preliminaries}\label{app:sec:preliminaries}

We first collect the principal notation, symbols and terminology.
We recall the compactness input from \cite{Compactness} and its
consequences from \cite{Polarized}, retaining the statements needed
below. For fixed-orbifold uniqueness we explain the extension to a
real metric class with an independent auxiliary polarization.
We then define cohomological test configurations and normalize their
Donaldson--Futaki invariant in arbitrary dimension. The subsequent
arguments concern complex surfaces.

\subsection{Notation, terminology and constants}
\label{cone:index:notation}
For the convenience of readers, we collect the principal notation,
symbols and terminology in the following table. Local notation is
specified where it occurs.

\noindent\emph{Flow and geometric limits.}
\begin{description}[style=nextline,leftmargin=2.7cm,labelwidth=2.4cm,labelsep=0.3cm,align=left,font=\normalfont,itemsep=0pt,parsep=0pt,topsep=4pt]
\item[$\Omega,c$] The real K\"ahler class $[\omega]/(2\pi)$ and
$c=c_1(M)$.
\item[$\mathcal U,\delta(\Omega)$] The controlled cone and its positive
numerical gap, defined in \eqref{vc:eq:cone}.
\item[$V,S,\bar S$] Total volume, complex scalar curvature and its average.
\item[$\Ca,\widehat{\Ca}$] The Calabi energy
$\Ca(\omega)=\int_M(S_\omega-\bar S)^2\,dV_\omega$ and its
volume-normalized form $\widehat{\Ca}=\Ca/V$.
The initial-energy barrier is \eqref{ru:eq:barrier}.
\item[$C_S(\omega)$] The scale-invariant Sobolev constant in
\eqref{tc:eq:sob}.
\item[cscK, GH] Constant scalar curvature K\"ahler, and
Gromov--Hausdorff.
\item[$Z,Z_{\reg}$] A limiting K\"ahler orbifold and its regular locus.
The same $Z$ denotes its coarse complex space when appropriate.
\item[Compactness property] Global GH convergence, regular smooth
convergence and volume convergence on translated time windows
(Definition~\ref{app:def:time-window}), with time-independent
comparisons. See Theorem~\ref{app:thm:compactness}.
\item[$F_\tau$] The limiting extremal motion
$\exp(\frac{\tau}{2}\nabla S_Z)$, satisfying
$\omega_Z(\tau)=F_\tau^*\omega_Z(0)$.
\end{description}

\par\addvspace{\baselineskip}
\noindent\emph{Markings and auxiliary polarizations.}
\begin{description}[style=nextline,leftmargin=2.7cm,labelwidth=2.4cm,labelsep=0.3cm,align=left,font=\normalfont,itemsep=0pt,parsep=0pt,topsep=4pt]
\item[$W,\beta_i,a_i$] The smallest rational real subspace containing
$\Omega$, an integral rational basis, and the coefficients
$\Omega=\sum_i a_i\beta_i$. See Subsection~\ref{ru:sec:marking}.
\item[$\mathcal B_i,\eta_i(t)$] Source line bundles with
$c_1(\mathcal B_i)=\beta_i$ and the harmonic representatives of $\beta_i$.
\item[$\mathfrak m_Z,B_Z$] The marking map and the marked image
$B_Z=\mathfrak m_Z(B)$. See \eqref{ru:eq:marking}.
\item[$\mathfrak L$] All extracted geometric limits equipped with
their cohomological markings (Lemma~\ref{ru:lem:marking}).
\item[$A_\nu,A_{\nu,Z}$] Fixed very ample auxiliary bundles and their
marked limiting bundles. Here $A_\nu=L_\nu^{r_\nu}$
(Lemma~\ref{ru:lem:embeddings}).
\item[$V_\nu$] The fixed section space $H^0(M,A_\nu)$. By contrast,
$V_0=\Omega^2$ denotes an intersection number in the marking argument.
\item[$\alpha_\nu$] Fixed positive coefficients expressing the original
class as $\Omega=\sum_\nu\alpha_\nu c_1(A_\nu)$. See
\eqref{rd:eq:positive-basis}.
\end{description}

\par\addvspace{\baselineskip}
\noindent\emph{Gram matrices and finite-dimensional dynamics.}
\begin{description}[style=nextline,leftmargin=2.7cm,labelwidth=2.4cm,labelsep=0.3cm,align=left,font=\normalfont,itemsep=0pt,parsep=0pt,topsep=4pt]
\item[$H_\nu(t)$] The section inner product
$\int_M h_\nu(t)(s,s')\,dV_t$. Its matrix in an ordered basis
is the Gram matrix.
\item[$\mathcal N(B)$] Determinant normalization
$(\det B)^{-1/N}B$ for a positive Hermitian matrix of size $N$.
\item[$C_{\nu,Z},\Lambda_\nu$] Limiting trace-free Hermitian generators
and their common unitary-conjugacy representatives. See
\eqref{ru:eq:gram} and Lemma~\ref{ru:lem:spectrum}.
\item[$\mathbb F$] The countable field
$\Q(\alpha_1,\ldots,\alpha_d,\pi)$ used to control possible spectra.
\item[$\mathbb G,K$] The frame groups
$\prod_\nu GL(V_\nu^*)$ and $\prod_\nu U(V_\nu^*)$.
Their centralizers of $\Lambda=(\Lambda_\nu)_\nu$ are
$G_\Lambda$ and $K_\Lambda$.
\item[{$[M],[\overline M],\mathscr C$}] The source Hilbert point,
the intermediate exponential limit, and the compatible Hilbert
limit set in the proof of Theorem~\ref{ru:thm:unique}(1).
\end{description}

\par\addvspace{\baselineskip}
\noindent\emph{Stability, equivariance and the example.}
\begin{description}[style=nextline,leftmargin=2.7cm,labelwidth=2.4cm,labelsep=0.3cm,align=left,font=\normalfont,itemsep=0pt,parsep=0pt,topsep=4pt]
\item[$G$] A finite group of holomorphic automorphisms. In
Section~\ref{cone:sec:fourpoint}, $G\cong S_4$.
\item[$(\mathcal X,\mathcal A)$] A normal cohomological test
configuration, carrying a real Bott--Chern $(1,1)$-class.
See Definition~\ref{real:def:test-configuration}.
\item[$\DF$] The Donaldson--Futaki invariant, normalized by
\eqref{rd:eq:cohom-df}.
\item[$K_{\mathcal Y/\mathbb P^1}$] The relative canonical divisor
$K_{\mathcal Y}-(\overline\pi\circ p)^*K_{\mathbb P^1}$
on a resolution of the compactified total space.
\item[Products and integral fibres] A product may carry a nontrivial action
on $M$. An integral fibre is reduced and irreducible.
\item[$J^{\NA}$] The non-Archimedean $J$-functional: the largest
asymptotic section weight minus the mean weight, with the actual
polarization exponent retained.
\item[$\mathcal S_G$] The classes satisfying the local restricted
equivariant stability gap \eqref{real:eq:local-gap}.
\item[$A_a,H,E$] On the four-point blow-up,
$A_a=aH-E$, where $H$ is the pullback of the line class and
$E=E_1+\cdots+E_4$.
\item[$J,a_+$] The parameter interval $J=[3,a_+)$, where
$a_+=8+10/\sqrt3$. See \eqref{cone:eq:symmetric-classes}.
\item[$\tau(a)$] The threshold for the volume-normalized Calabi energy,
$\tau(a)=4\delta(a)/(a^2-4)$. See \eqref{cone:eq:startup}.
\item[$\beta_A(F),j_A(F)$] The divisorial invariants in
\eqref{conealg:eq:beta}--\eqref{conealg:eq:j}, used to establish
the restricted gap of Proposition~\ref{cone:prop:gap}.
\end{description}

\subsection{Geometric compactness}\label{app:sec:geometric-compactness}

We use the orbifold conventions of \cite[Section~2.2]{Polarized}.
The local charts are $(\widetilde U,\Gamma)$, with
$\widetilde U\subset\C^2$ a ball and $\Gamma\subset U(2)$ finite
and free away from the origin.
The \emph{coarse space} is obtained by gluing the quotients
$\widetilde U/\Gamma$ with their invariant holomorphic functions.
We write $Z$ for either object when unambiguous: smooth forms are
orbifold forms, whereas sheaf cohomology and Hilbert schemes refer
to the coarse space. An orbifold line bundle is a compatible system
of equivariant line bundles on the charts. It descends to a
\emph{Cartier bundle}, meaning a line bundle on the coarse space,
exactly when the stabilizers act trivially on its fibres. A
sufficiently divisible tensor power always descends. The orbifold
canonical class agrees with the rational canonical class of the
coarse space, since there are no codimension-one stabilizers.

\begin{definition}[Translated time windows]
\label{app:def:time-window}
Let $g(t)$, $t\geq0$, be a family of metrics, and let
$t_j\to\infty$. For a fixed $T>0$, the interval
$[t_j-T,t_j+T]$ is called a \emph{time window centred at $t_j$}.
For sufficiently large $j$, this interval lies in $[0,\infty)$.
Introducing the translated time $\tau=t-t_j$, we regard the
restricted family as
\[
 g_j(\tau)\coloneqq g(t_j+\tau),
 \qquad \tau\in[-T,T].
\]
We call $[-T,T]$ the corresponding \emph{translated time window}.
The same convention applies to other time-dependent objects,
such as curvature forms, bundle metrics and Gram matrices.

Convergence on every bounded translated time window means
convergence on $[-T,T]$ for every fixed $T>0$, uniformly in
$\tau$ in the topology specified in the relevant statement.
Whenever such convergence is obtained by extracting a
subsequence, one subsequence is understood to work for every
fixed $T>0$. This can be arranged by diagonal extraction over
$T=1,2,\ldots$.
\end{definition}

\begin{theorem}[Geometric input from the compactness paper]\label{app:thm:compactness}\label{tc:hyp:wc}
The long-time compactness theorem of \cite[Theorem 1.1 and Definition 6.1]{Compactness}
has the following consequence. If the Calabi flow is smooth on
$[0,\infty)$ and satisfies \eqref{tc:eq:sob}, then any $t_j\to\infty$
has a subsequence converging on every bounded translated time interval
to an eternal extremal flow on one connected compact K\"ahler orbifold $Z$.
The orbifold has finitely many isolated quotient points, with effective
local groups free on the unit three-sphere. Global GH convergence is
realized by time-independent correspondences whose distortions tend to
zero uniformly on bounded time intervals. Convergence is mixed
$C^\infty$ on regular compact subsets under compatible embeddings
$\Phi_j$, independent of time, whose graphs lie in those correspondences.
These embeddings include collars and cover the corresponding regular
regions. The limiting metrics extend smoothly in orbifold charts and
\[
\omega_Z(\tau)=F_\tau^*\omega_Z(0),\qquad
F_\tau\coloneqq \exp\left(\tfrac\tau2\nabla S_Z\right).
\]
The natural volume measures converge and have total mass $V$.
\end{theorem}
We refer to this conclusion as the compactness property. It applies
to any fixed real K\"ahler class.

\begin{definition}[Geometric convergence modulo holomorphic isometry]
\label{app:def:geometric-convergence}
For the translated families of
Definition~\ref{app:def:time-window}, we use the spacetime
convergence recalled in Theorem~\ref{app:thm:compactness}. In particular, both the global correspondences and
the compatible embeddings of regular domains are independent
of $\tau$.

Two limiting families are identified if one time-independent
orbifold biholomorphism pulls back every metric of one family to
the corresponding metric of the other. We say that the whole
flow converges geometrically modulo holomorphic isometry if,
for every sequence $t_j\to\infty$ and every $T>0$, its translated
families converge on $[-T,T]$ in the preceding sense to a
representative of the same holomorphic-isometry class.
The comparison maps may depend on the sequence and on $T$, but
not on $\tau$ within that interval.
\end{definition}

We also record two consequences
proved in \cite[Lemmas~2.2--2.3]{Polarized}. Those proofs use only
the compactness property and the scalar $L^2$ bound
$\int_M S^2dV_t=\Ca(t)+\bar S^2V\leq\Ca(0)+\bar S^2V$.
They therefore apply here without a polarization and will be used
to pass scalar and section integrals to the limit.

\begin{lemma}[Uniform volume tails]\label{app:lem:tails}
For an exhaustion $K_\ell\Subset Z_{\reg}$ by smooth compact domains and
any bounded interval $I$, the compactness property gives
\begin{equation}\label{tc:eq:volume}
\lim_{\ell\to\infty}\limsup_{j\to\infty}\sup_{\tau\in I}
\Vol_{g(t_j+\tau)}(M\setminus\Phi_j(K_\ell))=0.
\end{equation}
\end{lemma}

\begin{lemma}\label{tc:lem:extremal}
Every limit supplied by the compactness property is extremal and
$\omega_Z(\tau)=F_\tau^*\omega_Z(0)$, where $F_\tau$ is a complete
one-parameter group of orbifold biholomorphisms. Moreover the scalar
average of the limit is $\bar S$.
\end{lemma}

The Sobolev bound also gives uniform local noncollapsing, by
\cite[Section~3.1, equation~(3.3)]{Polarized}, using the cutoff
argument of Y.~Wang \cite[Proposition~2.4, p.~1265]{Wang94}:
\begin{equation}\label{tc:eq:uniform-noncollapse}
 \Vol_{g(t)}B_{g(t)}(x,s)\geq\kappa s^4,
 \qquad \kappa\coloneqq(20C_S)^{-2},\qquad 0<s\leq V^{1/4}.
\end{equation}
This argument uses no polarization. The estimate passes to every
orbifold limit and bounds the local group orders: the volume density
at a quotient point is the Euclidean density divided by the group order.

\subsection{Uniqueness on a fixed orbifold}
\label{app:sec:orbifold-uniqueness}

For an auxiliary polarization $A$ and a holomorphic vector field $Y$,
write
\[
 \Aut(Z,A,Y)\coloneqq
 \{f\in\Aut(Z):f^*A\simeq A,\ f_*Y=Y\},
\]
where $\Aut(Z)$ consists of orbifold biholomorphisms. Its identity
component is denoted by $\Aut^0(Z,A,Y)$. The following real-class
form follows from the proof of \cite[Theorem~2.5]{Polarized}. We
explain why its auxiliary polarization need not represent the metric
class. For a compact K\"ahler orbifold $Z$, its irregularity is
$q(Z)\coloneqq\dim_{\mathbb C}H^1(Z,\mathcal O_Z)$.

\begin{theorem}[The centralizer form of extremal uniqueness]\label{un:thm:fixed}
Let $(Z,A)$ be a compact connected polarized K\"ahler orbifold with
isolated quotient points and $q(Z)=0$. Let $\omega_0,\omega_1$ be two
orbifold-smooth extremal metrics in the same real K\"ahler class.
The auxiliary polarization $A$ need not represent that class.
Suppose that their holomorphic extremal fields coincide:
\[
Y\coloneqq \nabla^{1,0}_{\omega_0}S_0
 =\nabla^{1,0}_{\omega_1}S_1.
\]
Then some $a\in\Aut^0(Z,A,Y)$ satisfies $a^*\omega_1=\omega_0$.
\end{theorem}

\begin{proof}[Reduction to the fixed-orbifold proof]
We check the real-class extension of
\cite[Theorem~2.5 and Section~2.3.2]{Polarized}.
Since $q(Z)=0$, orbifold Hodge theory gives $H^{0,1}(Z)=0$.
Thus $\iota_X\omega_i$ has a smooth complex holomorphy potential
for every holomorphic vector field $X$, regardless of the
rationality of $[\omega_i]$. If $X$ commutes with $Y$, its
normalized potential can be chosen invariant under the closure
$T$ of the imaginary extremal flow.

Connected automorphisms preserve every real cohomology class and
also preserve $A$, because $\operatorname{Pic}^0(Z)=0$.
Hence $G\coloneqq\Aut^0(Z,A,Y)$ is the same centralizer used
in the fixed-orbifold proof. The Hamiltonian centralizer calculation
in \cite[Lemma~4.4]{Polarized} applies to either metric and makes
$G$ reductive, with its Hamiltonian isometry subgroups maximal compact.
The auxiliary bundle makes $G$ linear algebraic without restricting
the common metric class.

The weak-geodesic, modified Mabuchi-convexity and fixed-measure
estimates of \cite[Lemmas~2.6--2.8]{Polarized} are stated for an
arbitrary smooth K\"ahler background on the fixed orbifold.
With the preceding facts, the Berman--Berndtsson perturbation
argument \cite{BB}, as implemented in
\cite[Section~2.3.2]{Polarized}, applies in the common
real class. All correcting automorphisms belong to $G$, so it
gives $a^*\omega_1=\omega_0$ for some $a\in G$, as required.
\end{proof}

\subsection{Cohomological stability for real classes}
\label{real:sec:stability-definitions}

We fix the stability convention used in the main theorem.
We retain the normalization $[\omega]=2\pi\Omega$ and formulate
the definition on the compactified family. For a normal complex
space $Y$, the real Bott--Chern group $H^{1,1}_{\mathrm{BC}}(Y,\R)$
consists of real $(1,1)$-forms with local smooth $\ddc$-potentials,
modulo $\ddc f$ for globally defined smooth real functions $f$.
Smooth functions and potentials on $Y$ are understood through
local embeddings into complex manifolds. A class is K\"ahler if
it has a representative with local strictly plurisubharmonic
potentials. See \cite{SDPoly}.

\begin{definition}[Normal cohomological test configuration]
\label{real:def:test-configuration}
Let $M$ be a compact connected K\"ahler manifold and let
$\Omega\in\mathcal K_M$. A \emph{normal cohomological test
configuration} for $(M,\Omega)$ consists of compactified data
$(\overline{\mathcal X},\overline\pi,\rho,\Phi,
\overline{\mathcal A})$ satisfying the following conditions:
\begin{enumerate}[label=\textup{(\roman*)}]
\item $\overline{\mathcal X}$ is a compact normal K\"ahler complex
space, and $\overline\pi:\overline{\mathcal X}\to\mathbb P^1$
is a flat surjective holomorphic map. Normality means that each
local analytic ring $\mathcal O_{\overline{\mathcal X},x}$ is an
integrally closed domain in its field of fractions.
\item $\rho$ is a holomorphic $\C^*$-action on
$\overline{\mathcal X}$ such that
$\overline\pi(\rho_\lambda(x))=\lambda\cdot\overline\pi(x)$,
where $\lambda\cdot[z_0:z_1]=[\lambda z_0:z_1]$ is the standard
action on $\mathbb P^1$.
\item $\Phi$ is a specified $\C^*$-equivariant biholomorphism
over $\mathbb P^1\setminus\{0\}$,
\[
 \Phi:M\times(\mathbb P^1\setminus\{0\})
 \longrightarrow
 \overline{\mathcal X}\setminus\overline\pi^{-1}(0),
\]
where the action on the source is trivial on $M$ and standard
on the base.
\item $\overline{\mathcal A}\in
H^{1,1}_{\mathrm{BC}}(\overline{\mathcal X},\R)$ satisfies
\[
 \rho_\lambda^*\overline{\mathcal A}=\overline{\mathcal A}
 \quad(\lambda\in\C^*),\qquad
 \Phi^*\overline{\mathcal A}=\operatorname{pr}_M^*\Omega.
\]
\end{enumerate}
We write $\mathcal X\coloneqq\overline\pi^{-1}(\C)$,
$\pi\coloneqq\overline\pi|_{\mathcal X}$ and
$\mathcal A\coloneqq\overline{\mathcal A}|_{\mathcal X}$,
where $\C=\mathbb P^1\setminus\{\infty\}$, and abbreviate the
configuration by $(\mathcal X,\mathcal A)$. Its central fibre is
$\mathcal X_0\coloneqq\pi^{-1}(0)$, with its induced complex-space
structure. The configuration is \emph{relatively K\"ahler} if
there is a K\"ahler class $\beta$ on $\mathbb P^1$ such that
\[
 \overline{\mathcal A}+\overline\pi^*\beta
 \quad\text{is K\"ahler on }\overline{\mathcal X}.
\]
\end{definition}

This is the compactified formulation of
\cite{SDPoly}. See also
\cite[Definition~2.10]{DRKahler}. Normality of the total space
does not require the central fibre to be normal, reduced or
irreducible. The stability notions below test normal,
relatively K\"ahler configurations in this sense.

For a resolution $p:\mathcal Y\to\overline{\mathcal X}$, set
$\mathcal B\coloneqq p^*\overline{\mathcal A}$. In complex dimension
$n\coloneqq\dim_{\C}M$, put
\[
 \mu_\Omega\coloneqq
 \frac{c_1(M)\cdot\Omega^{n-1}}{\Omega^n}.
\]
With our normalization $\DF=(a_1/a_0)b_0-b_1$, the
general-dimensional intersection formula is
\[
 \DF(\mathcal X,\mathcal A)
 =\frac{1}{2n!}\left(
       \frac{n}{n+1}\mu_\Omega\,\mathcal B^{n+1}
       +K_{\mathcal Y/\mathbb P^1}\cdot\mathcal B^n\right).
\]
Both intersection products have top degree on the
$(n+1)$-dimensional space $\mathcal Y$. Here
$K_{\mathcal Y/\mathbb P^1}\coloneqq
K_{\mathcal Y}-(\overline\pi\circ p)^*K_{\mathbb P^1}$
is the relative canonical divisor of the resolved total space.
To verify the normalization, consider a polarized test configuration
with a relatively ample $\C^*$-linearized line bundle $\mathcal L$,
where $c_1(\mathcal L)=\mathcal A$ and $c_1(L)=\Omega$.
Take $\overline{\mathcal A}=c_1(\overline{\mathcal L})$ for
its canonical compactification $\overline{\mathcal L}$.
For $k\gg1$, let $h(k)$ be the dimension of
$H^0(\mathcal X_0,\mathcal L_0^k)$ and let $w(k)$ be the total
weight of its induced $\C^*$-action on sections. Write
\[
 h(k)=a_0k^n+a_1k^{n-1}+O(k^{n-2}),\qquad
 w(k)=b_0k^{n+1}+b_1k^n+O(k^{n-1}).
\]
Riemann--Roch \cite[Example~1.1.27]{LazarsfeldI} and its equivariant form give
\[
 \begin{aligned}
 a_0&=\frac{\Omega^n}{n!},&
 a_1&=\frac{c_1(M)\cdot\Omega^{n-1}}{2(n-1)!},\\
 b_0&=\frac{\mathcal B^{n+1}}{(n+1)!},&
 b_1&=-\frac{K_{\mathcal Y/\mathbb P^1}\cdot\mathcal B^n}{2n!}.
 \end{aligned}
\]
Substitution yields the displayed formula, which also defines
the invariant for real cohomological classes. Compare
\cite[Definitions~2.13--2.14]{DRKahler}.
For surfaces, $n=2$, so $1/(2n!)=1/4$ and $n/(n+1)=2/3$.
Thus the formula used throughout this paper is
\begin{equation}\label{rd:eq:cohom-df}
 \DF(\mathcal X,\mathcal A)
 \coloneqq\frac14\left(
       \frac23\mu_\Omega\,\mathcal B^3
       +K_{\mathcal Y/\mathbb P^1}\cdot\mathcal B^2\right).
\end{equation}
This is independent of the resolution and of adding a base class.
See \cite[Definitions~2.13--2.14 and Lemmas~2.15, 2.24]{DRKahler}.
The overall positive normalization factor changes neither signs
nor the stability notion. 

Following \cite{SDPoly},
\emph{cohomological K-polystability} means that every normal
relatively K\"ahler cohomological test configuration has
$\DF\geq0$, with equality only for products over $\C$.
Here a product is equivariantly isomorphic to $M\times\C$,
with a diagonal action induced by a one-parameter subgroup
of $\Aut(M,\Omega)$ and the standard
action on $\C$. The action on $M$ need not be trivial, and the
canonical compactification need not be a product over
$\mathbb P^1$. Dropping the equality clause defines
cohomological K-semistability.

\section{Auxiliary classes and embeddings}
\label{ru:sec:marking-construction}

An irrational K\"ahler class has no positive integer multiple that
is the first Chern class of a line bundle. To obtain projective
embeddings, we therefore use nearby rational K\"ahler classes,
called \emph{auxiliary classes}, whose sufficiently divisible
multiples define ample line bundles.

To use corresponding auxiliary classes on the source and its
geometric limit, we construct a \emph{cohomological marking}.
This is a linear correspondence between selected cohomology
classes that recovers the limiting metric class and preserves
their intersection numbers. We first treat a sequence of cscK
metrics in varying classes. We then treat all limits of a fixed
Calabi flow, obtaining auxiliary polarizations and embedding
powers that work uniformly over those limits.

\subsection{Auxiliary classes on varying cscK limits}
\label{vc:sec:compactness}
We first obtain an orbifold limit by geometric compactness.
The cohomological marking and the auxiliary projective embeddings
are constructed afterwards on this same limit.

Let $M$ be a fixed smooth rational projective surface, put
$c\coloneqq c_1(M)$, and use the controlled cone $\mathcal U$
and the function $\delta$ defined in \eqref{vc:eq:cone}.
The metric class is $[\omega]=2\pi\Omega$.

\begin{theorem}[Static varying-class compactness]
\label{vc:thm:compactness}
Suppose that $\Omega_j$ belong to a compact subset of $\mathcal U$,
$\Omega_j\to\Omega_0$, and $\omega_j\in2\pi\Omega_j$ are cscK.
Then after passing to a subsequence they converge to a connected compact
cscK K\"ahler orbifold $(Z,\omega_\infty)$ with isolated quotient
singularities. The convergence is smooth away from a finite set $\Sigma$,
which may contain smooth points of $Z$, as well as its orbifold points.
The total volume is preserved, and the volumes of the omitted
neighbourhoods tend uniformly to zero. In particular, with
$\Omega_\infty\coloneqq[\omega_\infty]/(2\pi)$,
\[
 \Omega_\infty^2=\Omega_0^2,\qquad
 c_1(Z)\cdot\Omega_\infty=c\cdot\Omega_0.
\]
\end{theorem}

\begin{proof}
The scalar curvature in our convention, and the volume, are
\[
 S_j=\frac{2c\cdot\Omega_j}{\Omega_j^2}>0,
 \qquad v_j=2\pi^2\Omega_j^2.
\]
Write $s_j=2S_j$ for the Riemannian scalar curvature. For a unit-volume
Yamabe minimizer $q_j$ in the conformal class of the metric $g_j$,
the Chern--Gauss--Bonnet and signature identities give
\[
 4\pi^2c^2
 =\int_M\left(2|W^+_{q_j}|^2+\frac{s_{q_j}^2}{24}
               -\frac{|\operatorname{Ric}^0_{q_j}|^2}{2}\right)dV_{q_j}
 \leq \frac1{12}\int_Ms_j^2dV_j+\frac{Y(g_j)^2}{24}.
\]
We used conformal invariance of the Weyl integral and
$|W^+_{g_j}|^2=s_j^2/24$. Thus
\begin{equation}\label{vc:eq:yamabe}
 Y(g_j)^2\geq96\pi^2\delta(\Omega_j).
\end{equation}
Since $s_j$ is a positive constant, the conformal Laplacian is
positive and $Y(g_j)>0$. The Yamabe inequality therefore supplies a
uniform inequality
\[
 \|u\|_{L^4(g_j)}^2\leq C
       \bigl(\|\nabla u\|_{L^2(g_j)}^2
                     +v_j^{-1/2}\|u\|_{L^2(g_j)}^2\bigr).
\]
The volumes are bounded above and away from zero. The same curvature
identities give a uniform bound for $\int_M|\mathrm{Rm}_{g_j}|^2dV_j$.

Apply Tian--Viaclovsky's static compactness theorem after the constant
rescaling $g_j\mapsto v_j^{-1/2}g_j$. See
\cite[Theorem~1.3, Section~3.1, and the conclusion of Section~7
before Section~7.1]{TV}.
See also Anderson's general orbifold compactness theorem
\cite[Theorem~1.1]{AndersonOrbifold}.
For the application of Tian--Viaclovsky's theorem, we use the
Sobolev convention in their equation~(3.11).
Its weaker lower-order term does not by itself imply irreducibility
from $b_1(M)=0$. Irreducibility here follows from their K\"ahler
one-endedness argument. Their definition of convergence, equation~(1.7),
gives volume preservation and vanishing omitted volumes. Parallel
complex structures pass to the smooth regular limit and extend in
the orbifold charts. Constant scalar curvature extends as well.
The scalar-integral identity then gives the asserted canonical
intersection number.
\end{proof}

The classes $\Omega_0$ and $\Omega_\infty$ in
Theorem~\ref{vc:thm:compactness} belong to the cohomology of different
spaces. The numerical equalities in that theorem do not by
themselves specify a correspondence between their cohomology
classes. We next construct such a correspondence on the classes
needed to express $\Omega_0$. Here a \emph{rational subspace} means
a real linear subspace spanned by rational cohomology classes. Its
elements need not themselves be rational.

\begin{lemma}[Cohomological marking of the limiting class]
\label{vc:lem:marking}
In the setting of Theorem~\ref{vc:thm:compactness}, let $W$ be the
smallest rational linear subspace of $H^2(M,\mathbb R)$ containing
$\Omega_0$. Choose integral classes $\beta_1,\ldots,\beta_d$ forming
a rational basis of $W$, and write
\[
 \Omega_0=\sum_{i=1}^d a_i\beta_i.
\]
The coefficients $a_i$ are linearly independent over $\mathbb Q$.
Then there are holomorphic orbifold line bundles $\mathcal B_{i,\infty}$
on $Z$, with rational classes
$\beta_{i,\infty}\coloneqq c_1(\mathcal B_{i,\infty})$, such that
\begin{align}
 \Omega_\infty&=\sum_i a_i\beta_{i,\infty},\label{vc:eq:marking}\\
 \beta_{i,\infty}\cdot\beta_{\ell,\infty}
      &=\beta_i\cdot\beta_\ell,\qquad
 c_1(Z)\cdot\beta_{i,\infty}=c\cdot\beta_i.
 \label{vc:eq:intersection}
\end{align}
The connections and Hermitian structures on these bundles are
limits of those on fixed holomorphic line bundles on $M$.
\end{lemma}

\begin{proof}
Set $\vartheta_j\coloneqq\omega_j/(2\pi)$ and let $\eta_{i,j}$
be the real harmonic representative of $\beta_i$ for $g_j$.
Since $M$ is rational, all its real second cohomology is of type
$(1,1)$. The harmonic representatives are therefore of type $(1,1)$.
The K\"ahler identities
and the Hodge star formula on a surface give
\begin{equation}\label{vc:eq:harmonic}
 \Lambda_{\vartheta_j}\eta_{i,j}
   =\frac{2\beta_i\cdot\Omega_j}{\Omega_j^2},\qquad
 \|\eta_{i,j}\|_{L^2(g_j)}^2
   =\frac{2(\beta_i\cdot\Omega_j)^2}{\Omega_j^2}-\beta_i^2.
\end{equation}
The middle-degree norm is unchanged by constant rescaling of the
metric. In particular these norms are uniformly bounded. Interior
elliptic estimates on the smoothly converging regions give, after
one extraction, smooth limits $\eta_{i,\infty}$ on $Z\setminus\Sigma$.

These forms extend smoothly as orbifold forms at every point of
$\Sigma$. Indeed, pull a limiting form back to a punctured orbifold
four-ball. It is $L^2$, closed, and coclosed for a smooth metric on the
whole ball. Choose cutoffs $\chi_\varepsilon$ vanishing near the origin
and equal to one outside the ball of radius $2\varepsilon$.
Then $\|d\chi_\varepsilon\|_{L^2}=O(\varepsilon)$.
In the weak equations for $d$ and $d^*$, every cutoff error is bounded
by a constant times
$\|\eta\|_{L^2(B_{2\varepsilon})}
 \|d\chi_\varepsilon\|_{L^2}$, and tends to zero.
Thus both equations hold distributionally on the whole ball.
Elliptic regularity for $(d+d^*)^2$ gives a smooth extension, which
is invariant under the chart group. The same argument covers smooth
concentration points.

There is a fixed holomorphic line bundle $\mathcal B_i$ on $M$ with
$c_1(\mathcal B_i)=\beta_i$, by
\cite[Section~2.2, p.~798]{VoisinSurvey}. The $\partial\bar\partial$-lemma
\cite[Chapter~VI, Lemma~8.6]{DemaillyCADG}
gives Hermitian metrics with real Chern curvature
$2\pi\eta_{i,j}$.
On regular compact sets choose unitary local radial gauges for their
connections. Smooth curvature bounds give bounds for these local
connection forms. Transition functions are $U(1)$-valued. After
extraction their derivatives and values converge on overlaps.
This constructs a limiting Hermitian holomorphic line bundle on
$Z\setminus\Sigma$, with curvature $2\pi\eta_{i,\infty}$.
To extend it, lift to a punctured four-ball. The lifted line bundle
is topologically trivial because $H^2(B^4\setminus\{0\},\mathbb Z)=0$.
Choose an invariant smooth primitive for its extended curvature.
Its difference from the connection is closed and is exact since
$H^1(B^4\setminus\{0\},\mathbb R)=0$. A unitary gauge therefore
extends the connection. The deck transformations preserve this
connection and act on the fibre by constant unitary characters,
so they extend across the origin. This produces
$\mathcal B_{i,\infty}$, with
$[\eta_{i,\infty}]=\beta_{i,\infty}$.

The K\"ahler classes $\Omega_j=[\omega_j]/(2\pi)$ on $M$
need not lie in $W$, since $W$ was defined using their limiting
class $\Omega_0$. Nevertheless,
$\vartheta_j-\sum_i a_i\eta_{i,j}$ is the harmonic representative of
$\epsilon_j\coloneqq\Omega_j-\Omega_0$. Its squared norm is
\[
 \frac{2(\epsilon_j\cdot\Omega_j)^2}{\Omega_j^2}
                      -\epsilon_j^2\longrightarrow0.
\]
Smooth regular convergence gives
$\omega_\infty/(2\pi)=\sum_i a_i\eta_{i,\infty}$. The extensions
give \eqref{vc:eq:marking} on all of $Z$.
Passing the constant traces in \eqref{vc:eq:harmonic} to the limit,
and using volume preservation, gives
$\beta_{i,\infty}\cdot\Omega_\infty=\beta_i\cdot\Omega_0$.
Consequently
\[
 \sum_\ell a_\ell
 \bigl(\beta_{i,\infty}\cdot\beta_{\ell,\infty}
                       -\beta_i\cdot\beta_\ell\bigr)=0.
\]
The numbers in parentheses are rational: a finite tensor power of
each orbifold line bundle descends to a line bundle on the coarse
space. The minimality of $W$ implies that the $a_\ell$ are
$\mathbb Q$-linearly independent, proving the first identities in
\eqref{vc:eq:intersection}. The preserved canonical degree gives
\[
 \sum_i a_i
 \bigl(c_1(Z)\cdot\beta_{i,\infty}-c\cdot\beta_i\bigr)=0,
\]
and proves the remaining identities by the same rationality argument.
\end{proof}

With the notation of Lemma~\ref{vc:lem:marking}, we define
\[
 \mathfrak m_Z:W\longrightarrow H^2_{\mathrm{orb}}(Z,\mathbb R),
 \qquad
 \mathfrak m_Z\!\left(\sum_i x_i\beta_i\right)
       \coloneqq\sum_i x_i\beta_{i,\infty}.
\]
Here $H^2_{\mathrm{orb}}(Z,\mathbb R)$ denotes the de Rham
cohomology of smooth orbifold differential forms. We call
$\mathfrak m_Z$ the \emph{cohomological marking} associated with
the extracted subsequence, and $\mathfrak m_Z(B)$ the
\emph{marked class} corresponding to $B\in W$.
By \eqref{vc:eq:marking}--\eqref{vc:eq:intersection}, it sends
$\Omega_0$ to $\Omega_\infty$ and preserves intersection pairings
and degrees against the first Chern class. The marking may depend
on the extraction and is defined only on $W$.

The basis classes $\beta_i$ need not be K\"ahler. We therefore
choose nearby rational classes $B\in W$ and their marked images
$B_\infty\coloneqq\mathfrak m_Z(B)$ as auxiliary classes.
The next proposition establishes their positivity and the
resulting complete projective embeddings.

\begin{proposition}[Auxiliary complete linear series on the same limit]
\label{vc:prop:polarized-limit}
For every rational tuple $(b_i)$ sufficiently close to $(a_i)$, put
\[
 B\coloneqq\sum_i b_i\beta_i,\qquad
 B_\infty\coloneqq\sum_i b_i\beta_{i,\infty}.
\]
Then $B$ and $B_\infty$ are K\"ahler classes. After taking a
sufficiently divisible integer $m=m(B,Z)>0$, the complete linear
series of $mB$ on $M$ and $mB_\infty$ on the coarse space of $Z$
give embeddings with the same Hilbert polynomial, and the latter
embedded surface is a scheme-theoretic Hilbert limit of the former.
Moreover
$h^1(Z,\mathcal O_Z)=h^2(Z,\mathcal O_Z)=0$.
\end{proposition}

\begin{proof}
The smooth orbifold form
$\theta_\infty\coloneqq\sum_i b_i\eta_{i,\infty}$ tends to
$\omega_\infty/(2\pi)$ as $b_i\to a_i$, so it is positive for
nearby tuples. The class $B$ is likewise K\"ahler on $M$, since
$\Omega_0$ is an interior point of its K\"ahler cone. Choose $q>0$
such that all $qb_i$ are integers, and set
\[
 L_B\coloneqq\bigotimes_i\mathcal B_i^{qb_i},\qquad
 L_{B,\infty}\coloneqq
              \bigotimes_i\mathcal B_{i,\infty}^{qb_i}.
\]
Use the product Hermitian metrics provided by
Lemma~\ref{vc:lem:marking}. The source curvature
$2\pi q\sum_i b_i\eta_{i,j}$ need not be positive, but its trace
with respect to $\omega_j$ is the constant
\[
 \frac{2qB\cdot\Omega_j}{\Omega_j^2}.
\]
For each fixed $k$, the section Bochner identity
\cite[Chapter~VII, Section~1]{DemaillyCADG}, followed by Moser
iteration using the uniform Sobolev bound, gives
\[
 \sup_M|s|^2\leq C_{k,B}\int_M|s|^2dV_j,
 \qquad s\in H^0(M,L_B^k).
\]
Take one diagonal subsequence over the positive integers $k$.
Orthonormal bases then have smooth limits on regular compact
sets. The uniform upper bound and vanishing omitted volumes preserve
all inner products. On a uniformizing ball, the limiting holomorphic sections extend
across the origin by Hartogs' theorem
\cite[Chapter~I, Theorem~3.28]{DemaillyCADG}. Thus, for
every fixed $k$, there is a dimension-preserving injection of source
sections into $H^0(Z,L_{B,\infty}^k)$.

The positive limiting curvature and Baily's global embedding theorem
\cite[Section~7, Global Imbedding Theorem]{Baily} imply that a
sufficiently divisible power of $L_{B,\infty}$ is an ample Cartier
bundle on the coarse space. For all sufficiently large divisible $k$, Serre vanishing
\cite[Theorem~1.2.6(ii)]{LazarsfeldI} and Riemann--Roch on a
resolution give
\[
 h^0(Z,L_{B,\infty}^k)
 =\frac{q^2k^2}{2}B_\infty^2
       +\frac{qk}{2}c_1(Z)\cdot B_\infty+\chi(\mathcal O_Z).
\]
Indeed, let $p:Y\to Z$ be a resolution. Quotient singularities
are rational \cite[Proposition~5.15]{KM}, so
$p_*\mathcal O_Y=\mathcal O_Z$ and
$R^ip_*\mathcal O_Y=0$ for $i>0$, as in
\cite[Definition~5.8]{KM}.
For a Cartier line bundle $A$ on $Z$, the projection formula
and the Leray spectral sequence therefore give
\[
 H^i(Z,A)\cong H^i(Y,p^*A).
\]
Moreover, $p^*A$ has zero intersection with every exceptional
curve. The surface Riemann--Roch formula on $Y$ consequently
gives the displayed formula on $Z$. Alternatively, the same formula follows directly from Kawasaki's
orbifold Riemann--Roch theorem
\cite[formula~(IV), p.~153]{Kawasaki}.
For the sufficiently divisible powers considered here, the local
stabilizers act trivially on the bundle fibres. Since the
singularities are isolated, the singular contributions are
independent of $k$ and are included in the constant term
$\chi(\mathcal O_Z)$.

The quadratic and
linear coefficients coincide with those on $M$, by
\eqref{vc:eq:intersection}. The section injection and
$\chi(\mathcal O_M)=1$ imply $\chi(\mathcal O_Z)\geq1$.

For nearby $B$, its canonical degree is negative:
$K_Z\cdot B_\infty=-c\cdot B<0$. A nonzero section of the
dualizing sheaf $\mathcal O_Z(K_Z)$ would give an effective canonical Weil divisor of
negative degree with respect to an ample class, which is impossible.
Since rational singularities are Cohen--Macaulay
\cite[Theorem~5.10]{KM}, Serre duality gives
\[
 H^2(Z,\mathcal O_Z)\cong H^0(Z,\mathcal O_Z(K_Z))^*=0.
\]
See \cite[Theorem~5.71 and Proposition~5.75]{KM}. Connectedness and
normality imply
$\chi(\mathcal O_Z)=1-h^1(\mathcal O_Z)\leq1$.
Equality follows, so both higher cohomology groups vanish and the
two Hilbert polynomials coincide.

Take one sufficiently large common power that is very ample
and has vanishing higher cohomology on both spaces, using
\cite[Theorem~1.2.6(ii), (iv)]{LazarsfeldI}. The limiting orthonormal sections then form the complete
section space on $Z$.
We use the projective Hilbert scheme with the prescribed
Hilbert polynomial and its universal flat family, as constructed in
\cite[Sections~5.1 and~5.5]{NitsureHilbert}.
If $W'$ is a Hilbert-scheme limit of the embedded $M$, limiting
homogeneous equations vanish on the regular embedded limit and
hence on reduced $Z$. Thus $Z\subset W'$ as schemes. Let $\mathcal K$ be the kernel of
$\mathcal O_{W'}\to\mathcal O_Z$. Equality of Hilbert polynomials
gives $\chi(\mathcal K(k))=0$. For sufficiently large $k$,
Serre vanishing gives $H^0(W',\mathcal K(k))=0$, while
$\mathcal K(k)$ is globally generated
\cite[Theorem~1.2.6 (iii)]{LazarsfeldI}.
Hence $\mathcal K=0$ and $Z=W'$.
The required exponent is allowed to depend on this limit and on $B$.
\end{proof}

\subsection{Marked classes and uniform auxiliary ampleness}\label{ru:sec:marking}
We now apply the cohomological marking of
Subsection~\ref{vc:sec:compactness} to the geometric limits of a
fixed Calabi flow. The source class is fixed, but the marking must
be constructed for every extracted limit. We seek a common
neighborhood of auxiliary classes and common embedding powers that
work for all these marked limits. Bounds on curve intersections
will give the required uniform ampleness.

For the fixed-flow argument put $c\coloneqq c_1(M)$ and
$\vartheta(t)\coloneqq\omega(t)/(2\pi)$. The controlled-cone
conditions make $M$ rational. The long-time theorem supplies the
global flow and a uniform Sobolev constant, so
Theorem~\ref{app:thm:compactness} applies. Set $V_0\coloneqq\Omega^2$ and
$s_0\coloneqq c\cdot\Omega>0$.

Since $M$ is rational, $H^{2,0}(M)=H^{0,2}(M)=0$.
The Lefschetz theorem on $(1,1)$-classes therefore implies that
every integral second cohomology class is the first Chern class
of a holomorphic line bundle.

Let $W$ be the smallest rational real subspace of $H^2(M,\R)$
containing $\Omega$. As in Subsection~\ref{vc:sec:compactness}, we
choose integral classes $\beta_1,\ldots,\beta_d$ forming a
rational basis of $W$, and write
\[
 \Omega=\sum_{i=1}^d a_i\beta_i.
\]
The coefficients $a_i$ are linearly independent over $\Q$.
Let $\eta_i(t)$ be the harmonic representative of $\beta_i$,
and choose a holomorphic line bundle $\mathcal B_i$ with
$c_1(\mathcal B_i)=\beta_i$. We equip these bundles with smooth
families of Hermitian metrics whose real Chern curvatures are
$2\pi\eta_i(t)$. Multiplication by positive factors depending
only on time does not change their Chern connections or the
resulting markings.

For each extracted geometric limit $Z$, the construction of
Lemma~\ref{vc:lem:marking} gives limiting orbifold line bundles
$\mathcal B_{i,Z}$. We use the marking notation introduced above:
\[
 \beta_{i,Z}\coloneqq c_1(\mathcal B_{i,Z}),\qquad
 \mathfrak m_Z(\beta_i)=\beta_{i,Z},\qquad
 B_Z\coloneqq\mathfrak m_Z(B)\quad(B\in W).
\]
Let $\mathfrak L$ denote the collection of geometric limits
equipped with all markings obtained by these extractions.
The next lemma establishes the marking identities for flow
slices and the additional uniform bounds needed below.

\begin{lemma}[Uniform marking for the fixed flow]
\label{ru:lem:marking}
Every geometric limit admits a marking as described above.
For every marked limit $Z\in\mathfrak L$, the Chern classes
$\beta_{i,Z}$ have limiting harmonic representatives and satisfy
\begin{equation}\label{ru:eq:marking}
 \Omega_Z=\sum_i a_i\beta_{i,Z},\qquad
 \beta_{i,Z}\beta_{j,Z}=\beta_i\beta_j,\qquad
 c_Z\beta_{i,Z}=c\beta_i,
\end{equation}
where $\Omega_Z\coloneqq[\omega_Z]/(2\pi)$ and $c_Z\coloneqq c_1(Z)$.
There is an integer $Q>0$, independent of $Z$ and of the subsequence,
such that every $\mathcal B_{i,Z}^Q$ descends to a Cartier line bundle.
Moreover $c_Z^2$ has a uniform lower bound.
\end{lemma}

\begin{proof}
The proof of Lemma~\ref{vc:lem:marking} applies to these fixed-class
flow slices. We specify the inputs that replace its cscK assumption.
The harmonic identities are
\begin{equation}\label{ru:eq:harmonic}
 \Lambda_{\vartheta(t)}\eta_i(t)
 =\frac{2\beta_i\Omega}{V_0},\qquad
 \|\eta_i(t)\|_2^2
 =\frac{2(\beta_i\Omega)^2}{V_0}-\beta_i^2,
 \qquad \vartheta(t)=\sum_i a_i\eta_i(t).
\end{equation}
Local elliptic estimates give smooth limits on regular compact
sets. The removability and bundle-extension arguments in the
proof of Lemma~\ref{vc:lem:marking} apply here and produce the
limiting orbifold line bundles $\mathcal B_{i,Z}$.

Constant traces and volume preservation imply
$\beta_{i,Z}\Omega_Z=\beta_i\Omega$.
The Calabi energy decreases, so $\int S(t)^2dV_t$ is bounded.
Cauchy--Schwarz and Lemma~\ref{app:lem:tails} therefore preserve
the scalar integral and give $c_Z\Omega_Z=s_0$.
All intersections of orbifold Chern classes are rational.
Substituting $\Omega_Z=\sum a_i\beta_{i,Z}$ in these identities
and using the rational independence of the $a_i$ proves
\eqref{ru:eq:marking}. 

The Sobolev inequality passes to each limit. Applying its local
noncollapsing consequence \eqref{tc:eq:uniform-noncollapse} at an
orbifold point, we bound the order
of its local group uniformly: the density there is
$\Vol(B^4_1)/|\Gamma|$. Thus $|\Gamma|\leq N_0$, with $N_0$
depending only on the fixed Sobolev bound and normalization.
Take $Q=\operatorname{lcm}(1,\ldots,\lfloor N_0\rfloor)$.
Every local fibre character is killed by $Q$, and the descent
criterion in Subsection~\ref{app:sec:geometric-compactness} applies.

Finally the K\"ahler surface Chern identities and the scalar
$L^2$ bound give a uniform bound for $\int_M|\operatorname{Ric}|^2$.
By smooth regular convergence and Fatou's lemma, we get the same bound
on each limit. Since $c_Z$ is represented by $\rho_Z/(2\pi)$,
$|\int_Z\rho_Z\wedge\rho_Z|\leq C\int_Z|\operatorname{Ric}|^2$.
This gives the asserted lower bound for $c_Z^2$.
\end{proof}

\begin{lemma}[Common auxiliary ample neighborhood]
\label{ru:lem:ample}
There is a neighborhood $U$ of $\Omega$ in $W$ such that every
rational $B=\sum_i b_i\beta_i\in U$ is ample on $M$ and its
marked class $B_Z\coloneqq\sum_i b_i\beta_{i,Z}$ is ample on
the coarse space of every $Z\in\mathfrak L$.
\end{lemma}

\begin{proof}
\textit{Step 1. Projectivity of each individual limit.}
For a fixed $Z$, the smooth harmonic representatives in
\eqref{ru:eq:marking} sum to $\omega_Z/(2\pi)$.
A sufficiently close rational combination is still positive on
this compact orbifold. A divisible integral multiple is the
curvature class of an orbifold line bundle. Baily's global
embedding theorem \cite[Section~7]{Baily} makes $Z$ projective.
This choice is allowed to depend on $Z$. 

\textit{Step 2. Adjunction on a quotient surface.}
Every integral curve $C\subset Z$ satisfies
\begin{equation}\label{ru:eq:adjunction}
 C^2-c_ZC\geq-2.
\end{equation}
Indeed, let $p:Y\to Z$ be the minimal resolution. Its exceptional curves $E_j$ are smooth rational curves by
\cite[Theorem~4.7 and Proposition~4.18]{KM}.
Negative definiteness and minimality of the resolution imply
$E_j^2\leq-2$ \cite[Section~I, p.~6]{MumfordSurface}.
The numerical pullback has the form
$p^*C=\widetilde C+D$, where $D$ is effective and exceptional,
as in \cite[Section~II(b), pp.~17--18]{MumfordSurface}.
It is characterized by $(p^*C)\cdot E_j=0$ for every exceptional
curve $E_j$. The exceptional intersection matrix is negative
definite with nonnegative
off-diagonal entries. Its inverse has nonpositive entries, so
the equations defining $D$ give its effectiveness. Orthogonality to the
exceptional curves gives
\[
 C^2=\widetilde C^{\,2}-D^2\geq\widetilde C^{\,2}.
\]
Also $K_YE_j=-2-E_j^2\geq0$, and the discrepancy divisor
has zero intersection with $p^*C$. Consequently
\[
 K_ZC=K_Yp^*C\geq K_Y\widetilde C.
\]
The arithmetic-genus form of adjunction on the smooth surface $Y$
\cite[Chapter~V, Theorem~1.6 and Exercise~1.3(a)]{HartshorneAG}
gives
\[
 \widetilde C^{\,2}+K_Y\cdot\widetilde C
 =2p_a(\widetilde C)-2\geq-2,
\]
and proves \eqref{ru:eq:adjunction}.

\textit{Step 3. Bound all marked curve degrees by their area.}
On the projective resolution $p:Y\to Z$ used above, the numerical
pullbacks of $\Omega_Z$, $c_Z$ and the $\beta_{i,Z}$ belong to
$N^1(Y)_{\mathbb R}$, the real space of divisor classes modulo
numerical equivalence. For $\Omega_Z$, this follows from
$\Omega_Z=\sum_i a_i\beta_{i,Z}$.
The Hodge index theorem gives the intersection form signature
$(1,\dim N^1(Y)_{\mathbb R}-1)$ on this space
\cite[Chapter~V, Theorem~1.9 and Remark~1.9.1]{HartshorneAG}.
Since $(p^*\Omega_Z)^2=V_0>0$, its orthogonal complement is
negative definite, even though $p^*\Omega_Z$ need not be ample.
We use the norm $\|\xi\|_-^2\coloneqq-\xi^2$ on this complement.
Write $A\coloneqq\Omega_ZC>0$ and
$u^2\coloneqq A^2/V_0-C^2\geq0$.
Thus $u$ is the norm of the perpendicular component of $p^*C$.
The squared norm of the pullback of
$c_Z-(s_0/V_0)\Omega_Z$ is $s_0^2/V_0-c_Z^2$, which is
uniformly bounded by Lemma~\ref{ru:lem:marking}.
Let $C_c$ be a uniform bound for this norm. Cauchy--Schwarz on
the negative definite orthogonal complement and
\eqref{ru:eq:adjunction} give
\[
 c_ZC\geq s_0A/V_0-C_cu,
 \qquad
 u^2\leq A^2/V_0-s_0A/V_0+C_cu+2
       \leq A^2/V_0+C_cu+2.
\]
It follows that $u\leq C_1(A+1)$.
The corresponding perpendicular norm of each $\beta_{i,Z}$
is fixed by \eqref{ru:eq:marking}. Therefore, with the Euclidean
norm on the degree vector $v(C)\coloneqq(\beta_{i,Z}C)_{i=1}^d$
and $a\coloneqq(a_1,\ldots,a_d)$,
\begin{equation}\label{ru:eq:lattice}
 |v(C)|\leq C_2(A+1),\qquad
 v(C)\in Q^{-1}\Z^d,\qquad A=a\cdot v(C).
\end{equation}
The lattice assertion holds because $Q\beta_{i,Z}$ is Cartier.
For $0<A\leq1$ the vectors in \eqref{ru:eq:lattice} belong to
one finite set. Since $a$ is fixed, the positive values of
$a\cdot v$ on this finite set have a positive minimum.
For $A\geq1$ the first inequality gives $|v(C)|\leq2C_2A$.
Together these observations give one constant $C_3$ such that
\begin{equation}\label{ru:eq:linear-degree}
 |v(C)|\leq C_3\Omega_ZC
 \quad\text{for every curve on every limit.}
\end{equation}
If the finite set with $0<A\leq1$ is empty, only the second
case is needed.

\textit{Step 4. Choose a single neighborhood.}
Write $b\coloneqq(b_1,\ldots,b_d)$ for the coefficient vector
of $B=\sum_i b_i\beta_i$. If $|b-a|<(2C_3)^{-1}$, then
\[
 B_ZC=b\cdot v(C)\geq\tfrac12\Omega_ZC>0.
\]
Shrink this neighborhood so that $B$ is K\"ahler on $M$,
$B^2>0$, and $cB>0$. The marking gives $B_Z^2=B^2$.
For rational $B$, a divisible multiple is Cartier on $Z$.
Applying the Nakai--Moishezon criterion to a Cartier multiple of
$B_Z$ shows that $B_Z$ is ample by
\cite[Theorem~1.2.23 and Definition~1.3.6]{LazarsfeldI}.
This proves the lemma. This proof uses the
fixed coefficients of $\Omega$ and does not assert a uniform
neighborhood as $\Omega$ varies.
\end{proof}

\begin{lemma}[Common complete auxiliary embeddings]
\label{ru:lem:embeddings}
There are fixed line bundles $L_1,\ldots,L_d$ on $M$ whose first
Chern classes span $W$, and fixed positive integers $r_\nu$,
such that every limit has ample Cartier bundles $L_{\nu,Z}$
with $c_1(L_{\nu,Z})$ the marked image of $c_1(L_\nu)$.
The bundles $A_\nu\coloneqq L_\nu^{r_\nu}$ and
$A_{\nu,Z}\coloneqq L_{\nu,Z}^{r_\nu}$ are very ample, have
vanishing higher cohomology, and have matching dimensions of
sections. Along any geometrically convergent sequence, after further
extraction and compatible bundle identifications, complete orthonormal
bases of the fixed spaces $V_\nu\coloneqq H^0(M,A_\nu)$ converge
to complete bases on $Z$. The resulting embeddings in
$\prod_\nu\mathbb P(V_\nu^*)$ have a fixed Hilbert polynomial
and converge in the Hilbert scheme. Moreover
$H^1(Z,\mathcal O_Z)=H^2(Z,\mathcal O_Z)=0$.
The auxiliary classes can be chosen so that, with fixed
$\alpha_\nu>0$,
\begin{equation}\label{rd:eq:positive-basis}
 \Omega=\sum_{\nu=1}^d\alpha_\nu c_1(A_\nu),\qquad
 \Omega_Z=\sum_{\nu=1}^d\alpha_\nu c_1(A_{\nu,Z})
 \quad\text{for every }Z\in\mathfrak L.
\end{equation}
\end{lemma}

\begin{proof}
Choose rational classes $B_1,\ldots,B_d$ in the common neighborhood
$U$ of Lemma~\ref{ru:lem:ample} whose positive cone contains
$\Omega$ in its interior. Here is an explicit justification.
Choose a rational linear functional $\ell$ on $W$ with
$\ell(\Omega)>0$. In the rational affine hyperplane $\ell=1$,
choose a small simplex with $d$ rational vertices $v_\nu$
containing $\Omega/\ell(\Omega)$ in its interior. Choose a
positive rational $q$ sufficiently close to $\ell(\Omega)$,
so that $B_\nu\coloneqq qv_\nu$ all belong to $U$.
The $B_\nu$ are linearly independent. If $\gamma_\nu>0$ are
the barycentric coordinates of $\Omega/\ell(\Omega)$, then
$\Omega=\sum_\nu(\ell(\Omega)\gamma_\nu/q)B_\nu$.
For $d=1$, the same assertion follows by taking any positive
rational class on the ray sufficiently close to $\Omega$.
We clear their
denominators and multiply by $Q$ from Lemma~\ref{ru:lem:marking}.
This gives the fixed source bundles and Cartier limiting bundles.
For each fixed $k\geq1$, let $h_{\nu,k}(t)$ be the induced harmonic
Chern metric on $L_\nu^k$. Its curvature has constant trace against
$\omega(t)$, independent of time. The section Bochner inequality
and the uniform Sobolev bound give
\begin{equation}\label{ru:eq:moser}
 \sup_M|s|_{h_{\nu,k}(t)}^2
 \leq C_{\nu,k}\int_M|s|_{h_{\nu,k}(t)}^2dV_t,
 \qquad s\in H^0(M,L_\nu^k).
\end{equation}
The constant may depend on $\nu$ and $k$, but not on time or the
extracted limit.
The extraction, Hartogs extension and volume-tail argument of
Proposition~\ref{vc:prop:polarized-limit} give injections of
source section spaces into limiting section spaces, for every
fixed power after diagonal extraction.

For a sufficiently large power on one given limit, Riemann--Roch
and Serre vanishing apply. The quadratic and linear coefficients
agree by \eqref{ru:eq:marking}, so these injections imply
$\chi(\mathcal O_Z)\geq\chi(\mathcal O_M)=1$.
The negative canonical degree against an ample auxiliary class
gives $H^2(Z,\mathcal O_Z)=0$. By connectedness,
$\chi(\mathcal O_Z)\leq1$. Thus both higher cohomology groups
vanish and each auxiliary Hilbert polynomial is fixed.

The quotient singularities are rational. Matsusaka's theorem for normal
projective varieties with rational singularities and a fixed
Hilbert polynomial gives a common very ample power with vanishing
higher cohomology
\cite[Theorem~2.4, pp.~206--207, and p.~176]{Matsusaka}.
For each of the finitely many $L_\nu$, we apply the theorem to the
family of all limiting polarized spaces together with the original
polarized surface $(M,L_\nu)$. The resulting integers $r_\nu$
are chosen before any subsequence and work on the source and every limit.
Dimension equality makes each limiting basis complete.
Clearing denominators and taking the powers $r_\nu$, we divide the
coefficients in the positive expression for $\Omega$ by positive
integers. They therefore give \eqref{rd:eq:positive-basis}.
The second equality follows from the marking.

Compose the product embedding with the Segre embedding
\cite[Chapter~I, Exercise~2.14]{HartshorneAG}.
The hyperplane bundle of the resulting projective embedding
restricts to $\bigotimes_\nu A_\nu$ on $M$ and to
$\bigotimes_\nu A_{\nu,Z}$ on $Z$. Its restriction to $M$
therefore has first Chern class
\[
 \sum_\nu c_1(A_\nu)=\sum_\nu r_\nu c_1(L_\nu).
\] The mixed intersections are
fixed by \eqref{ru:eq:marking}, so its Hilbert polynomial is also
fixed. Any Hilbert limit contains the reduced embedded geometric
limit, since limiting homogeneous equations vanish on its regular
locus. The equality of Hilbert polynomials makes this inclusion an
equality of schemes, as in Proposition~\ref{vc:prop:polarized-limit}.
We retain the product factors throughout: preserving their line
bundles preserves $\Omega$, whereas preserving only their sum
need not do so.
\end{proof}

\section{Convergence and class-preserving degeneration}
\label{ru:app:proof}\label{ru:sec:unique}

We prove Theorem~\ref{ru:thm:unique} using the auxiliary bundles
constructed above. By convergence of their Gram matrices, we get common
finite-dimensional dynamics. A countability argument and extremal
uniqueness identify the limit of each flow. A comparison of the
auxiliary Gram matrices for two flows then proves independence of
the initial metric. We finally construct a degeneration preserving
the original real class and apply its stability condition.

\subsection{Convergence of Gram matrices}

For each auxiliary line bundle $A_\nu$ constructed in
Lemma~\ref{ru:lem:embeddings}, choose a Hermitian metric
$h_\nu(t)$ whose real Chern curvature is $2\pi\eta_\nu(t)$,
where $\eta_\nu(t)$ is the $\omega(t)$-harmonic representative
of $c_1(A_\nu)$. We call these the \emph{auxiliary Hermitian
metrics}. They are fibre metrics on the line bundles $A_\nu$,
while $\omega(t)$ is the evolving K\"ahler metric on $M$.
Each $\eta_\nu(t)$ is a fixed integral linear combination
of the forms in \eqref{ru:eq:harmonic}.
The section estimate \eqref{ru:eq:moser} applies with $k=r_\nu$
to $A_\nu=L_\nu^{r_\nu}$ and is unchanged by a scalar rescaling
of $h_\nu(t)$.

We take Hermitian products to be linear in the first entry and
conjugate-linear in the second. These metrics induce $L^2$ inner products on
$V_\nu=H^0(M,A_\nu)$:
\[
 H_\nu(t)(s,s')\coloneqq
       \int_M h_\nu(t)(s,s')\,dV_t.
\]
We prove convergence of the associated Gram matrices using
the harmonic curvature forms $\eta_\nu(t)$.

For an ordered basis $s=(s_1,\ldots,s_{N_\nu})$ of
$V_\nu=H^0(M,A_\nu)$, where $N_\nu\coloneqq\dim V_\nu$,
the \emph{Gram matrix} is the positive Hermitian matrix with entries
\[
 \bigl(\operatorname{Gram}_{H_\nu(t)}(s)\bigr)_{ab}
 \coloneqq H_\nu(t)(s_a,s_b)
 =\int_M h_\nu(t)(s_a,s_b)\,dV_t,
 \qquad 1\leq a,b\leq N_\nu.
\]
We regard an ordered basis as a column of sections. If $s'=Ps$,
then
\[
 \operatorname{Gram}_{H}(s')
 =P\operatorname{Gram}_{H}(s)P^*,
\]
where $P^*$ denotes the conjugate transpose. The same matrix $P$
changes the evaluation coordinates of the complete embedding.
At each time, the prescribed curvature determines $h_\nu(t)$
up to multiplication by a positive factor $c_\nu(t)$ that is
constant over $M$. This factor may depend on both $\nu$ and $t$.
Replacing $h_\nu(t)$ by $c_\nu(t)h_\nu(t)$ multiplies
$H_\nu(t)$ and its Gram matrices by the same factor.

For a positive Hermitian matrix $G$ of size $N$, we define
\[
 \mathcal N(G)\coloneqq(\det G)^{-1/N}G.
\]
Then
\[
 \det\mathcal N(G)=1,\qquad
 \mathcal N(cG)=\mathcal N(G)\quad(c>0).
\]
Thus this normalization removes the common scale of the Gram
matrix while retaining the relative sizes of its eigenvalues.
We choose a fixed source basis in each $V_\nu$ and use the scalar
freedom in $h_\nu(t)$ to impose $\det H_\nu(t)=1$, with the
determinant computed in that basis.

Throughout this subsection, translated time windows are
understood in the sense of Definition~\ref{app:def:time-window}.

In the following lemma, the basis $s_{\nu,j}$ is chosen to be
$H_\nu(t_j)$-orthonormal and is then kept fixed as $\tau$ varies.
Consequently,
$\operatorname{Gram}_{H_\nu(t_j+\tau)}(s_{\nu,j})$
records the changing inner product on these fixed sections.
It equals $I_{N_\nu}$ at $\tau=0$, but need not be the identity
at other times. Convergence on time windows means convergence
uniformly for $|\tau|\leq T$, for every fixed $T>0$.

\begin{lemma}[Auxiliary Gram convergence]
\label{ru:lem:gram}
Let $t_j\to\infty$ and let $s_{\nu,j}$ be any ordered
$H_\nu(t_j)$-orthonormal bases. After passing to a subsequence compatible with the geometric
convergence, there are trace-free Hermitian matrices $C_{\nu,Z}$
such that
\begin{equation}\label{ru:eq:gram}
 \mathcal N\!\left(
  \operatorname{Gram}_{H_\nu(t_j+\tau)}(s_{\nu,j})\right)
       \longrightarrow e^{-2\tau C_{\nu,Z}}
\end{equation}
uniformly on every bounded interval. The projective action of each
$C_{\nu,Z}$ on the corresponding complete embedding is the
same extremal motion $F_\tau$ on $Z$.
\end{lemma}

\begin{proof}
\textit{Step 1. Harmonic curvatures on one time window.}
Use the time-independent comparison maps of
Theorem~\ref{app:thm:compactness}. The harmonic curvature forms
have fixed $L^2$ bounds by \eqref{ru:eq:harmonic}, hence uniform
spatial elliptic bounds on regular compact sets throughout a
bounded time window. Their periods on every fixed regular cycle
are independent of $\tau$, since the source class and comparison
maps are fixed. Removing finitely many points with spherical
space-form links does not change real $H^2$. This follows from
the relative cohomology sequence and
$H^1(S^3/\Gamma,\R)=H^2(S^3/\Gamma,\R)=0$.
Thus every extended limit has the same marked cohomology class.
The unique harmonic representative for $F_\tau^*\omega_Z(0)$ is
$F_\tau^*\eta_{\nu,Z}(0)$, because $F_\tau$ is isotopic to the
identity. Extracting also for any sequence $\tau_j\to\tau$
proves convergence to these forms, uniformly in $\tau$ on compact
intervals in every spatial $C^k$ norm. 

\textit{Step 2. Control the relative Hermitian metrics.}
At time zero of the window, Lemma~\ref{ru:lem:embeddings} gives
bundle identifications and convergence
$s_{\nu,j}\to\sigma_\nu$ to a complete orthonormal basis.
Fix $\nu$ and suppress it from the notation. For any
$\tau_j\to\tau$, rescale $h(t_j+\tau_j)$ temporarily so that
its Gram matrix in $s_j$ has trace one. Denote the rescaled metric
by $\widehat h_j$, and put
\[
 q_j\coloneqq\widehat h_j/h(t_j),\qquad
 B_j\coloneqq\sum_{a=1}^N|s_{j,a}|_{h(t_j)}^2.
\]
The ratio is a function on the same source holomorphic bundle.
We pull it back by the time-independent maps. Equation
\eqref{ru:eq:moser} and trace normalization imply
\begin{equation}\label{ru:eq:trace-control}
 q_jB_j\leq C,\qquad
 \int_M q_jB_j\,dV_{t_j+\tau_j}=1.
\end{equation}
On each regular compact set $B_j$ converges to the Bergman sum
of a complete very ample system, which is strictly positive.
Consequently $q_j$ has a uniform upper bound there.
Choose a fixed regular compact set $K_0$ with omitted volumes
less than $1/(2C)$. Equation~\eqref{ru:eq:trace-control} puts at
least half of the integral on $K_0$. Since $B_j$ is uniformly
bounded there and its volume is bounded, there is a point
$x_j\in K_0$ with $q_j(x_j)\geq c_0>0$.

The complex Hessian of $\log q_j$ is the difference of the two
Chern curvatures, with sign determined by the Chern convention.
Step~1 bounds it in every spatial $C^k$ norm on regular compact
sets. On a connected relatively compact domain containing $K_0$
and a prescribed compact set, choose $A$ with
$u_j\coloneqq A-\log q_j\geq0$. Its Laplacian is uniformly
bounded. The inhomogeneous Harnack inequality, applied along a
finite chain of balls from $K_0$, bounds $u_j$ from above on
the prescribed compact set. The point $x_j$ supplies the bound
for its infimum. Thus $q_j$ is bounded above and below there.
Interior elliptic estimates for $\log q_j$ give smooth positive
subsequential limits. This argument prevents a time-dependent
complex gauge from destroying convergence of the fixed sections.

\textit{Step 3. Identify the limit and its scalar ambiguity.}
The imaginary extremal motion preserves the limiting K\"ahler
metric and every harmonic auxiliary curvature. Since $q(Z)=0$, the action of an isometry torus containing
the extremal Killing field lifts, after passing to a finite
cover, to all the auxiliary bundles
\cite[Lemma~4.2, Step~1]{Polarized}. Averaging a Hermitian metric
with the prescribed invariant curvature gives unitary lifts.
Complexification lifts $F_\tau$ holomorphically. Its pullback
of $h_Z(0)$ has exactly the curvature obtained in Step~1.
The metric produced in Step~2 has that same curvature on the
regular locus. Their logarithmic ratio is pluriharmonic there.
On a punctured uniformizing two-ball, it is the real part of a
holomorphic function because the punctured ball is simply
connected. Hartogs extension extends this function across the
origin. The logarithmic ratio therefore extends to a smooth
orbifold pluriharmonic function on compact $Z$, and is constant.
It follows that the limiting metric is
$c(\tau)\widetilde F_\tau^*h_Z(0)$ for a positive scalar.

The section upper bound and volume tails pass all inner products
to the limit. Indeed each squared section norm is bounded by
the trace-one sum in \eqref{ru:eq:trace-control}. Its integral on
an omitted set is at most $C$ times the omitted volume. This
also controls mixed pairings. The resulting Gram matrix is
positive definite. The pullback group calculation of
\cite[Lemma~4.1]{Polarized} now gives
\[
 G_Z(\tau)=c(\tau)e^{-2\tau\widetilde C_Z},\qquad
 \widetilde C_Z^*=\widetilde C_Z.
\]
By determinant normalization, $c(\tau)$ and the scalar part
of $\widetilde C_Z$ are removed. This proves \eqref{ru:eq:gram} for any
sequence $\tau_j\to\tau$. By compactness of the time interval and
this sequential characterization, we have uniform convergence.
\end{proof}

\subsection{Countable spectra and the common limit}

The extremal generator need not have rational eigenvalues in a
real class. For uniqueness, it suffices to put all its possible
spectra in one countable set. We then use the discrete positive
frame updates, so avoid using the differentiability statement about the
limiting Gram matrices.

Write the original class in terms of the fixed auxiliary bundles:
\[
 \Omega=\sum_\nu\alpha_\nu c_1(A_\nu),\qquad
 \mathbb F\coloneqq\Q(\alpha_1,\ldots,\alpha_d,\pi).
\]
Here $\mathbb F\subset\R$ is the smallest field containing $\Q$,
the coefficients $\alpha_1,\ldots,\alpha_d$, and $\pi$.
Its elements are exactly the numbers
\[
 \frac{P(\alpha_1,\ldots,\alpha_d,\pi)}
      {R(\alpha_1,\ldots,\alpha_d,\pi)},
 \qquad R(\alpha_1,\ldots,\alpha_d,\pi)\ne0,
\]
where $P$ and $R$ are polynomials with rational coefficients.
There are only countably many such polynomial pairs, so
$\mathbb F$ is countable even if some generators are
transcendental. 
Adjoining $\pi$ accommodates the $2\pi$ factors in our curvature
normalization.

The auxiliary bundles $A_\nu$ and the positive coefficients
$\alpha_\nu$ are fixed once and for all in
Lemma~\ref{ru:lem:embeddings}, before any subsequence is extracted.
For every geometric limit $Z$, their marked limiting classes
satisfy $\Omega_Z=\sum_\nu\alpha_\nu c_1(A_{\nu,Z})$ with
these same coefficients. Thus $\mathbb F$ is independent of
the subsequence and its limit. The proof below places all
limiting spectra in this one fixed countable field.

\begin{lemma}[A common real Hermitian generator]
\label{ru:lem:spectrum}
All tuples $(C_{\nu,Z})_\nu$ in Lemma~\ref{ru:lem:gram} are
conjugate to one tuple $(\Lambda_\nu)_\nu$ by
$K\coloneqq\prod_\nu U(V_\nu^*)$.
\end{lemma}

\begin{proof}
We choose on each limit an effective Hamiltonian isometry torus
containing its extremal Killing field, and lift a finite cover
to every $A_{\nu,Z}$. Its section weights are integral.
For an integral torus basis $\xi_1,\ldots,\xi_\ell$, let $f_a$ be Hamiltonians
for $\omega_Z$, and put $u_a=f_a-\Vol(Z)^{-1}\int_Zf_a\,dV_Z$.
The multivariable version of the equivariant
Riemann--Roch calculation in
\cite[Lemma~4.2, Steps~2--3]{Polarized}, based on the
equivariant orbifold index formula
\cite[Section~3, Theorem~3.4]{Meinrenken}, gives
\begin{equation}\label{ru:eq:field-system}
 Q_{ab}\coloneqq\int_Zu_au_b\,dV_Z\in\mathbb F,
 \qquad
 F_a\coloneqq\int_Zu_a(S_Z-\bar S)\,dV_Z\in\mathbb F.
\end{equation}
Here is why evaluation at a real class is legitimate. 

We apply the
equivariant index formula to
$\bigotimes_\nu A_{\nu,Z}^{k_\nu}$, with all $k_\nu$ large
positive integers, and take zero, one and two torus derivatives.
The alternating weight sums in this equivariant Euler
characteristic are polynomials in the $k_\nu$ with rational
coefficients, by repeated integer interpolation. The auxiliary
bundles are Cartier, so their finite isotropy characters are trivial.
No periodic dependence on $k_\nu$ remains. This arithmetic
step does not require simultaneous higher-cohomology vanishing.
We use the degree-two and degree-one parts without torus
derivatives, the degree-three and degree-two parts with one
derivative, and the degree-four part with two derivatives.
These are the mixed versions of $a_0,a_1,b_{0a},b_{1a},d_{ab}$
in \cite[Lemma~4.2]{Polarized}. In particular the Futaki vector
uses the subleading degree-two part of the first derivative,
including its Todd contribution. The isolated nonidentity
sectors have degree at most the number of torus derivatives,
as in \cite[equation~(4.2)]{Polarized}, and do not affect any of
these specified coefficients. Thus all their mixed coefficients
are rational.

Equivariant Chern--Weil theory
\cite[Section~7, Examples~7.1--7.2]{MeinrenkenCartan}
expresses the uncentered integrals for
$\sum_\nu x_\nu c_1(A_{\nu,Z})$
through the corresponding mixed polynomial coefficients.
The constructions use invariant connections and apply
in the orbifold charts. This identity holds on the open cone of positive
combinations, and hence as a polynomial identity. It can
therefore be evaluated at the positive real coefficients
$x_\nu=\alpha_\nu$ of \eqref{rd:eq:positive-basis}.
Centering divides only by the nonzero volume.
We keep the Todd contribution throughout this polynomial
calculation. After evaluation at the actual K\"ahler form
$\omega_Z$, its additional term is the integrated Laplacian
of the corresponding Hamiltonian and is zero, exactly as in
\cite[Lemma~4.2]{Polarized}. This proves
\eqref{ru:eq:field-system}. The fixed $2\pi$ conventions are
covered by adjoining $\pi$ to $\mathbb F$.

The matrix $Q$ is positive definite because the torus is
effective. If $v$ denotes the coefficients of the extremal
Killing field, its Hamiltonian is $\kappa(S_Z-\bar S)$ for a fixed
$\kappa\in\Q(\pi)\subset\mathbb F$.
With the Hamiltonians normalized to have zero mean, this means
\[
 \sum_a v_a u_a=\kappa(S_Z-\bar S).
\]
Multiplying by $u_b$ and integrating over $Z$ gives
$Qv=\kappa F$. Since $Q$ is invertible and has entries in
$\mathbb F$, we obtain
\[
 v=\kappa Q^{-1}F\in\mathbb F^\ell.
\]
On every auxiliary section space, choose a simultaneous weight
basis for the torus action. The integral torus generators have
integral diagonal weights in this basis. Subtracting their mean
weights gives rational diagonal entries. The trace-free generator
is a linear combination of these diagonal matrices with
coefficients $v_a\in\mathbb F$, so its ordered spectrum belongs
to $\mathbb F^{\dim V_\nu}$. A trivial torus gives the
zero tuple and needs no linear system.

For arbitrary real $t$, we take the ordered logarithmic generalized
eigenvalues of $H_\nu(t+1)$ relative to $H_\nu(t)$, multiply
by $-1/2$, and subtract their mean. We combine the finitely many
lists into a vector $\sigma(t)$. This is a continuous,
basis-independent function of $t$. Lemma~\ref{ru:lem:gram}
implies that it is relatively compact at infinity and that each
cluster point is a tuple of the preceding ordered spectra.
Its cluster set is connected: it is the intersection of the
nested compact connected closures of its continuous tails.
It lies in the countable set
$\prod_\nu\mathbb F^{\dim V_\nu}$, so it is a single point.
The spectral theorem for Hermitian matrices proves the assertion.
\end{proof}

\begin{proof}[Proof of Theorem~\ref{ru:thm:unique}(1) for a fixed flow]
\textit{Step 1. Verify the finite-dimensional dynamics.}
In each block, we choose an $H_\nu(0)$-orthonormal basis and update
at integer times by the positive inverse square root of its
next Gram matrix. Write $s_n$ for the basis at time $n$ in one
block and $G_n(\tau)$ for its Gram matrix at time $n+\tau$.
Write $s_{\nu,n}=A_{\nu,n}s_{\nu,0}$ and
$A_n\coloneqq(A_{\nu,n})_\nu$ for the simultaneous changes of
evaluation coordinates. Scalar matrices do not change the
projective embeddings.
Our determinant normalization makes $\det G_n(\tau)=1$.
Along any subsequence, Lemmas~\ref{ru:lem:gram} and
\ref{ru:lem:spectrum} give, after extraction,
\[
 G_n(1)\longrightarrow e^{-2C},\qquad
 G_n(2)\longrightarrow e^{-4C},\qquad
 C=k\Lambda_\nu k^{-1}.
\]
Since $s_{n+1}=G_n(1)^{-1/2}s_n$, its next Gram matrix is
\begin{equation}\label{ru:eq:two-step-gram}
 G_n(1)^{-1/2}G_n(2)G_n(1)^{-1/2}
       \longrightarrow e^{-2C}.
\end{equation}
Thus the two consecutive positive frame updates both tend to
$e^C$. This is condition $(*)$ of \cite[Section~2]{CSW},
simultaneously in the product group
$\mathbb G\coloneqq\prod_\nu GL(V_\nu^*)$.

We use a projective Hilbert embedding for the product ambient space
and its Segre polarization. The finite-dimensional arguments
of \cite[Lemma~2.6, Proposition~2.7 and Proposition~2.10]{CSW}
then apply. More explicitly, we perform the filtration splittings
in each $V_\nu^*$ separately. Their direct sums, and the
near-identity corrections aligning their summands, belong to
$\mathbb G$. Coincident eigenvalues in different factors do
not require mixing those factors. The convexity and growth
arguments apply to any finite-dimensional representation of
this product. We obtain a nonempty compact connected compatible
limit set $\mathscr C$, invariant under
\[
 K_\Lambda\coloneqq
 \prod_\nu\operatorname{Cent}_{U(V_\nu^*)}(\Lambda_\nu),
 \qquad
 G_\Lambda\coloneqq
 \prod_\nu\operatorname{Cent}_{GL(V_\nu^*)}(\Lambda_\nu),
\]
with the following explicit meaning. Let $[M]$ be the product
Hilbert point defined by the initial complete bases, so that
the integer-time point is $A_n[M]$. We choose the block-unitary
gauges $u_n$ in condition $(*)'$ of \cite[Section~2]{CSW},
for which
\[
 u_{n+1}^{-1}u_n\longrightarrow I,\qquad
 u_{n+1}^{-1}A_{n+1}A_n^{-1}u_n\longrightarrow e^\Lambda.
\]
The first condition is part of $(*)'$ and is used in the
connectedness argument of \cite[Lemma~2.6]{CSW}. The conjugated
frame-update limit alone is not the hypothesis being applied.
Then $\mathscr C$ is the $K_\Lambda$-saturation of the set of
subsequential limits of $u_n^{-1}A_n[M]$.
The blockwise splitting gives $C_0\in\mathbb G$ and the common
intermediate point
\[
 [\overline M]\coloneqq
       \lim_{s\to+\infty}e^{s\Lambda}C_0[M],\qquad
 \mathscr C\subset\overline{G_\Lambda[\overline M]}.
\]
This is a real exponential limit. Every point of $\mathscr C$ represents a geometric
limit with all its auxiliary
bundles, embedded by complete orthonormal section spaces.

\textit{Step 2. Verify the orbit-uniqueness hypotheses.}
The stabilizer of a limit in $G_\Lambda$, modulo the scalar
torus of the factors, is the group of automorphisms preserving
all auxiliary bundles and commuting with the extremal field.
Its identity component is the extremal centralizer: connected
automorphisms preserve all the bundles because
$\operatorname{Pic}^0(Z)=0$. The analytic centralizer argument
of \cite[Lemma~4.4]{Polarized} applies in the real K\"ahler class,
using $H^{0,1}(Z)=0$ to obtain holomorphy potentials. It makes
this identity component reductive. The full stabilizer is an
algebraic group with finitely many components. Its unipotent
radical is therefore trivial as well. The scalar extension
does not change reductivity.

Suppose $x_2=gx_1$ are points of $\mathscr C$, with
$g\in G_\Lambda$. The induced orbifold isomorphism preserves
every auxiliary class, and hence preserves
$\Omega_Z=\sum_\nu\alpha_\nu c_1(A_{\nu,Z})$.
It also intertwines the extremal fields. We apply
Theorem~\ref{un:thm:fixed} in this common real K\"ahler class.
It corrects $g$ by an automorphism in the identity component of
the extremal centralizer to a holomorphic isometry $F$.
That correction preserves every auxiliary bundle. The pullback
by $F$ identifies their harmonic curvatures, and therefore
their Hermitian metrics up to positive constants. Rescaling
the bundle identifications makes the pullback maps on all
complete section spaces unitary. Their projective actions
commute with the extremal generators. Taking traces removes
any scalar commutators. We have proved
\[
 \mathscr C\cap G_\Lambda x=K_\Lambda x.
\]
\cite[Lemma~4.8]{Polarized}, equivalently the slice argument in
\cite[proof of Theorem~2.12]{CSW}, gives one $K_\Lambda$-orbit.
The same metric correction compares two metrics represented
by a single Hilbert point. Hence all integer-time limits are
holomorphically isometric as extremal orbifolds.

\textit{Step 3. Return to arbitrary times.}
For $t_j=n_j+\tau_j$, with $n_j$ integral and $0\leq\tau_j<1$,
extract $\tau_j\to\tau$ and apply
Theorem~\ref{app:thm:compactness} on translated time intervals
centred at $n_j$. The real-time limit is
$F_\tau^*\omega_Z$, which is holomorphically isometric to
$\omega_Z$. Thus every subsequence has the same limit up to
holomorphic isometry. Relative compactness then gives
convergence of the whole flow modulo holomorphic isometry in
the sense of Definition~\ref{app:def:geometric-convergence}.
Indeed, a sequence remaining outside a neighborhood of the
limiting family and all its holomorphically isometric copies
would have a subsequence converging into that neighborhood.
This proves the
fixed-flow assertion, including uniqueness of the translated
self-similar windows up to the same holomorphic isometries.
Independence of the initial metric is proved next.
\end{proof}

\subsection{Independence of the initial metric}
\label{ru:sec:two-flow}

We now compare two flows in the same real class. The bounded-geometry
polarized case was proved in \cite[Section~4, concluding Remark~(1),
p.~3170]{CSW}, using Mabuchi-distance contraction and estimates for the
difference of the K\"ahler potentials. Here a positive combination of
auxiliary classes converts an endpoint integral estimate for that
difference into bounds for all relative Gram matrices. This avoids
higher derivative estimates comparing the two metrics near their
exceptional regions.

\begin{theorem}[Comparison of two flows]\label{ru:thm:two-flow}
Fix $(M,J)$ and $\Omega\in\mathcal U$. Let $\omega_i(t)\in2\pi\Omega$,
$i=1,2$, be global smooth Calabi flows, each satisfying
\eqref{tc:eq:sob} and the compactness property of
Theorem~\ref{app:thm:compactness}. Their limiting extremal orbifolds
are holomorphically isometric. In particular, the limit in
Theorem~\ref{ru:thm:unique}(1) is independent of the initial metric
among all metrics satisfying \eqref{ru:eq:barrier}.
\end{theorem}

\begin{proof}
\textit{Step 1. An endpoint estimate for the potentials.}
Choose $\omega_*\in2\pi\Omega$ and write
$\omega_i(t)=\omega_*+\ddc\varphi_i(t)$. We normalize the initial
potentials to have Aubin--Mabuchi energy zero, and choose
$\partial_t\varphi_i=S_i-\bar S$. This normalization is preserved,
because its derivative is $V^{-1}\int_M(S_i-\bar S)dV_{i,t}=0$.
Let $d_p$, $p=1,2$, denote the path metric induced by
$\|\xi\|_{p,\varphi}=(V^{-1}\int_M|\xi|^p dV_\varphi)^{1/p}$.
The distance contraction of Calabi--Chen \cite[Section~1.2]{CalabiChen}
gives
\[
 d_2(\varphi_1(t),\varphi_2(t))
 \leq d_2(\varphi_1(0),\varphi_2(0)).
\]
For a formulation requiring only the existence of these two smooth
flows, we use the contractive weak gradient flow and its agreement with
the smooth flow in \cite[Theorem~2.12 and Proposition~6.1]{BDLFlow}.
From H\"older's inequality, we get $d_1\leq d_2$. Darvas's endpoint comparison
\cite[Theorem~3]{DarvasFinite}, also stated in
\cite[Theorem~2.4, equation~(8)]{BDLFlow}, therefore implies
\begin{equation}\label{ru:eq:two-flow-endpoint}
 \int_M|\delta_t|(dV_{1,t}+dV_{2,t})\leq D,
 \qquad \delta_t\coloneqq\varphi_2(t)-\varphi_1(t),
\end{equation}
where $D$ is independent of $t$. These metric-space results apply to
arbitrary real K\"ahler classes. The fixed volume factor is absorbed
into $D$.

\textit{Step 2. Common auxiliary bundles and the scalar ambiguity.}
We apply Lemmas~\ref{ru:lem:ample}--\ref{ru:lem:embeddings} to both flows.
We intersect their ample neighborhoods of $\Omega$, choose the rational
simplex inside that intersection, and take common divisible powers.
The resulting fixed bundles $A_\nu$ have complete-section convergence
for both flows and satisfy
\[
 \Omega=\sum_\nu\alpha_\nu c_1(A_\nu),\qquad \alpha_\nu>0.
\]
All the required bounds hold for the union of the two flows' limit
sets by taking the larger of their constants. Taking powers divides
the coefficients by positive integers and preserves their positivity.

Let $h_{\nu,i}(t)$ be the harmonic Chern metric on $A_\nu$ for
$\omega_i(t)$. Write its real curvature as $\Theta_{\nu,i}(t)$.
Uniqueness of harmonic representatives gives
\[
 \Theta_{\nu,i}(t)\coloneqq-\ddc\log h_{\nu,i}(t),\qquad
 \sum_\nu\alpha_\nu\Theta_{\nu,i}(t)=\omega_i(t).
\]
Define $H_{\nu,i}(t)(s,s')\coloneqq\int_Mh_{\nu,i}(t)(s,s')dV_{i,t}$.
Fix one source basis of each $H^0(M,A_\nu)$ and use the scalar
freedom in $h_{\nu,i}(t)$ to normalize $\det H_{\nu,i}(t)=1$
in that basis, for both flows. For the ratios of these Hermitian
metrics set
\[
 u_{\nu,t}\coloneqq-\log\frac{h_{\nu,2}(t)}{h_{\nu,1}(t)}.
\]
The curvature identity implies
\begin{equation}\label{ru:eq:two-flow-scalar}
 \sum_\nu\alpha_\nu u_{\nu,t}=\delta_t+c_t
\end{equation}
for a spatial constant $c_t$: the difference of the two sides is
pluriharmonic on compact connected $M$. We will cancel this constant,
without assuming it is bounded.

\textit{Step 3. Bound the relative Gram matrices.}
Diagonalize $H_{\nu,2}(t)$ relative to $H_{\nu,1}(t)$. In an
$H_{\nu,1}(t)$-orthonormal eigenbasis $(s_{\nu,a})$, we write
\[
 H_{\nu,2}(t)(s_{\nu,a},s_{\nu,b})=d_{\nu,a}(t)\delta_{ab},
 \qquad \lambda_{\nu,a}(t)\coloneqq d_{\nu,a}(t)^{-1}.
\]
Let $\rho_{\nu,i}(t)$ denote the Bergman sum for $A_\nu$ with
$h_{\nu,i}(t)$ and $dV_{i,t}$. The section estimate
\eqref{ru:eq:moser} bounds these sums above by constants $B_\nu$.
The second sum, expressed in the first eigenbasis, is
\begin{equation}\label{ru:eq:two-flow-bergman}
 \rho_{\nu,2}(t)
 =e^{-u_{\nu,t}}\sum_a\lambda_{\nu,a}(t)
             |s_{\nu,a}|_{h_{\nu,1}(t)}^2\leq B_\nu.
\end{equation}
We claim that every relative condition number is uniformly bounded:
\begin{equation}\label{ru:eq:two-flow-condition}
 \sup_{t\geq0}
 \frac{\lambda_{\nu,\max}(t)}{\lambda_{\nu,\min}(t)}<\infty.
\end{equation}

If this failed, choose $t_j$ along which one ratio tends to infinity.
The Gram forms vary smoothly on finite time intervals, so after
extraction $t_j\to\infty$. Extract geometric limits of both flows
and their complete section spaces along this same sequence.
For each $\nu$, choose a unit eigensection of the first Gram form
corresponding to $\lambda_{\nu,\max}(t_j)$. Complete-section
convergence, and a subsequence of its unit coefficient vectors,
give a nonzero limiting section on the first limit $Z_1$. Its
limiting $L^2$ norm is one. The connected normal space $Z_1$ is
irreducible. Thus a regular ball avoids the zero divisors of all
these finitely many sections. On the corresponding source sets $E_j$
we have, for some $\varepsilon,v>0$ and all sufficiently large $j$,
\[
 |s_{\nu,\max}|_{h_{\nu,1}(t_j)}^2\geq\varepsilon
 \quad\text{for every }\nu,\qquad
 \Vol_{\omega_1(t_j)}(E_j)\geq v.
\]
Equation~\eqref{ru:eq:two-flow-bergman} gives
$u_{\nu,t_j}\geq\log\lambda_{\nu,\max}(t_j)-C$ on $E_j$.
Since every $\alpha_\nu$ is positive, summing and using
\eqref{ru:eq:two-flow-scalar} yields
\[
 \delta_{t_j}+c_{t_j}
 \geq\sum_\nu\alpha_\nu\log\lambda_{\nu,\max}(t_j)-C
 \quad\text{on }E_j.
\]
Integrate over $E_j$ and use \eqref{ru:eq:two-flow-endpoint} to obtain
\begin{equation}\label{ru:eq:two-flow-max}
 \sum_\nu\alpha_\nu\log\lambda_{\nu,\max}(t_j)-c_{t_j}
 \leq C+D/v.
\end{equation}
We reverse the flows and repeat on a regular ball of the second limit.
The relative eigenvalues are reciprocated, while $\delta_t$ and
$c_t$ both change sign. Consequently,
\begin{equation}\label{ru:eq:two-flow-min}
 -\sum_\nu\alpha_\nu\log\lambda_{\nu,\min}(t_j)+c_{t_j}
 \leq C'.
\end{equation}
The two regular balls need not correspond. Adding these inequalities
cancels the scalar ambiguity and gives
\[
 \sum_\nu\alpha_\nu
 \log\frac{\lambda_{\nu,\max}(t_j)}{\lambda_{\nu,\min}(t_j)}
 \leq C''.
\]
Each summand is nonnegative and each coefficient is positive. This
contradicts the chosen sequence and proves
\eqref{ru:eq:two-flow-condition}. The nonvanishing constants above
need only hold on the extracted contradictory sequence. 

Determinant normalization makes the product of the relative
eigenvalues equal to one. The condition-number bounds therefore give
one $C_{\mathrm{rel}}\geq1$, independent of time, such that
\begin{equation}\label{ru:eq:two-flow-forms}
 C_{\mathrm{rel}}^{-1}H_{\nu,1}(t)\leq H_{\nu,2}(t)
 \leq C_{\mathrm{rel}}H_{\nu,1}(t)
 \quad\text{for every }\nu\text{ and }t\geq0.
\end{equation}
These are inequalities of Hermitian forms. 

\textit{Step 4. Identify the limiting spaces and extremal fields.}
Take simultaneous translated limits at $t_j\to\infty$, using a
diagonal extraction over all bounded time intervals. In one auxiliary
factor choose ordered orthonormal bases $a_j$ and $b_j$ for the
respective forms at time $t_j$. Regard them as columns of sections,
as above, and define $T_j$ by $b_j=T_ja_j$.
Equation~\eqref{ru:eq:two-flow-forms} bounds $T_j$ and its inverse.
After extraction $T_j\to T\in GL(N_\nu,\C)$.
Doing this in every factor, complete-section and Hilbert convergence
show that these projective changes of evaluation coordinates identify the two limits
and every auxiliary bundle. Thus the resulting biholomorphism
$F:Z_1\to Z_2$ preserves their marked real metric classes.

This is also an orbifold biholomorphism. Near an isolated quotient
point the punctured uniformizing ball is the universal cover of
the punctured quotient, since the local group acts freely away
from the origin. A local coarse-space isomorphism lifts on these
covers. Hartogs extension of the lift and its inverse extends the
isomorphism across the origins. These are precisely our effective
orbifold charts.

The same $T$ intertwines the extremal fields. Lemma~\ref{ru:lem:gram}
and the determinant normalization give
\[
 \operatorname{Gram}_{H_{\nu,1}(t_j+\tau)}(a_j)
       \longrightarrow e^{-2\tau C_1},\qquad
 \operatorname{Gram}_{H_{\nu,2}(t_j+\tau)}(b_j)
       \longrightarrow e^{-2\tau C_2},
\]
where $C_1,C_2$ are trace-free Hermitian. The frames are fixed at
the base time throughout each translated window. Passing
\eqref{ru:eq:two-flow-forms} to the limit in the basis $b_j$ gives
\begin{equation}\label{ru:eq:two-flow-windows}
 C_{\mathrm{rel}}^{-1}e^{-2\tau C_2}
 \leq T e^{-2\tau C_1}T^*
 \leq C_{\mathrm{rel}}e^{-2\tau C_2},\qquad \tau\in\R.
\end{equation}
Here we used $b_j=T_ja_j$ and the Gram transformation rule above.
The same $C_{\mathrm{rel}}$ works for all $\tau$, because the original estimate
holds at every time. Conjugating this inequality by $e^{\tau C_2}$
bounds the positive matrix $U(\tau)U(\tau)^*$ above and below, where
$U(\tau)\coloneqq e^{\tau C_2}Te^{-\tau C_1}$.
Thus $U(\tau)$ and its
inverse are bounded for all positive and negative $\tau$.
In eigenbases its $(a,b)$-entry is
\[
 T_{ab}e^{\tau(\gamma_{2,a}-\gamma_{1,b})},
\]
where $\gamma_{i,a}$ is the corresponding eigenvalue of $C_i$.
Boundedness in both directions forces $T_{ab}=0$ unless
$\gamma_{1,b}=\gamma_{2,a}$. Hence
\begin{equation}\label{ru:eq:two-flow-intertwine}
 C_2T=TC_1.
\end{equation}
The projective actions of these matrices are the extremal motions.
Their intertwining in the complete embeddings shows that $F$
intertwines the extremal fields on $Z_1$ and $Z_2$.

\textit{Step 5. Apply fixed-orbifold uniqueness.}
The metrics $\omega_{Z_1}$ and $F^*\omega_{Z_2}$ are now extremal
metrics in the same real class with the same extremal field.
Lemma~\ref{ru:lem:embeddings} gives $q(Z_1)=0$ and an auxiliary
polarization, so Theorem~\ref{un:thm:fixed} supplies
$a\in\Aut^0(Z_1)$ with
$a^*F^*\omega_{Z_2}=\omega_{Z_1}$. Thus $F\circ a$ is the
required holomorphic isometry.

By the fixed-flow argument above, we identify the simultaneous limits
with the unique limits of the respective flows. It applies under
the present global existence and Sobolev hypotheses: the energy
barrier was used only to supply those hypotheses. This proves the
comparison theorem and completes Theorem~\ref{ru:thm:unique}(1).
\end{proof}

The comparison constants may depend on the pair of initial metrics.
Pairwise comparison is sufficient for independence of the initial
metric. 
The class and complex structure remain fixed throughout this proof.
We treat the continuity under variation of the class separately in
Section~\ref{real:sec:existence}.

\subsection{Degeneration preserving the original real class}
\label{rd:sec:degeneration}

The common auxiliary embeddings retain more than the complex
limit: they retain every class in the positive expression
\eqref{rd:eq:positive-basis}. We use this observation to extend
the stability argument to the original, possibly irrational,
class. The rational approximation is used only for a generator in
a fixed finite-dimensional torus. The class of the flow remains
unchanged.

\begin{proposition}[Class-preserving degeneration and its invariant]
\label{rd:prop:degeneration}
Under the assumptions of Theorem~\ref{ru:thm:unique}, let
$(Z,\omega_Z,\Omega_Z)$ be a geometric limit, with
$[\omega_Z]=2\pi\Omega_Z$. There is a normal relatively
K\"ahler cohomological test configuration $(\mathcal X,\mathcal A)$
for $(M,\Omega)$ whose central fibre, including its restricted
class, is $(Z,\Omega_Z)$. Its central fibre is reduced,
irreducible and normal. If $\omega_Z$ is not cscK, the
configuration can be chosen with $\DF(\mathcal X,\mathcal A)<0$.
If $\omega_Z$ is cscK, it can be chosen with
$\DF(\mathcal X,\mathcal A)=0$.
\end{proposition}

\begin{proof}
\textit{Step 1. Keep all auxiliary polarizations.}
We use the bundles $A_\nu$ from Lemma~\ref{ru:lem:embeddings},
chosen with the positive coefficients in
\eqref{rd:eq:positive-basis}. The common Hilbert scheme is that
of subschemes of $\prod_\nu\mathbb P(V_\nu^*)$ with the fixed
Segre Hilbert polynomial. Its universal family retains the
restrictions of each factor bundle $\mathcal O(1)$.
Choose an equivariant projective Hilbert embedding into
$\mathbb P(E)$, for a finite-dimensional representation $E$
of $\mathbb G=\prod_\nu GL(V_\nu^*)$.
The proof of Theorem~\ref{ru:thm:unique}(1) supplies
\begin{equation}\label{rd:eq:orbit-data}
 x\coloneqq C_0[M],\qquad
 y\coloneqq\lim_{t\to+\infty}e^{t\Lambda}x,\qquad
 p\coloneqq[Z]\in\overline{G_\Lambda y}.
\end{equation}
The stabilizer of $p$ in $G_\Lambda$ is reductive, as verified
in Step~2 of that proof. It suffices to treat integer-time limits. Step~3 of the
fixed-flow proof, together with the bundle lifts constructed
in Lemma~\ref{ru:lem:gram}, identifies every arbitrary-time
limit with an integer-time limit by a holomorphic isometry
preserving every auxiliary bundle.

\textit{Step 2. Replace the real generator without changing its limit.}
Suppose first that $Z$ is not cscK. Let
\[
 T\coloneqq\overline{\{e^{\sqrt{-1}s\Lambda}:s\in\R\}}
       \subset\prod_\nu U(V_\nu^*).
\]
This is a compact torus. Its complexification $T^{\C}$ is an
algebraic torus contained in the center of $G_\Lambda$:
each of its elements is scalar on each eigenspace of each
$\Lambda_\nu$. We identify the real cocharacter space of $T^{\C}$
with the Hermitian generators $H$ for which
$e^{\sqrt{-1}sH}\in T$. It contains $\Lambda$, and its rational
cocharacters are dense.

Then we decompose a lift of $x$ into the finitely many occurring
$T^{\C}$-weight spaces,
\[
 v=\sum_{\chi\in\mathcal S}v_\chi,\qquad v_\chi\ne0.
\]
We see that distinct characters of $T^{\C}$ have distinct values on
$\Lambda$. Indeed, if $(\chi-\psi)(\Lambda)=0$, the character
$\chi\psi^{-1}$ is trivial on $e^{\sqrt{-1}\R\Lambda}$,
hence on its closure $T$, and consequently on $T^{\C}$.
Thus there is a unique character $\chi_0\in\mathcal S$
maximizing $\chi(\Lambda)$, and $y=[v_{\chi_0}]$.
The inequalities
\begin{equation}\label{rd:eq:weight-chamber}
 \chi_0(H)<\chi(H)\quad
       (\chi\in\mathcal S\setminus\{\chi_0\})
\end{equation}
define an open rational chamber containing $-\Lambda$.
Every rational $H$ in this chamber, after multiplication by
an integer $b>0$ clearing its denominators, defines
\[
 \lambda(z)\coloneqq\exp(bH\log z),\qquad
 \lim_{z\to0}\lambda(z)x=y.
\]
The limit is exactly $[v_{\chi_0}]$, not merely a point in its
orbit. The entire weight component, with its coefficients, is
retained. The subgroup is single-valued because $bH$ is an
integral cocharacter. It fixes $y$ and, being central in
$G_\Lambda$, fixes its entire orbit closure, including $p$.

We also choose $H$ to preserve the strict Futaki sign. For a
Hermitian generator $D$ in this torus let $\theta_D$ be the
real potential, defined up to a constant, satisfying
\[
 \left.\frac{d}{ds}\right|_{s=0}
          \exp(-sD)^*\omega_Z=\ddc\theta_D.
\]
These Hamiltonian potentials exist since
$H^{0,1}(Z)=0$ and the imaginary-time torus acts by isometries.
Define the central Futaki character, with our normalization, by
\begin{equation}\label{rd:eq:real-futaki}
 \mathfrak F_{\Omega_Z}(D)
 \coloneqq-\frac{1}{4(2\pi)^2}
       \int_Z\theta_D(S_Z-\bar S)\,dV_Z.
\end{equation}
It is linear and continuous in $D$. Changing the additive
constant in $\theta_D$ has no effect. The action of
$e^{t\Lambda}$ on $Z$ is its extremal Calabi flow. Therefore
$\theta_{-\Lambda}=S_Z-\bar S$ up to a constant, and
\begin{equation}\label{rd:eq:negative-real}
 \mathfrak F_{\Omega_Z}(-\Lambda)
       =-\frac{\Ca(\omega_Z)}{4(2\pi)^2}<0.
\end{equation}
Choose a rational $H$ sufficiently close to $-\Lambda$ to
satisfy both \eqref{rd:eq:weight-chamber} and
$\mathfrak F_{\Omega_Z}(H)<0$. Multiplication by $b>0$
preserves this sign. We write
$\mathfrak F_{\Omega_Z}(\lambda)$ for the value at $bH$.
This step approximates the generator, not $\Omega_Z$.

\textit{Step 3. Reach the metric limit by one algebraic subgroup.}
Apply the reductive-stabilizer accessibility argument of
\cite[Proposition~1(2), proof in Section~5.2]{Donaldson} to
$p\in\overline{G_\Lambda y}$. As explained in the paragraph preceding
\cite[Lemma~4.5]{Polarized}, this slice argument does not require
the stabilizer of $y$ to be finite. It supplies
$g\in G_\Lambda$ and an algebraic one-parameter subgroup
$\mu:\C^*\to G_\Lambda$ such that
\[
 \lim_{z\to0}\mu(z)gy=p.
\]
If $p\in G_\Lambda y$, take $\mu$ trivial. Since $\lambda$
is central in $G_\Lambda$, it commutes with $g$ and $\mu$.
The commuting-degeneration lemma
\cite[Lemma~4.5]{Polarized}, applied in the representation $E$, gives,
for every sufficiently large integer $m$,
\begin{equation}\label{rd:eq:combined}
 \nu_m(z)\coloneqq\lambda(z)^m\mu(z),\qquad
 \lim_{z\to0}\nu_m(z)gx=p.
\end{equation}
Both $\lambda$ and $\mu$ fix $p$, so both induce actions on
$Z$ and its auxiliary bundles.

If $Z$ is cscK, its extremal field vanishes. All trace-free
auxiliary generators are then zero, so $\Lambda=0$ and
$y=x$. Apply the same accessibility argument directly to
$p\in\overline{\mathbb Gx}$ and use $\nu=\mu$. 

\textit{Step 4. Construct the flat family and its exact real class.}
Pull back the Hilbert universal family along the algebraic
orbit map in \eqref{rd:eq:combined}, extended to $z=0$ by
$p$. This yields a flat projective family
\[
 \mathcal X_m\subset
       \prod_\nu\mathbb P(V_\nu^*)\times\mathbb A^1,
 \qquad \pi_m:\mathcal X_m\to\mathbb A^1.
\]
Its central fibre is the scheme $Z$, and every nonzero fibre
is the image of $M$ by $\nu_m(z)gC_0$. The action
\[
 a\cdot(q,z)=(\nu_m(a)q,az)
\]
preserves the family and lifts to each restricted factor
bundle, denoted by $\mathcal L_{\nu,m}$. The equivariant
trivialization over $\C^*$ identifies these bundles with
$\operatorname{pr}_M^*A_\nu$. At zero they restrict to
$A_{\nu,Z}$. Every fibre is geometrically normal, and the base
is smooth. Flatness therefore implies normality of the total
space, by \cite[Lemma~10.163.8, Tag~0C22]{Bertini}. Alternatively,
the $R_1+S_2$ argument in \cite[Proposition~4.6]{Polarized} applies.
In particular no normalization changes the central fibre.

Set
\begin{equation}\label{rd:eq:total-class}
 \mathcal A_m\coloneqq
       \sum_\nu\alpha_\nu c_1(\mathcal L_{\nu,m}).
\end{equation}
The marking and \eqref{rd:eq:positive-basis} give the exact
restrictions
\[
 \mathcal A_m|_{\mathcal X_{m,1}}=\Omega,\qquad
 \mathcal A_m|_{\mathcal X_{m,0}}=\Omega_Z,
\]
under the specified identifications. The linearizations make
$\mathcal A_m$ invariant. Since every $\alpha_\nu>0$, the
weighted sum of the factor Fubini--Study forms is positive on
the ambient product. Its restriction is a relative K\"ahler
representative on this family. This uses the usual notion of
a K\"ahler form on a complex space, obtained from ambient local
embeddings. It does not identify orbifold-smooth potentials
with smooth potentials on the coarse space.

For completeness, compactify over $\mathbb P^1$ using the
action trivialization, and glue each auxiliary bundle to
$\operatorname{pr}_M^*A_\nu$ on the chart at infinity. Write
$\overline{\mathcal L}_{\nu,m}$ for the resulting bundle.
Its restriction to every fibre is ample. The fibrewise
criterion for relative ampleness
\cite[Lemma~37.50.3, Tag~0D2S]{Bertini}, equivalently EGA
IV$_3$, Corollary~9.6.4, makes each
$\overline{\mathcal L}_{\nu,m}$ relatively ample. For a
sufficiently large integer $N$, all
$\overline{\mathcal L}_{\nu,m}\otimes\pi_m^*\mathcal O(N)$
are ample on the compactified total space. Their positive
real combination is a K\"ahler class. Consequently
$\overline{\mathcal A}_m+N(\sum_\nu\alpha_\nu)
\pi_m^*c_1(\mathcal O(1))$ is K\"ahler. This proves the
required relative K\"ahlerness of the exact real class.

\textit{Step 5. Identify the invariant on this fixed family.}
We justify the real-class version of the central-weight
calculation, including its sign. Keep $m$, the family and
all its linearized bundles fixed. For positive rational
coefficients $a_\nu$ close to $\alpha_\nu$, put
\[
 B(a)\coloneqq\sum_\nu a_\nu c_1(A_\nu),\qquad
 \mathcal A_m(a)\coloneqq
       \sum_\nu a_\nu c_1(\mathcal L_{\nu,m}).
\]
Write $B_Z(a)\coloneqq\sum_\nu a_\nu c_1(A_{\nu,Z})$
for the corresponding central class.
Choose an integer $q>0$ clearing the denominators. The bundle
$\bigotimes_\nu\mathcal L_{\nu,m}^{qa_\nu}$ is relatively
ample, and its central restriction is the corresponding
ample Cartier bundle on $Z$. The equivariant orbifold Riemann--Roch calculation below
identifies its Donaldson--Futaki invariant with the central
Futaki character. Dividing by $q^2$ removes the
denominator clearing in complex dimension two.
The marked intersection identities \eqref{ru:eq:marking} give
$B_Z(a)^2=B(a)^2$ and
$c_1(Z)\cdot B_Z(a)=c_1(M)\cdot B(a)$, so the central
and general-fibre slopes and average scalar curvatures agree.

More explicitly, for an integral polarization $B$ on a surface,
the coefficients are $a_0=B^2/2$ and $a_1=c_1(M)\cdot B/2$.
Let $\overline{\mathcal B}$ denote its extended polarization
class on the compactified family. The leading weight coefficients
are $b_0=\overline{\mathcal B}^{3}/6$ and
$b_1=-K_{\mathcal Y/\mathbb P^1}\cdot
(p^*\overline{\mathcal B})^2/4$, where intersections may be
computed on a resolution. Hence
$(a_1/a_0)b_0-b_1$ is precisely
\eqref{rd:eq:cohom-df}. On the central orbifold, the
weight-to-Futaki calculation
\cite[equations~(4.8)--(4.9)]{Polarized}, together with the
equivariant index expansion used in the proof of
Lemma~\ref{ru:lem:spectrum}, gives
\eqref{rd:eq:real-futaki}. As in that proof, the isolated
nonidentity sectors do not affect the Hilbert coefficients
of degrees two and one or the weight coefficients of
degrees three and two.

Now let the rational vectors $a$ tend to $\alpha$ on this
\emph{same family}. The mixed intersection numbers are fixed.
The expression in \eqref{rd:eq:cohom-df} is continuous in $a$, its only
denominator being the positive volume $B(a)^2$. The central
equivariant Riemann--Roch formula has the same continuity,
by its mixed Chern--Weil expression. This can also be seen
using smooth orbifold representatives for nearby classes and
their Hamiltonian potentials: positivity persists near
$\omega_Z$, and the invariant is independent of the chosen
representative. Thus the rational identity passes to the
exact real class. For an algebraic central action $\eta$,
we denote the resulting character by
$\mathfrak F_{\Omega_Z}(\eta)$. It extends
\eqref{rd:eq:real-futaki}. For general holomorphic generators
one uses their complex holomorphy potentials. Such potentials
exist since $H^{0,1}(Z)=0$.
In particular $\mu$ need not act by isometries of $\omega_Z$.
For its rational-class weight calculation, average a positive
orbifold representative over the compact circle in $\mu$.
Metric independence of the Futaki character
\cite[Theorem~4.21]{SzekelyhidiExtremal} permits its evaluation
at $\omega_Z$ after passing to the real class.
The proof uses local differential identities and integration
by parts, so the same argument applies on compact
K\"ahler orbifolds.

The central actions $\lambda$ and $\mu$ commute. At rational
coefficients their section generators add, so their total
weights add at each power. Passing to the real coefficients
therefore gives the exact identity
\begin{equation}\label{rd:eq:df-sum}
 \DF(\mathcal X_m,\mathcal A_m)
   =\mathfrak F_{\Omega_Z}(\nu_m)
   =m\mathfrak F_{\Omega_Z}(\lambda)
          +\mathfrak F_{\Omega_Z}(\mu).
\end{equation}
Although the family changes with $m$, the two characters on
the right are computed on the fixed $(Z,\Omega_Z)$ and are
independent of $m$. If $Z$ is not cscK, the first is strictly
negative by Step~2. Choose $m$ sufficiently large both for
\eqref{rd:eq:combined} and for negativity in
\eqref{rd:eq:df-sum}. If $Z$ is cscK, every central Futaki
character vanishes because $S_Z-\bar S=0$. The direct
configuration constructed with $\mu$ consequently has
$\DF=0$. This proves both assertions.
\end{proof}

\begin{remark}
The rational approximation in
Proposition~\ref{rd:prop:degeneration} changes a finite-dimensional
generator while keeping the real class $\Omega$ fixed. The
resulting test configuration has the metric-limit complex space
as its central fibre, preserves all auxiliary classes, and has
the Donaldson--Futaki sign required for the stability argument.
It need not reproduce the flow-induced filtration or preserve
an optimality property. The stability argument therefore does
not require approximation of $\Omega$ by stable rational classes.
\end{remark}

\subsection{Stability and the original complex surface}

The class-preserving configuration has the sign required to apply
the stability hypothesis. The product equality case recovers both
the original complex surface and the exact real class.

\begin{proof}[Proof of Theorem~\ref{rd:thm:convergence}(2)--(3)]
Proposition~\ref{rd:prop:degeneration} constructs a configuration
preserving the exact real class and having central fibre
$(Z,\Omega_Z)$. Its total space is normal and its central fibre
is reduced, irreducible and normal. If $Z$ is not cscK, the
configuration can be chosen with negative Donaldson--Futaki
invariant. Either the full or the stated restricted semistability
condition excludes this possibility.

When $Z$ is cscK, the same proposition gives a configuration
with $\DF=0$. The product equality condition identifies its
central complex space with $M$. It also preserves the real
class: compose a product identification over $\C$ with a
constant automorphism of $M$ so that its identification at
$t=1$ agrees with the specified general-fibre marking.
The inclusions
of the central and general fibres are homotopic, so the
restriction of $\mathcal A$ to each is the same class $\Omega$.
On the compact K\"ahler fibre, the $\partial\bar\partial$-lemma
identifies this de Rham equality with equality of real
$(1,1)$-classes.
Thus there is a biholomorphism $F$ with
$F^*\Omega_Z=\Omega$. The local quotient groups act freely away
from their isolated fixed points. The Hartogs--Jacobian argument
in the proof of \cite[Corollary~4.10]{Polarized} consequently
shows that every such group is trivial when the coarse space
is smooth. Hence $F^*\omega_Z$ is a smooth cscK metric in
$2\pi\Omega$. The whole-flow convergence is exactly
Theorem~\ref{ru:thm:unique}(1).
\end{proof}

The restricted stability conditions suffice because
Proposition~\ref{rd:prop:degeneration} constructs a test configuration
whose central fibre is normal, reduced and irreducible. The
argument does not require or establish a general reduction of
arbitrary test configurations to this form.

The resulting convergence is geometric modulo holomorphic isometry,
in the sense of Definition~\ref{app:def:geometric-convergence}.
The comparison maps may depend on the chosen sequence, so this
statement does not assert convergence in a fixed coordinate gauge.

\section{Equivariant convergence}\label{eqv:sec:extensions}

We retain a prescribed finite group of symmetries in the real-class
constructions. Equivariant compactness is supplied by the compactness
paper. The linearizations and finite-dimensional centralizer arguments
were established in the polarized companion. We recall the geometric
input and then carry the auxiliary bundles and their positive real
combination through those arguments.

\subsection{Geometric compactness and fixed-orbifold uniqueness}
\label{eqv:sec:geometric}

We shall also use the finite-group equivariant refinement of the
compactness theorem, whose geometric proof is given in the companion
paper.

\begin{lemma}[Finite-group equivariant compactness]\label{eqv:lem:compactness}
Suppose $(M_j,g_j(\tau),J_j)$ converges with the time-independent
GH comparisons and regular smooth convergence of
Theorem~\ref{app:thm:compactness} on
every bounded time interval. If a fixed finite group $\Gamma$ acts
by time-independent biholomorphisms preserving every $g_j(\tau)$,
then, after a further subsequence, its actions converge to one
time-independent action on the limiting orbifold $Z$ by
biholomorphisms preserving every limiting metric. The GH
correspondences can be chosen $\Gamma$-equivariantly, and the local
action maps converge smoothly on regular compact subsets.
\end{lemma}
This is \cite[Theorem 6.12]{Compactness}. Its application to Calabi
flow is \cite[Corollary 6.13]{Compactness}. It requires no polarization.
Smooth convergence of the local action maps uses the original
comparison embeddings. The embeddings themselves are not asserted
to be exactly equivariant.

We also need the finite-group refinement of fixed-orbifold
uniqueness from Subsection~\ref{app:sec:orbifold-uniqueness}.

\begin{lemma}[Fixed-orbifold uniqueness with finite symmetry]
\label{eqv:lem:fixed}
In Theorem~\ref{un:thm:fixed}, suppose that a finite group $\Gamma$
preserves both metrics and the polarization. Then the isometry can be
chosen in $\operatorname{Cent}_{\Aut^0(Z,A,Y)}(\Gamma)^0$.
\end{lemma}
\begin{proof}
Put $H\coloneqq \Aut^0(Z,A,Y)$, acting effectively on $Z$, and let
$K\coloneqq \operatorname{Stab}_H(\omega_0)$. Theorem~\ref{un:thm:fixed} places $\omega_1$ in the metric
orbit $H/K$. As established in its proof, $H$ is complex
reductive and $K$ is a maximal compact subgroup.
Choose a faithful matrix realization of $H$ in which $K$ is
unitary and $H$ is stable under conjugate transpose.
The global Cartan decomposition
\cite[Proposition~7.14]{KnappLie}
identifies $\sqrt{-1}\Lie K$ with $H/K$.
Thus there is a unique $X\in\Lie K$ such that, for the
pushforward action on forms,
\[
\omega_1=(\exp(\sqrt{-1}X))_*\omega_0.
\]
Since $\Gamma$ preserves $\omega_0$ and $Y$, it normalizes $H$ and $K$.
Invariance of $\omega_1$ and uniqueness of these polar coordinates imply
$\operatorname{Ad}_\gamma X=X$ for every $\gamma\in\Gamma$.
Thus $a\coloneqq \exp(\sqrt{-1}X)$ centralizes $\Gamma$, is joined to the identity
in that centralizer, and satisfies $a^*\omega_1=\omega_0$.
\end{proof}

\subsection{Finite-group equivariance in the real class}
\label{eqv:sec:real-proof}

We now prove the equivariant real-class theorem. The auxiliary classes and
their linearizations can be chosen invariantly, so the argument
takes place in the representation centralizer while preserving
the exact class. The resulting statement is as follows.

\begin{theorem}[Equivariant convergence in a real class]
\label{eqv:thm:real}
Assume the cone and initial-energy hypotheses of
Theorem~\ref{ru:thm:unique}. Let $\Gamma$ be a finite group of
holomorphic automorphisms of $M$ preserving $\omega_0$.
For this fixed $\Gamma$-invariant flow, the convergence and stability
conclusions of Theorem~\ref{ru:thm:unique} have the following
equivariant refinements.
\begin{enumerate}[label=\textnormal{(\arabic*)}]
\item All subsequential limits of this flow are identified by
$\Gamma$-equivariant holomorphic isometries, and the whole flow
converges geometrically modulo $\Gamma$-equivariant holomorphic
isometry, in the sense of
Definition~\ref{app:def:geometric-convergence} with equivariant
identifications.
\item The class-preserving configuration of
Proposition~\ref{rd:prop:degeneration} can be chosen
$\Gamma$-equivariantly. Its total space is normal, its central
fibre is normal, reduced and irreducible, and its invariant is
negative if the limit is not cscK and zero if it is cscK.
\item Suppose every $\Gamma$-equivariant normal relatively
K\"ahler cohomological test configuration with normal, reduced,
irreducible central fibre has nonnegative Donaldson--Futaki
invariant. Then the limiting orbifold is cscK.
If, in addition, equality occurs only for products, the flow
converges geometrically to a smooth cscK metric on $M$ in
$2\pi\Omega$. The identification with the prescribed
$\Gamma$-action is equivariant if the zero-invariant products
are equivariant products. This holds in particular if
$\Aut^0(M)=\{1\}$.
\end{enumerate}
The configuration in (2) carries finitely many linearized
relatively ample bundles, each inducing complete embeddings on
the general and central fibres, whose positive real combination
is its cohomological class. It suffices in (3) to test these
configurations. 
\end{theorem}

\begin{proof}[Proof of Theorem~\ref{eqv:thm:real}]
\textit{Invariant auxiliary classes and linearizations.}
Let $W=W_\Omega$ be the minimal rational subspace containing
$\Omega$. For each $\gamma\in\Gamma$, the subspace
$\ker(\gamma^*-I)$ is rational and contains $\Omega$.
Consequently $W\subset H^2(M,\R)^\Gamma$, pointwise.
All the auxiliary classes in Lemma~\ref{ru:lem:embeddings}
are therefore invariant. Since $M$ is rational,
$\operatorname{Pic}^0(M)=0$, so the corresponding holomorphic
line bundles are preserved as well. The norm-bundle construction in the proof of
\cite[Lemma~5.1]{Polarized} gives genuine
$\Gamma$-linearizations on their $|\Gamma|$-th powers. Taking these powers divides the
positive coefficients in \eqref{rd:eq:positive-basis} by
$|\Gamma|$. It does not change $\Omega$ or rescale the flow:
these bundles are auxiliary polarizations.

The harmonic curvatures are invariant along the invariant flow.
Averaging logarithms gives invariant Hermitian metrics with
these curvatures. Equivariant compactness and complete-section
convergence retain the genuine representations on every limit,
by the proof of \cite[Corollary~5.3]{Polarized}.

\textit{Dynamics in the representation centralizer.}
By Maschke's theorem and Schur's lemma
\cite[Theorem~4.1.1, Proposition~2.3.9 and Corollary~2.3.10]{EtingofRep},
write
\[
 V_\nu^*=
 \bigoplus_\rho W_\rho\otimes\mathbb C^{m_{\nu\rho}},
\]
where the $W_\rho$ are pairwise nonisomorphic irreducible
complex $\Gamma$-representations and $m_{\nu\rho}$ are their
multiplicities. Under this decomposition,
\[
 \mathcal R\coloneqq
 \prod_\nu\operatorname{Cent}_{GL(V_\nu^*)}(\Gamma)
 \cong\prod_{\nu,\rho}GL(m_{\nu\rho},\mathbb C).
\]
All Gram matrices and positive frame updates belong to
$\mathcal R$. Following the multiplicity-space formulation of
\cite[Lemma~5.4]{Polarized}, apply
Lemma~\ref{ru:lem:spectrum}, recording the ordered spectrum
in each multiplicity space separately.
Its countable-field argument and connectedness give one
$\Lambda\in\Lie\mathcal R$, up to the compact centralizer.
Normalize one trace for each auxiliary factor $V_\nu^*$.
Subtracting separate traces in its isotypic blocks would
incorrectly alter their relative projective weights.
The blockwise argument in \cite[Lemma~5.5]{Polarized} now gives
\[
 x=C_0[M],\quad C_0\in\mathcal R,\qquad
 y=\lim_{t\to\infty}e^{t\Lambda}x,
 \qquad [Z]\in\overline{\mathcal R_\Lambda y}.
\]
The target stabilizer in $\mathcal R_\Lambda$ is reductive.
Indeed, $\Gamma$ acts by limiting isometries and preserves
the compact real form of the extremal centralizer. Taking
its fixed Lie algebra preserves the complexification property.
The linear lifts add a scalar torus. The possible character
obstructions are finite and do not change the identity component.

\textit{Equivariant degeneration and uniqueness.}
The torus $\overline{e^{\sqrt{-1}\R\Lambda}}$ used in
Proposition~\ref{rd:prop:degeneration} is contained in
$\mathcal R$, and its complexification is central in
$\mathcal R_\Lambda$. The rational replacement of the generator
therefore stays in this group. Accessibility, the combined
subgroup and the initial coordinate change also belong to
$\mathcal R$. Their universal Hilbert family is
$\Gamma$-equivariant, with all factor bundles linearized.
Its positive real class and invariant are exactly those in
Steps~4--5 of that proposition. This proves (2).

For (1), use Lemma~\ref{eqv:lem:fixed} in the common real
class. Its correcting automorphism centralizes $\Gamma$ and
the extremal field and preserves every auxiliary bundle,
since $\operatorname{Pic}^0(Z)=0$. The correcting path lifts
commuting with $\Gamma$: a scalar commutator is a character
of the finite group and is constant along the path.
Harmonic uniqueness makes the lifted isometry unitary in
every auxiliary section space. Thus the orbit-intersection
argument in the fixed-flow proof of Theorem~\ref{ru:thm:unique}(1) works in
$\mathcal R_\Lambda$, and gives equivariant uniqueness,
including the translated time windows.

Finally apply the sign and product argument of
Theorem~\ref{rd:thm:convergence}(2)--(3) to (2). A product identifies
the exact real class and removes the quotient singularities.
An equivariant product preserves the group identification.
If $\Aut^0(M)=\{1\}$, its $\C^*$-action on $M$ is trivial,
and the commuting $\Gamma$-action is constant over the base,
as in \cite[Corollary~5.6]{Polarized}. This proves (3).
\end{proof}

\subsection{Equivariant auxiliary embeddings for cscK sequences}

The continuation argument in Section~\ref{real:sec:existence}
also uses the static counterpart of the preceding constructions. We record it separately so that
the varying-class results of Subsection~\ref{vc:sec:compactness}
can first be read without a group action.

\begin{corollary}[Finite symmetry in the static construction]
\label{eqv:cor:static}
In Theorem~\ref{vc:thm:compactness}, suppose a fixed finite group
$G$ acts holomorphically on $M$ and preserves every $\omega_j$.
Then the limiting cscK orbifold has a holomorphic isometric
$G$-action, and the convergence can be extracted equivariantly.
The auxiliary complete embeddings of
Proposition~\ref{vc:prop:polarized-limit} can be chosen
$G$-equivariantly after taking a fixed finite tensor power of
each source bundle. Their Hilbert limits lie in the orbit
closures for the centralizers of the resulting representations.
\end{corollary}
\begin{proof}
We use the notation of Theorem~\ref{vc:thm:compactness},
Lemma~\ref{vc:lem:marking} and the proof of
Proposition~\ref{vc:prop:polarized-limit}.
We adapt the proof of \cite[Theorem~6.12]{Compactness} to the
static convergence of Theorem~\ref{vc:thm:compactness}.
Here smooth convergence is only asserted on $Z\setminus\Sigma$,
and $\Sigma$ may contain smooth points of $Z$.
First use the global GH correspondences to extract the actions.
Since $G$ is finite, one subsequence gives, for every $\gamma\in G$,
an isometry $\gamma_\infty:Z\to Z$.
Indeed, the transported graphs of the source isometries converge
to graphs of distance-preserving maps. Following $p_j$,
$\eta p_j$ and $\gamma\eta p_j$ through the same GH
correspondences gives
\[
 (\gamma\eta)_\infty=\gamma_\infty\eta_\infty.
\]
The identity and inverse relations pass to the limit in the same
way, so these maps define an action on $Z$.

Replace $\Sigma$ by the finite invariant set
$\bigcup_{\gamma\in G}\gamma_\infty(\Sigma)$.
Each limiting isometry preserves the regular locus, since
regular points are characterized by Euclidean tangent cones.
For a compact subset of $Z\setminus\Sigma$, choose a larger
relatively compact comparison domain in $Z\setminus\Sigma$
containing both this subset and its image under $\gamma_\infty$.
The interior coverage argument in that proof applies on these domains.
It uses compatibility with GH convergence and smooth convergence
on a slightly larger domain disjoint from $\Sigma$.
Thus the local maps $\Phi_j^{-1}\gamma\Phi_j$ are well-defined
near the compact subset for large $j$.
The isometry identity bounds their first derivatives, and
preservation of the Levi-Civita connection gives bounds for
all higher derivatives. Their smooth limits are the maps
$\gamma_\infty$. Passing to the limit in complex linearity
shows that the action is holomorphic on $Z\setminus\Sigma$.

At a smooth point belonging to $\Sigma$, each $\gamma_\infty$
is already smooth as an isometry between Riemannian manifolds.
Its differential is complex linear on the dense complement of
$\Sigma$, hence everywhere by continuity.
At an orbifold point, use the lifting of orbifold isometries
employed in the proof of \cite[Theorem~6.12]{Compactness}.
The resulting smooth lifts to uniformizing charts are holomorphic
away from their centres. Complex linearity extends across the
centres by continuity. Thus the limiting action is by orbifold
biholomorphisms.

Finally, take the union of the translates of each GH
correspondence under the source and limiting $G$-actions.
As in the same proof, these
correspondences are exactly equivariant and their distortions
still tend to zero. They retain the original comparison
embeddings. This gives equivariant GH convergence and smooth
convergence of the local action maps on compact subsets of
$Z\setminus\Sigma$.

For the equivariant assertion, the minimal rational subspace $W$
is contained in $H^2(M,\mathbb R)^G$, since $\Omega_0$ is invariant
and this fixed subspace is rational. The bundles $\mathcal B_i$
are therefore $G$-invariant. The norm-bundle construction used in the proof of
Theorem~\ref{eqv:thm:real} gives genuine $G$-linearizations
on the fixed powers $\mathcal B_i^{|G|}$.
Choose $G$-invariant Hermitian metrics with the harmonic curvatures,
and extract the connections also on these linearized powers.
Their underlying limiting bundles are the corresponding tensor
powers of the bundles already constructed. This only introduces a
further fixed power of $L_B$ before choosing the final exponent $m$.
It changes none of the marked real classes. All these extractions
retain the finite group action. In each section
space the varying inner products are $G$-invariant, so
representation-compatible orthonormal bases differ by elements of
the centralizer of $G$. The Hilbert limit is in that orbit closure.
\end{proof}

\section{Closedness and existence}\label{real:sec:existence}

We first prove closedness and continuation under cohomological
K-polystability of the real classes themselves. The essential input
is a degeneration preserving the exact limiting class. We then give a separate
criterion using a local restricted gap in nearby rational auxiliary
classes. The next section verifies that criterion for the family $A_a$.

\subsection{Closedness and continuation under cohomological K-polystability}
\label{ex:sec:closedness}

Let $M$ be a compact K\"ahler surface with $\Aut^0(M)=\{1\}$,
let $G$ be a finite group of holomorphic automorphisms, and put
$c\coloneqq c_1(M)$. Use the controlled cone
$\mathcal U$ from \eqref{vc:eq:cone}, with metric class
$[\omega]=2\pi\Omega$. If $\mathcal U$ is nonempty, $M$ is
rational and projective, as explained in the introduction.
The finite group $G$ may be trivial.
For a real class $\Omega$, write $W_\Omega$ for the smallest real
linear subspace of $H^2(M,\mathbb R)$ defined over $\mathbb Q$
and containing $\Omega$.
For a smooth K\"ahler metric $\omega_0\in2\pi\Omega$, the
initial-energy condition used below is
\begin{equation}\label{real:eq:energy}
 \Ca(\omega_0)<8\pi^2\delta(\Omega),
 \qquad\text{equivalently}\qquad
 \widehat{\Ca}(\omega_0)<\frac{4\delta(\Omega)}{\Omega^2}.
\end{equation}
Here $\widehat{\Ca}\coloneqq\Ca/\Vol$ and
$\Vol(M,\omega_0)=2\pi^2\Omega^2$.

\begin{theorem}[Closedness and continuation under cohomological K-polystability]
\label{real:thm:poly-continuation}
In the preceding setting the following statements hold.
\begin{enumerate}[label=\textnormal{(\arabic*)}]
\item \textbf{Closedness at a polystable real class.}
Let $\Omega_j,\Omega_0\in\mathcal U\cap H^2(M,\mathbb R)^G$
with $\Omega_j\to\Omega_0$. Suppose that each $2\pi\Omega_j$
contains a smooth $G$-invariant cscK metric and that
$(M,\Omega_0)$ is cohomologically K-polystable. Then
$2\pi\Omega_0$ contains a smooth $G$-invariant cscK metric.
\item \textbf{Continuation from one class.}
Let $D\subset\mathcal U\cap H^2(M,\mathbb R)^G$ be nonempty
and connected, and suppose that $(M,\Omega)$ is cohomologically
K-polystable for every $\Omega\in D$. If one class $2\pi\Omega$
with $\Omega\in D$ contains a smooth cscK metric, then every
class $2\pi B$, $B\in D$, does. It suffices instead to have one
smooth $G$-invariant initial metric satisfying
\eqref{real:eq:energy} in a class in $D$.
\end{enumerate}
The resulting cscK metrics are unique and $G$-invariant, and depend
locally smoothly on the real class. For a set $D$ that is not a
manifold, this means that they locally extend to smooth families
on open subsets of the invariant real cohomology space.
In both assertions, it suffices to impose cohomological
K-polystability only on normal relatively K\"ahler
$G$-equivariant configurations with normal, reduced, irreducible
central fibre.
\end{theorem}

In part~(1), stability is assumed at the limiting class $\Omega_0$.
It is not inferred by taking a limit of stable classes.
We first construct the
exact-class degeneration needed in the proof.

\begin{proposition}[An exact-class degeneration of a static cscK limit]
\label{vc:prop:real-degeneration}
In the setting of Theorem~\ref{vc:thm:compactness}, suppose
$\Omega_j\to\Omega_0$ and the cscK metrics are $G$-invariant.
Their marked limit $(Z,\omega_\infty,\Omega_\infty)$ is the
central fibre of a $G$-equivariant normal relatively K\"ahler
cohomological test configuration for $(M,\Omega_0)$ with
restricted class $\Omega_\infty$ and Donaldson--Futaki
invariant zero. It is induced by an algebraic one-parameter
subgroup in a product of complete auxiliary embeddings. The
factor bundles $\mathcal L_\nu$ are genuinely $G$-linearized
and relatively ample, induce complete embeddings on both the
general and central fibres, and satisfy
\[
 \mathcal A=\sum_\nu\alpha_\nu c_1(\mathcal L_\nu),
 \qquad \alpha_\nu>0,\qquad
 c_1(\mathcal L_\nu|_{\mathcal X_1})\in W_{\Omega_0}.
\]
\end{proposition}

\begin{proof}
Theorem~\ref{vc:thm:compactness}, Lemma~\ref{vc:lem:marking}
and Corollary~\ref{eqv:cor:static} give the equivariant cscK
orbifold and its marked real class. Choose finitely many rational classes in
$W_{\Omega_0}$ near $\Omega_0$ whose positive cone contains
$\Omega_0$, using only the simplex construction in
Lemma~\ref{ru:lem:embeddings}.
They are ample on both $M$ and this fixed limit, by
Proposition~\ref{vc:prop:polarized-limit}. Apply that proposition
and its equivariant refinement, Corollary~\ref{eqv:cor:static},
to each class and take one common subsequence over this finite
list. After the indicated powers and genuine $G$-linearizations,
we obtain complete embeddings by $A_\nu$ and $A_{\nu,Z}$, with
\[
 \Omega_0=\sum_\nu\alpha_\nu c_1(A_\nu),\qquad
 \Omega_\infty=\sum_\nu\alpha_\nu c_1(A_{\nu,Z}),\qquad
 \alpha_\nu>0.
\]
The mixed intersection identities \eqref{vc:eq:intersection}
fix the Segre Hilbert polynomial. The product Hilbert limit is
exactly $Z$ as a scheme, by the same complete-section and
Hilbert-polynomial argument. All factors use the same
geometric marking. In particular, with
$\mathcal R=\prod_\nu\operatorname{Cent}_{GL(H^0(M,A_\nu)^*)}(G)$,
we have $[Z]\in\overline{\mathcal R[M]}$.

The cscK stabilizer in $\mathcal R$ is reductive.
Its geometric identity component is the centralizer of $G$
in the complexification of the compact isometry group. The
vanishing $H^1(Z,\mathcal O_Z)=0$ is supplied by
Proposition~\ref{vc:prop:polarized-limit}. Preservation of every
auxiliary bundle follows from $\operatorname{Pic}^0(Z)=0$.
The linear stabilizer adds a scalar torus, and the finite
character obstruction to commuting lifts does not affect
reductivity, as in Theorem~\ref{eqv:thm:real}.
Accessibility gives $g\in\mathcal R$ and an algebraic
one-parameter subgroup $\mu\subset\mathcal R$ with
$\lim_{z\to0}\mu(z)g[M]=[Z]$. If the target is already in the
orbit, take $\mu$ trivial.

Pull back the universal Hilbert family. The factor bundles
extend over it, and their positive combination defines
$\mathcal A$ with the two stated class restrictions.
Flatness and normal fibres give normal total space. Fibrewise
ampleness and compactification give relative K\"ahlerness.
These are exactly Step~4 of Proposition~\ref{rd:prop:degeneration}.
Its fixed-family Riemann--Roch argument identifies
$\DF(\mathcal X,\mathcal A)$ with the Futaki character of the
central action in $\Omega_\infty$. Since $\omega_\infty$ is
cscK, this character is zero. 
\end{proof}

\begin{proof}[Proof of Theorem~\ref{real:thm:poly-continuation}]
\textit{Step 1. Closedness at the limiting class.}
After discarding finitely many terms, the classes $\Omega_j$ lie
in a compact subset of the controlled cone.
Theorem~\ref{vc:thm:compactness} gives a marked cscK orbifold limit
$(Z,\omega_\infty)$. Proposition~\ref{vc:prop:real-degeneration}
realizes it as the central fibre of a normal relatively K\"ahler
$G$-equivariant configuration for the exact pair $(M,\Omega_0)$,
with $\DF=0$. Cohomological K-polystability of this pair makes
the configuration a product.

Normalize the product identification at $t=1$. Homotopy of the
fibre inclusions then identifies the central class with
$\Omega_0$, as in Theorem~\ref{rd:thm:convergence}. The smooth
coarse central surface has trivial local orbifold groups by
\cite[Corollary~4.10]{Polarized}. Its cscK metric pulls back to a
smooth metric in $2\pi\Omega_0$. Since $\Aut^0(M)=\{1\}$,
the product identification is $G$-equivariant, as in
Theorem~\ref{eqv:thm:real}. This proves part~(1).

\textit{Step 2. Openness and smooth dependence.}
Fixed-class uniqueness and $\Aut^0(M)=\{1\}$ imply that a cscK
metric in a given class is unique, and hence invariant under every
element of $G$ preserving that class. Choose a linear map $T$
from $H^2(M,\mathbb R)^G$ to smooth $G$-invariant closed real
$(1,1)$-forms representing the corresponding classes. Choose
representatives of a basis and average over $G$.
At a known cscK metric $\omega_\Omega$, put
\[
 \rho(B)\coloneqq\omega_\Omega+2\pi T(B-\Omega).
\]
For $B$ near $\Omega$ these are positive forms in $2\pi B$.
Let $\Pi_\Omega$ project onto
functions of zero $\omega_\Omega^2$-mean and define
\[
 \mathscr F(B,\varphi)\coloneqq
 \Pi_\Omega S_{\rho(B)+\ddc\varphi}.
\]
For $k\geq0$ and $0<\alpha<1$, this is a smooth map from an
open set in $H^2(M,\R)^G\times C^{k+4,\alpha}_0(M)$ to
$C^{k,\alpha}_0(M)$. Its derivative in $\varphi$ at
$(\Omega,0)$ is the negative Lichnerowicz operator. This
elliptic self-adjoint operator has only constants in its kernel
because $H^0(M,T^{1,0}M)=0$. It is therefore an isomorphism
between the two zero-mean spaces.

The implicit-function theorem and elliptic regularity give a
smooth local family of smooth cscK metrics for real classes near
$\Omega$. Fixed-class uniqueness makes the family invariant and
makes all such local families agree wherever they overlap. This
proves the asserted local smooth dependence. For a connected set
$D$ that is not a manifold, the statement means that the metrics
are locally restrictions of these smooth families on open subsets
of the invariant real cohomology space.

\textit{Step 3. Continuation.}
Let $\mathcal E\subset D$ consist of classes containing a smooth
cscK metric. It is relatively open by Step~2 and relatively closed
by Step~1, since each limiting class in $D$ is assumed
cohomologically K-polystable. It is nonempty by hypothesis, so
connectedness gives $\mathcal E=D$. If the starting point is
instead an invariant metric satisfying \eqref{real:eq:energy},
Theorem~\ref{eqv:thm:real}(3) supplies a cscK metric in that
class and starts the same argument.

The configurations in Step~1 and in this initial flow all have
normal, reduced, irreducible central fibre. Thus the stated
restricted cohomological polystability also suffices. Taking $G$
trivial proves the cohomological K-polystability assertion in
part~(4) of Theorem~\ref{ru:thm:unique}.
\end{proof}

\begin{remark}\label{vc:rem:scope}
Closedness uses compactness together with polystability of the
exact limiting real class. It is not a consequence of the
controlled-cone inequalities alone. The proof is static: it applies
compactness to cscK metrics already present in the approximating
classes and does not interchange a varying-class limit with an
infinite-time limit of Calabi flows. 
\end{remark}

\subsection{A local restricted gap criterion for the example}
\label{real:sec:gap-criterion}

The preceding theorem assumes stability of the real class itself.
For the example, we instead verify a quantitative estimate in
nearby rational auxiliary classes. On each fixed constructed
family, this estimate passes to real coefficients and rules out
a nonproduct degeneration of zero Donaldson--Futaki invariant.
Strict positivity at nearby rational coefficients alone would
only give nonnegativity at the limit. The gap below supplies the
strict inequality needed in that case. It is not an additional
requirement in Theorem~\ref{real:thm:poly-continuation}.

We use the non-Archimedean $J$-functional with the actual
general-fibre polarization $L$. If $\mathrm{DH}$ denotes the
probability measure obtained as the limit of the section-weight
distributions of $\mathcal L_0^k$, with each weight divided by $k$, then
\[
 J^{\NA}(\mathcal X,\mathcal L)\coloneqq
 \sup\operatorname{supp}\mathrm{DH}
       -\int_{\R}\lambda\,d\mathrm{DH}(\lambda).
\]
This is \cite[Definition~3.10 and Proposition~7.8]{BHJ}.
In particular, tensoring $\mathcal L$ by a positive integer $q$
multiplies $J^{\NA}$ by $q$. We do not divide this norm by
the polarization exponent.

We say that a $G$-invariant class $\Omega\in\mathcal U$ satisfies the
\emph{local restricted stability gap} if there are a neighborhood
$U_\Omega$ of $\Omega$ in $W_\Omega\cap\mathcal K_M$ and
$\gamma_\Omega>0$ such that the following holds. For every rational
$B\in U_\Omega$, every $G$-linearized integral polarization
$L$ with $c_1(L)=mB$, and every nonproduct normal ample
$G$-equivariant test configuration with actual general-fibre
polarization $L$ and reduced, irreducible central fibre,
\begin{equation}\label{real:eq:local-gap}
 J^{\NA}>0,\qquad
 \frac{\DF(\mathcal X,\mathcal L)/L^2}
      {J^{\NA}(\mathcal X,\mathcal L)}
 \geq\frac{\gamma_\Omega}{2m}.
\end{equation}
The positive integer $m$ includes all denominator-clearing,
linearization and test-configuration powers. Both $U_\Omega$ and
$\gamma_\Omega$ may depend on $\Omega$. Denote by
$\mathcal S_G\subset\mathcal U\cap H^2(M,\R)^G$ the set of
classes satisfying this condition. The rational approximations in
this definition need only lie in $W_\Omega$. The class $\Omega$
itself need not be rational.

The local gap is an additional hypothesis beyond the controlled-cone
inequalities. Its use below tests only the constructed projective
families, and does not assert K-polystability against all
cohomological test configurations.

\begin{lemma}[The rational gap on a fixed real-class family]
\label{eqv:lem:real-gap}
Let $\Omega\in\mathcal S_G$. Suppose a normal flat
$G$-equivariant degeneration is obtained by an algebraic
one-parameter subgroup in a product of complete embeddings by
genuinely $G$-linearized bundles $A_\nu$ on $M$, with
$c_1(A_\nu)\in W_\Omega$. Suppose its central fibre $Z$ is a
normal integral orbifold surface, and every factor still gives
a complete embedding on $Z$. Let $\mathcal L_\nu$ be the
restricted factor bundles on the family. If
\[
 \alpha_\nu>0,\qquad
 \Omega=\sum_\nu\alpha_\nu c_1(A_\nu),\qquad
 \mathcal A=\sum_\nu\alpha_\nu c_1(\mathcal L_\nu),
\]
then $\DF(\mathcal X,\mathcal A)>0$ unless the family is a
product. A product has invariant zero.
\end{lemma}

\begin{proof}
Assume the underlying family is not a product. Write
$B_b=\sum_\nu b_\nu c_1(A_\nu)$ and
$\mathcal A_b=\sum_\nu b_\nu c_1(\mathcal L_\nu)$.
For positive rational vectors $b$ sufficiently close to $\alpha$,
$B_b$ belongs to the neighborhood $U_\Omega$ in the local gap.
Choose an integer $q>0$ with all $qb_\nu$ integral.
The bundle $\bigotimes_\nu\mathcal L_\nu^{qb_\nu}$ is genuinely
$G$-linearized and relatively ample, with general class $qB_b$.
The central fibre and nonproduct underlying family are unchanged.

Here is the norm needed to pass the gap to the real coefficients.
Choose circle-invariant Hermitian inner products in the fixed
complete factors, and use their Fubini--Study metrics on $Z$.
Let $\omega_\nu$ be their curvature forms and $h_\nu$ their
Hamiltonians in the normalization of section weights. For $b>0$ put
\[
 \omega_b=\sum_\nu b_\nu\omega_\nu,\qquad
 h_b=\sum_\nu b_\nu h_\nu,\qquad
 J(b)\coloneqq\max_Zh_b-
       \frac{\int_Zh_b\,\omega_b^2}{\int_Z\omega_b^2}.
\]
These forms and functions extend smoothly in orbifold charts.
The forms can be degenerate at isolated quotient points, but
$\omega_b$ is positive on $Z_{\reg}$ and its volume measure has
full support on $Z$. Equivariant orbifold Riemann--Roch identifies
the central weight distribution with the pushforward of the
normalized measure $\omega_b^2/\int_Z\omega_b^2$ by $h_b$.
Indeed its moments are the leading coefficients of the successive
torus derivatives of the index. The isolated nonidentity sectors
have lower degree, by the orbifold index calculation of
\cite[formula~(IV), p.~153]{Kawasaki} and
\cite{Meinrenken}.
Equality of all moments identifies these compactly supported
measures. Full support of the volume measure also identifies the
supremum of the pushforward's support with $\max_Zh_b$.
We therefore obtain
\[
 J^{\NA}(\mathcal X,q\mathcal A_b)=qJ(b)
\]
for rational $b$, by \cite[Definition~3.10 and Proposition~7.8]{BHJ}.
The displayed formula defines a continuous, degree-one homogeneous
function $J$ also for real $b>0$.

Our surface invariant is homogeneous of degree two. Thus
\eqref{real:eq:local-gap}, applied with full exponent $q$, gives
\[
 \frac{q^2\DF(\mathcal X,\mathcal A_b)}
      {q^2B_b^2\,qJ(b)}
       \geq\frac{\gamma_\Omega}{2q},
 \qquad\text{hence}\qquad
 \frac{\DF(\mathcal X,\mathcal A_b)}{B_b^2}
       \geq\frac{\gamma_\Omega}{2}J(b).
\]
The fixed-family continuity of the intersection formula
\eqref{rd:eq:cohom-df} lets $b\to\alpha$ and yields
\begin{equation}\label{eqv:eq:real-gap}
 \frac{\DF(\mathcal X,\mathcal A)}{\Omega^2}
       \geq\frac{\gamma_\Omega}{2}J(\alpha).
\end{equation}

If $J(\alpha)=0$, full support makes $h_\alpha$ constant.
The Hamiltonian vector field vanishes on $Z_{\reg}$, so the
central $\C^*$-action fixes every point of $Z$. In each complete
factor embedding the irreducible image lies in a single
projective weight space. Nondegeneracy makes that weight space
the entire factor, so the one-parameter subgroup is scalar in
every factor. Its projective orbit is constant and the universal
Hilbert family is a product, a contradiction.
Consequently $J(\alpha)>0$, proving strict positivity.
If the underlying family is a product, $\Aut^0(M)=\{1\}$
makes its induced action on $M$ trivial. Its invariant is zero.
\end{proof}

This lemma applies to the projective families with auxiliary
bundles constructed here. It does not claim that the rational
gap tests every analytic cohomological test configuration.

\begin{theorem}[Closedness under the local restricted gap]
\label{vc:thm:closedness}
In the setting of Section~\ref{real:sec:existence}, let
$\Omega_0\in\mathcal S_G$. If $\Omega_j\to\Omega_0$ and each
class $2\pi\Omega_j$ contains a smooth $G$-invariant cscK metric,
then $2\pi\Omega_0$ contains a smooth $G$-invariant cscK metric.
\end{theorem}
\begin{proof}
Proposition~\ref{vc:prop:real-degeneration} supplies a configuration
in the exact class $\Omega_0$ with $\DF=0$. Its auxiliary bundles
satisfy the hypotheses of Lemma~\ref{eqv:lem:real-gap}, so that
lemma forces the family to be a product. The product identification
then gives a smooth $G$-invariant cscK metric in $2\pi\Omega_0$,
exactly as in the closedness proof of
Theorem~\ref{real:thm:poly-continuation}. This uses the gap only
on the constructed family and does not assert full cohomological
K-polystability of $(M,\Omega_0)$.
\end{proof}

\begin{theorem}[Existence and continuation under the local restricted gap]
\label{real:thm:existence}
With the preceding notation, the following statements hold.
\begin{enumerate}[label=\textnormal{(\arabic*)}]
\item \textbf{A criterion in one real class.}
For every $\Omega\in\mathcal S_G$, the class $2\pi\Omega$
contains a smooth cscK metric if and only if it contains a smooth
$G$-invariant K\"ahler metric $\omega_0$ satisfying
\eqref{real:eq:energy}.
\item \textbf{Continuation from one class.}
Let $D\subset\mathcal S_G$ be a nonempty connected subset.
If one class $2\pi\Omega$, $\Omega\in D$, admits a smooth
$G$-invariant K\"ahler metric satisfying
\eqref{real:eq:energy}, then every class $2\pi\Omega$, $\Omega\in D$,
contains a smooth cscK metric. In particular, one known cscK
metric in $D$ is a sufficient starting point.
\end{enumerate}
Whenever these metrics exist, they are unique and $G$-invariant.
They are locally smooth functions of the real K\"ahler class.
In particular, their restrictions to smooth paths in $D$ are smooth.
No rationality assumption on the classes in $D$, and no density
assumption on its rational points, is required.
\end{theorem}

\begin{proof}
\textit{Step 1. Existence in one real class.}
Fix $\Omega\in\mathcal S_G$ and an invariant initial metric
satisfying \eqref{real:eq:energy}. The long-time theorem
\cite[Corollary~1.3]{LongTime} gives a global smooth invariant
Calabi flow with a uniform Sobolev bound.
Theorem~\ref{eqv:thm:real}(2) produces a configuration in its
exact class with normal integral central fibre and $\DF\leq0$.
Its auxiliary bundles satisfy Lemma~\ref{eqv:lem:real-gap}, which
forces the configuration to be a product. Products have invariant
zero because $\Aut^0(M)=\{1\}$. The strict negative sign for a
non-cscK limit therefore excludes that case. The product
identification gives a smooth invariant cscK metric in
$2\pi\Omega$, and the whole flow converges to its equivariant
holomorphic-isometry class. Conversely, a cscK metric has zero
Calabi energy and is $G$-invariant by uniqueness, proving (1).

\textit{Step 2. Continuation.}
Part~(1) supplies one cscK class in $D$. The cscK locus is
relatively open by the implicit-function argument in the proof
of Theorem~\ref{real:thm:poly-continuation}, and relatively closed
by Theorem~\ref{vc:thm:closedness}. Connectedness proves (2).
The same implicit-function argument and fixed-class uniqueness
give the asserted smooth dependence. This reuses the openness
argument, not an implication from the gap to full cohomological
K-polystability. This proves the alternative continuation criterion
described after Theorem~\ref{ru:thm:unique}.

Only the coefficients of the auxiliary bundles on each fixed
family are approximated by rational numbers. These rational
classes lie in $W_\Omega$ near $\Omega$ and need not lie in $D$.
The constants $\gamma_\Omega$ and neighborhoods $U_\Omega$ need
not be uniform over $D$.

For the initial-energy condition near a known cscK metric, use
the forms $\rho(B)$ from the proof of
Theorem~\ref{real:thm:poly-continuation}. Smooth dependence of
scalar curvature gives, after shrinking the neighborhood,
\[
 \widehat{\Ca}(\rho(B))\leq C_\Omega\|B-\Omega\|^2,
 \qquad
 \frac{4\delta(B)}{B^2}>
 \frac12\frac{4\delta(\Omega)}{\Omega^2}>0.
\]
Here any fixed norm on the finite-dimensional cohomology space
may be used. Thus all sufficiently close classes have invariant
initial metrics satisfying \eqref{real:eq:energy}.
\end{proof}

\section{Continuation on a symmetric four-point blow-up}
\label{cone:sec:fourpoint}

We apply the local restricted gap criterion of
Subsection~\ref{real:sec:gap-criterion} to a one-parameter family
on a fixed rational surface. The starting point is an anticanonical
K\"ahler--Einstein metric. A restricted equivariant stability estimate
then prevents the cscK metrics from leaving the original complex
surface at a limiting parameter. This gives existence throughout
an explicit real interval of symmetric classes by continuation
from the anticanonical class.

\subsection{Geometry and the initial-energy condition}

Let $X$ be the blow-up of $\mathbb P^2$ at four points in general
position. Write $H$ for the pullback of the line class and
$E\coloneqq E_1+\cdots+E_4$. The projective transformations permuting
the four points lift to a group $G\cong S_4$ on $X$.
After a projective change of coordinates we may take the points to be
\[
 [1:1:1],\quad [1:-1:-1],\quad [-1:1:-1],\quad [-1:-1:1].
\]
In these coordinates, $G$ is generated by coordinate permutations
and diagonal sign changes with product $+1$. Its action lifts to
$X$ and preserves the relevant line bundles. Put
\begin{equation}\label{cone:eq:symmetric-classes}
 c\coloneqq c_1(X)=3H-E,\qquad A_a\coloneqq aH-E,
 \qquad J\coloneqq[3,a_+),\quad a_+\coloneqq8+\frac{10}{\sqrt3}.
\end{equation}
These are real classes. No rationality of $a$ is imposed.
The conic pencil through the four points gives the nef class
$Q\coloneqq2H-E$. For $a>2$, the class
$A_a=(a-2)H+Q$ has positive square and positive intersection with
every curve: a curve of positive $H$-degree is covered by the first
term, and the only curves contracted by $H$ are the $E_i$, with
$Q\cdot E_i=1$. Thus $A_a$ is ample. Also
$\chi(\mathcal O_X)=1$ and $\Aut^0(X)=\{1\}$. Indeed, the identity
component fixes the finitely many $(-1)$-curves and descends after
contracting the $E_i$ to a projective transformation fixing the
four-point projective frame.

The existence of cscK metrics on this surface is classical.
Tian--Yau's theorem \cite[Theorem~3.3(ii)]{TianYauKE} gives the
anticanonical K\"ahler--Einstein metric in $2\pi A_3$.
Since $\Aut^0(X)=\{1\}$, local openness also gives cscK metrics
in nearby classes. Arezzo--Pacard's gluing construction gives
cscK metrics on the four-point blow-up with sufficiently small
exceptional divisors. See \cite[Theorem~1.4, Corollary~1.2 and
Section~7.4]{APBlowupII}.
In the present family, division by $a$ gives
$a^{-1}A_a=H-a^{-1}E$, so small equal exceptional coefficients
correspond to large $a$. The smallness threshold in those gluing
statements is not numerical. Those statements alone therefore do
not verify the fixed coefficient $1/12$ or the entire interval $J$.

The contribution here is continuation from the anticanonical
metric throughout the explicit interval $J$, including irrational
classes, using compactness and the restricted equivariant gap.
Theorem~\ref{cone:thm:fourpoint} also gives smooth dependence on
$a$ and convergence of the admissible invariant Calabi flows.
The class $A_{12}=12H-E$ serves as a concrete illustration, with
an explicit initial-energy threshold computed below.

The relevant intersection numbers and energy threshold are
\begin{align}
 A_a^2&=a^2-4,&c\cdot A_a&=3a-4,&c^2&=5,\notag\\
 \delta(a)&\coloneqq5-\frac23\frac{(3a-4)^2}{a^2-4},&
 \tau(a)&\coloneqq\frac{4\delta(a)}{a^2-4}.
 \label{cone:eq:threshold}
\end{align}
Both $c\cdot A_a$ and $\delta(a)$ are positive for $a\in J$.
Consequently these classes lie in the controlled cone. For a metric
$\omega\in2\pi A_a$, its average complex scalar curvature is
$\bar S_a=2(3a-4)/(a^2-4)$ and its volume is
$2\pi^2(a^2-4)$. Thus
\cite[Corollary~1.3 and its proof]{LongTime} applies under
the strict inequality
\begin{equation}\label{cone:eq:startup}
 \widehat{\Ca}(\omega)\coloneqq
 \frac{\Ca(\omega)}{\Vol(X,\omega)}
 =\int_X(S_\omega-\bar S_a)^2\frac{\omega^2}{\int_X\omega^2}
 <\tau(a),
 \qquad\text{equivalently}\qquad
 \Ca(\omega)<8\pi^2\delta(a).
\end{equation}
It gives a global smooth Calabi flow and a uniform Sobolev bound
for each such initial metric.

\subsection{An explicit gap for equivariant integral degenerations}
\label{conealg:sec:gap}

We prove the restricted equivariant gap used in the continuation
argument. The estimate concerns configurations with reduced,
irreducible central fibre. It does not assert full K-stability
or construct an admissible initial metric.

\begin{lemma}\label{conealg:lem:curves}
There is no $G$-fixed point on $X$.  Every $G$-invariant prime curve
on $X$ has class $D=dH-nE$, where $n\geq0$ and $d\geq2n$.
The only such curve of plane degree at most three is the strict
transform $C_0$ of $x^2+y^2+z^2=0$.  It avoids the four blown-up
points and has class $2H$.  All other invariant prime curves have
$d\geq4$.
\end{lemma}
\begin{proof}
The common eigenlines of the diagonal Klein four subgroup are the
three coordinate axes, and coordinate permutations act transitively
on them.  Thus $G$ has no fixed point on $\mathbb P^2$, and the
equivariant blow-down excludes a fixed point on $X$.

No exceptional curve is individually invariant.  An invariant prime
curve therefore descends to an irreducible plane curve.  Its four
point multiplicities agree by transitivity, giving $dH-nE$. Nefness
of $Q$ gives $2d-4n\geq0$.

The defining polynomial of an invariant plane curve spans a
one-dimensional representation of $S_4$.  Every character of $S_4$
is trivial on its Klein four subgroup.  Consequently all exponents
in a monomial of a semi-invariant have the same parity.  In degree
one there is no such monomial.  In degree two the possible
polynomials span $\langle x^2,y^2,z^2\rangle$, whose coordinate
permutation representation has only the invariant one-dimensional
summand $\langle x^2+y^2+z^2\rangle$.  In degree three the only
possibility is a multiple of $xyz$, which is reducible.  The conic
$C_0$ is smooth and its equation takes value $3$ at each $p_i$.
\end{proof}

The following two numerical invariants express the stability test
in terms of divisors: $\beta$ determines the Donaldson--Futaki sign,
whereas $j$ measures the size of the degeneration. Let $A$ be an
ample rational class, and let $F$ be a prime divisor on a smooth
birational model $\pi:Y\to X$. Write $A-xF$ for $\pi^*A-xF$.
Here $\Vol$ denotes algebraic volume, normalized by
$\Vol(A)=A^2$. More explicitly,
\[
 \Vol(A-xF)\coloneqq
 \limsup_{k\to\infty}\frac{2}{k^2}
 h^0\bigl(Y,k\pi^*A-\lfloor kx\rfloor F\bigr),
\]
where $k$ ranges over multiples clearing the denominator of $A$.
The log discrepancy and pseudoeffective threshold are
\[
 A_X(F)\coloneqq1+\operatorname{ord}_F(K_Y-\pi^*K_X),\qquad
 \tau_A(F)\coloneqq\sup\{x\geq0:\Vol(A-xF)>0\}.
\]
With $c=-K_X$, we define
\begin{align}
 \beta_A(F)&\coloneqq A_X(F)A^2+
       2\frac{c\cdot A}{A^2}\int_0^{\tau_A(F)}\Vol(A-xF)\,dx
       +\int_0^{\tau_A(F)}\Vol'(A-xF)\cdot K_X\,dx,
       \label{conealg:eq:beta}\\
 j_A(F)&\coloneqq A^2\tau_A(F)-
       \int_0^{\tau_A(F)}\Vol(A-xF)\,dx.
       \label{conealg:eq:j}
\end{align}
The derivative in \eqref{conealg:eq:beta} is in the canonical
direction, rather than with respect to $x$:
\[
 \Vol'(A-xF)\cdot K_X
 \coloneqq\left.\frac{d}{dt}\right|_{t=0}
 \Vol(\pi^*A-xF+t\pi^*K_X).
\]
These definitions are \cite[Definition~2.8 and the formula preceding
Proposition~2.11]{DL}. On a surface, if
$\pi^*A-xF=P_x+N_x$ is the Zariski decomposition for
$0\leq x<\tau_A(F)$, then
\[
 \Vol(A-xF)=P_x^2,\qquad
 \Vol'(A-xF)\cdot K_X=2P_x\cdot\pi^*K_X.
\]
Consequently $j_A(F)=\int_0^{\tau_A(F)}(A^2-P_x^2)\,dx$
is the accumulated loss of volume as $xF$ is subtracted. It is
strictly positive for a nontrivial divisorial valuation
\cite[Proposition~2.11]{DL}. For a prime curve $D\subset X$,
$A_X(D)=1$, and the formula used below is
\[
 \beta_A(D)=A^2+
 2\frac{c\cdot A}{A^2}\int_0^{\tau_A(D)}P_x^2\,dx
 -2\int_0^{\tau_A(D)}P_x\cdot c\,dx.
\]
Thus both invariants can be computed by intersection theory on
the successive Zariski chambers.
For $v\coloneqq \kappa\operatorname{ord}_F$, $\kappa>0$, extend both invariants by
positive homogeneity in $v$.  Rescaling the polarization gives
\begin{equation}\label{conealg:eq:scaling}
 \beta_{\ell A}(v)=\ell^2\beta_A(v),\qquad
 j_{\ell A}(v)=\ell^3j_A(v).
\end{equation}

We use the non-Archimedean $J$-functional
$J^{\NA}(\mathcal X,\mathcal L)\coloneqq
\lim_{k\to\infty}k^{-1}(\ell_{\max,k}-w_k/d_k)$
\cite[Definition 2.21]{DL}, where $d_k$, $w_k$, and $\ell_{\max,k}$
are the dimension, total weight, and largest weight of
$H^0(\mathcal X_0,\mathcal L_0^k)$ in the weight convention below.
The polarization $L$ on the general fibre includes all exponents.

\begin{proposition}\label{cone:prop:gap}
For every compact interval $I\subset[3,\infty)$ there is
$\gamma_I>0$ with the following property.  Let $a\in I\cap\Q$ and
let $L$ be an integral polarization with $c_1(L)=\ell A_a$,
where $\ell>0$ includes every
denominator-clearing and test-configuration exponent.  Every
nonproduct normal ample $G$-equivariant test configuration with
general polarized fibre $(X,L)$ and integral central fibre satisfies
\begin{equation}\label{cone:eq:polarized-gap}
 J^{\NA}>0,\qquad \frac{\DF/L^2}{J^{\NA}}\geq\frac{\gamma_I}{4\ell}>0.
\end{equation}
The polarization $L$ includes all exponents, as in
\eqref{real:eq:local-gap}. The weight convention is
$\DF=(a_1b_0-a_0b_1)/a_0$, for the Hilbert and total-weight
polynomials $a_0k^2+a_1k+O(1)$ and $b_0k^3+b_1k^2+O(k)$.
One may take
\begin{equation}\label{cone:eq:gamma}
 \gamma_I\coloneqq \min_{a\in I}\min\left\{
 \frac{9a^2+28a+40}{3a(a+2)^2},\quad
 \frac{3a^2+28a+4}{2(a^2-4)(a+1)}\right\}.
\end{equation}
``Integral'' means that the central scheme is reduced and
irreducible.  The assertion is restricted to these equivariant
test configurations.
\end{proposition}
\begin{proof}
Let $(\mathcal X,\mathcal L)$ be such a configuration.  Its central
prime divisor, using the standard trivialization over $\C^*$,
restricts to a valuation $v$ of $\C(X)$.  The central divisor and
the trivialization are $G$-equivariant, so
$v(g^*f)=v(f)$ for every $g\in G$. A nontrivial restriction is
a positive multiple of a divisorial valuation
\cite[Section~3.3, pp.~376--377]{DL}. Its centre on the
proper surface $X$ is irreducible and $G$-invariant.  It cannot be a
point by Lemma~\ref{conealg:lem:curves}.  Thus its centre is an
invariant prime curve $D$, and the discrete valuation ring at the
generic point of $D$ gives $v=\kappa\operatorname{ord}_D$ with $\kappa>0$.

If $v$ is trivial, the one-valuation filtration formula of
\cite[Lemma 3.13]{DL} assigns a common weight to all nonzero
sections in each degree.  This weight is linear in the degree.
After a scalar change of linearization its Rees algebra is the
ordinary section algebra tensored with $\C[t]$. Relative Proj
recovers the product configuration.  Ampleness permits this
reconstruction after a Veronese.  Thus this case is excluded by
the nonproduct assumption.  Notice also that $\Aut^0(X)$ is trivial: it fixes
each of the finitely many exceptional curves, descends to
$\mathbb P^2$, and fixes the four-point projective frame.

It remains to bound $\beta_{A_a}(D)/j_{A_a}(D)$.  Write
$D=dH-nE$ and $r\coloneqq n/d\in[0,1/2]$.  On the big interval of
$A_a-xD=u_xH-v_xE$, with $u_x\coloneqq a-dx$ and
$v_x\coloneqq 1-nx$, the positive part of
the Zariski decomposition is
\begin{equation}\label{conealg:eq:zariski}
 P_x=\begin{cases}
 u_xH-v_xE,&u_x\geq2v_x\geq0,\\
 u_xH,&u_x\geq0,\ v_x\leq0.
 \end{cases}
\end{equation}
Indeed the first class is $(u_x-2v_x)H+v_xQ$, hence nef. In the
second case the negative part is $(-v_x)E$, a sum of disjoint
negative curves orthogonal to $u_xH$. There are no other big
chambers: if $u_x<0$, intersection with the nef class $H$ is
negative. If $v_x>0$ and $u_x<2v_x$, intersection with $Q$ is negative.
In particular $\tau_{A_a}(D)\leq a/d$.  Surface volume and its
directional derivative are $P_x^2$ and $2P_xK_X$, respectively.

Set $h\coloneqq a-2$ and $V_a\coloneqq A_a^2=a^2-4$.
Here $V_a$ is an algebraic intersection number, and the metric
volume in $2\pi A_a$ is $2\pi^2V_a$.
Since $D$ is a prime divisor on the
smooth surface, $A_X(D)=1$.  Substituting $y\coloneqq dx$ into
\eqref{conealg:eq:beta} gives
\begin{equation}\label{conealg:eq:F}
 \beta_{A_a}(D)=V_a+\frac{F(a,r)}d,
\end{equation}
where $F(a,r)=F_1(a,r)$ for $0\leq r\leq1/a$ and
$F(a,r)=F_2(a,r)$ for $1/a\leq r\leq1/2$.
Direct polynomial integration yields
\begin{align}
 F_1(a,r)&=\frac{h}{3(a+2)}
 \left[-\frac{8(2a-1)}{1-2r}
       +\frac{h(14-3a)}{(1-2r)^2}\right],
 &&0\leq r\leq1/a,\label{conealg:eq:F1}\\
 F_2(a,r)&=-\frac{a^2(3a^2+8a-36)}{3V_a}
       +\frac{4(3a^2-6a-4)}{3rV_a},
 &&1/a\leq r\leq1/2.\label{conealg:eq:F2}
\end{align}
For completeness, in the first regime the upper integration
endpoint is $h/(1-2r)$, and the volume and canonical-derivative
integrands are
\[
 (a-y)^2-4(1-ry)^2,\qquad -6(a-y)+8(1-ry).
\]
In the second regime use these integrands up to $y=1/r$, then
$(a-y)^2$ and $-6(a-y)$ up to $y=a$.  Thus
\eqref{conealg:eq:F1} includes $r=0$ without a limiting convention.
Both expressions take the value $-a(3a-4)/3$ at $r=1/a$.

Let $P(a)\coloneqq 3a^2-6a-4$, which is positive for $a\geq3$, and write
$w\coloneqq 1-2r$.  Differentiation gives
\[
 \partial_rF_1=
 \frac{4h}{3(a+2)w^3}
 \bigl[h(14-3a)-4(2a-1)w\bigr]<0,
 \qquad
 \partial_rF_2=-\frac{4P(a)}{3V_a r^2}<0.
\]
For the first inequality, $w\geq h/a$ makes the bracket at most
$-hP(a)/a<0$.  Consequently
\begin{equation}\label{conealg:eq:Fmin}
 F(a,r)\geq F(a,1/2)=-D(a),\qquad
 D(a)\coloneqq \frac{(a-2)(3a^2+20a+8)}{3(a+2)}.
\end{equation}

Except for $C_0$, Lemma~\ref{conealg:lem:curves} gives $d\geq4$.
Equations~\eqref{conealg:eq:F}--\eqref{conealg:eq:Fmin} and the
threshold estimate therefore imply
\[
 \beta_{A_a}(D)\geq V_a-\frac{D(a)}d
 \geq\frac{(a-2)(9a^2+28a+40)}{12(a+2)}>0,
 \qquad 0<j_{A_a}(D)\leq\frac{aV_a}{d}.
\]
Here positivity of $j$ follows, for example, from strict decrease
of the volume for sufficiently small positive $x$, since
$A_aD>0$.  Retaining $d$ when taking the ratio gives
\begin{equation}\label{conealg:eq:large-degree-gap}
 \frac{\beta_{A_a}(D)}{j_{A_a}(D)}
 \geq\frac{dV_a-D(a)}{aV_a}
 \geq\frac{9a^2+28a+40}{3a(a+2)^2}.
\end{equation}

For $C_0=2H$, the threshold is $(a-2)/2$ and the positive part
is $(a-2x)H-E$ throughout.  Integrating gives
\begin{align*}
 \int_0^\tau\Vol(A_a-xC_0)\,dx&=\frac{a^3-12a+16}{6},&
 \int_0^\tau2P_xc\,dx&=\frac{3a^2}{2}-4a+2,\\
 j_{A_a}(C_0)&=\frac{(a-2)^2(a+1)}3,&
 \beta_{A_a}(C_0)&=\frac{(a-2)(3a^2+28a+4)}{6(a+2)}.
\end{align*}
Hence
\begin{equation}\label{conealg:eq:conic-gap}
 \frac{\beta_{A_a}(C_0)}{j_{A_a}(C_0)}
 =\frac{3a^2+28a+4}{2(a^2-4)(a+1)}.
\end{equation}

Finally, we check the conversion from these divisorial invariants
to our normalization of $\DF$. The filtration formula in
\cite[Lemma~3.13]{DL} shows that, after a scalar change of
linearization, the weights for $v=\operatorname{ord}_D$ satisfy
\[
 w(k)=\sum_{j\geq1}h^0(X,kL-jD).
\]
Such a scalar change adds a constant multiple of $kh(k)$ to
$w(k)$ and leaves both $\DF$ and $J^{\NA}$ unchanged.
For this calculation let $P_x$ be the positive part of $L-xD$,
and set
\[
 V_L\coloneqq L^2,\qquad
 I_L\coloneqq\int_0^{\tau_L(D)}P_x^2\,dx,\qquad
 K_L\coloneqq\int_0^{\tau_L(D)}P_x\cdot K_X\,dx.
\]
Riemann--Roch and the coefficient calculation of
\cite[Proposition~3.6 and Lemmas~3.7--3.8]{DL} give
\[
 \begin{aligned}
 a_0&=\frac{V_L}{2},& a_1&=\frac{c\cdot L}{2},\\
 b_0&=\frac{I_L}{2},&
 b_1&=-\frac12\int_0^{\tau_L(D)}P_x\cdot(K_X+D)\,dx
      =-\frac{K_L}{2}-\frac{V_L}{4}.
 \end{aligned}
\]
The last equality follows by integrating
$\frac{d}{dx}\Vol(L-xD)=-2P_x\cdot D$ from $0$ to
$\tau_L(D)$, where the endpoint volumes are $V_L$ and $0$.
Substitution in $\DF=(a_1b_0-a_0b_1)/a_0$ yields
\[
 4\DF=V_L+2\frac{c\cdot L}{V_L}I_L+2K_L=\beta_L(D).
\]
For $v=\kappa\operatorname{ord}_D$, the shifted weights and both
divisorial invariants are multiplied by $\kappa$. Together with
\cite[Lemma~3.16]{DL}, this proves
\begin{equation}\label{conealg:eq:DL}
 4\DF(\mathcal X,\mathcal L)=\beta_L(v),\qquad
 L^2J^{\NA}(\mathcal X,\mathcal L)=j_L(v).
\end{equation}
The factor $4$ is fixed by the displayed coefficient calculation
and agrees with \eqref{rd:eq:cohom-df}. The printed factor
$2(n-1)!$ in \cite[Propositions~3.9 and~3.15]{DL} should be
$2n!$ with their Definition~2.17, which uses the same Hilbert-weight
normalization: the integral before the last equality in the proof
of their Proposition~3.9 is $\beta/n$.
In particular,
\[
 \frac{\DF/L^2}{J^{\NA}}=\frac{\beta_L(v)}{4j_L(v)}.
\]
The valuation's positive scale cancels from this ratio. Combining
\eqref{conealg:eq:scaling}, \eqref{conealg:eq:large-degree-gap},
\eqref{conealg:eq:conic-gap} and \eqref{conealg:eq:DL} proves
\eqref{cone:eq:polarized-gap}.  Both expressions in \eqref{cone:eq:gamma} are
strictly positive and continuous on $[3,\infty)$.

Only the nonequivariant identities from \cite{DL} have been used.
The invariance of the central valuation was proved directly. Thus
no connected-group hypothesis is being applied to the finite group
$S_4$.
\end{proof}

\begin{remark}\label{conealg:rem:range-and-normalization}
If a configuration is initially described with exponent $r$ for
$qA_a$, then its actual general-fibre polarization is $\ell A_a$ with
$\ell\coloneqq rq$, and this is the $\ell$ in \eqref{cone:eq:polarized-gap}.
Both functions in \eqref{cone:eq:gamma} decrease on $[3,\infty)$:
their derivatives are, respectively,
\[
 -\frac{9a^3+38a^2+120a+80}{3a^2(a+2)^3},\qquad
 -\frac{3a^4+56a^3+52a^2+32a+96}
 {2((a^2-4)(a+1))^2}.
\]
Thus on $[3,12]$ one may take $\gamma_I=193/910$, giving
$193/(3640\ell)$ in \eqref{cone:eq:polarized-gap}.  This explicit restricted
gap does not control test configurations with reducible central
fibre or without the $G$-action. Reducedness alone is not the
essential restriction. See \cite[Lemma 2.29]{DL}.
\end{remark}

\subsection{Fixed rational classes and the local gap}

\begin{corollary}[A fixed rational class with an admissible start]
\label{cone:cor:rational-start}
Let $a\in[3,a_+)\cap\Q$. If a smooth $G$-invariant form
$\omega_0\in2\pi A_a$ satisfies \eqref{cone:eq:startup}, its Calabi
flow has a smooth cscK geometric limit on $X$ in the original class.
The geometric limit is unique up to $G$-equivariant
holomorphic isometry.
\end{corollary}
\begin{proof}
Choose an integer $q>0$ so that $qA_a$ is an integral $G$-linearized
class and rescale by $\widetilde\omega(t)\coloneqq q\omega(t/q^2)$,
as in \cite[Lemma~5.1]{Polarized}. The strict energy inequality and
\cite[Corollary 1.3]{LongTime} give a global flow with bounded
Sobolev constants. Since $X$ is rational and $c\cdot A_a>0$, the
numerical hypotheses hold. Proposition~\ref{cone:prop:gap} verifies
the restricted stability condition in
\cite[Corollary~5.6]{Polarized}, which gives a unique smooth
cscK limit on $X$. Since $\Aut^0(X)=\{1\}$, the identification
is $G$-equivariant. Dividing the limiting metric by $q$ returns
the class $2\pi A_a$.
\end{proof}

\begin{lemma}\label{cone:lem:local-gap}
Every $A_b$, $b\in J$, belongs to $\mathcal S_G$.
\end{lemma}
\begin{proof}
The smallest rational subspace $W_{A_b}$ is contained in
$\operatorname{span}_{\mathbb R}\{H,E\}$. If $b>3$, choose a compact
interval $I\subset(3,\infty)$ with $b$ in its interior. A rational
class in $W_{A_b}$ sufficiently close to $A_b$ has the form
\[
 B=tH-sE=sA_{t/s},\qquad t,s\in\mathbb Q,
 \quad\tfrac12<s<2,\quad t/s\in I.
\]
If $c_1(L)=mB$ is integral, then $ms$ is a positive integer:
it is the negative of the coefficient of any $E_i$ in the integral
cohomology basis. Applying \eqref{cone:eq:polarized-gap} with
$\ell=ms$ gives
\begin{equation}\label{cone:eq:auxiliary-gap}
 \frac{\DF/(mB)^2}{J^{\NA}}
 \geq\frac{\gamma_I}{4ms}
 \geq\frac{\gamma_I}{8m}.
\end{equation}
This is the local gap \eqref{real:eq:local-gap} with constant
$\gamma_{A_b}\coloneqq\gamma_I/4>0$. The integer $m$ includes all
powers needed for a genuine $G$-linearization and for the test configuration, so the
estimate applies to every auxiliary bundle used in closedness.
For $b=3$, the minimal rational subspace is the line
$\mathbb R A_3$. Its nearby rational classes are $sA_3$,
$s\in\mathbb Q$, and the same argument uses $I=\{3\}$.
\end{proof}

\subsection{Existence throughout the real interval}

The following application uses varying-class compactness only for
sequences of cscK metrics. Their Calabi energies are zero. The
uniform Sobolev bounds on compact parameter intervals therefore
come from the positive cohomological gap in
\eqref{cone:eq:threshold}.

\begin{theorem}\label{cone:thm:fourpoint}
For every real $a\in[3,8+10/\sqrt3)$, the class $2\pi A_a$ contains
a unique smooth cscK metric $\omega_a$. It is $G$-invariant, and
$a\mapsto\omega_a$ is smooth on this interval, with a smooth local
extension at its left endpoint.

In particular, $2\pi A_{12}$ contains such a metric. Every smooth
$G$-invariant initial metric in $2\pi A_a$ satisfying
\eqref{cone:eq:startup} has a global Calabi flow whose unique
geometric limit is the holomorphic-isometry class of $\omega_a$.
Non-cscK initial metrics satisfying this inequality exist in every
one of these real classes.
\end{theorem}

\begin{proof}
Since $A_3=c_1(X)$, Tian--Yau's theorem
\cite[Theorem~3.3(ii), p.~192]{TianYauKE} supplies a
K\"ahler--Einstein metric $\omega_3\in2\pi A_3$.
Fixed-class cscK uniqueness and $\Aut^0(X)=\{1\}$ make it unique
and hence $G$-invariant. Its Calabi energy is zero, so it satisfies
\eqref{cone:eq:startup}. Lemma~\ref{cone:lem:local-gap} places the
connected family $\{A_a:a\in J\}$ in $\mathcal S_G$.
Theorem~\ref{real:thm:existence}(2) now gives cscK existence
throughout that family.

For clarity, the open-and-closed step in this application is as
follows. Let $\mathcal A\subset J$ consist of parameters admitting
a smooth $G$-invariant initial form satisfying
\eqref{cone:eq:startup}. It contains $3$. Choose a $G$-invariant
Fubini--Study form on $\mathbb P^2$, normalized in $2\pi H$, and
let $\chi$ be its pullback to $X$. If $\omega$ witnesses
$a\in\mathcal A$, then $\omega+(b-a)\chi$ remains positive for
$b$ close to $a$, and its strict energy inequality persists by
smooth dependence of scalar curvature. Thus $\mathcal A$ is
relatively open.

If $a_j\in\mathcal A$ converges to $b\in J$, the fixed-class
existence criterion gives cscK metrics in $2\pi A_{a_j}$.
On a compact parameter interval containing the tail of this
sequence, $A_a^2$, $c\cdot A_a$ and $\delta(a)$ have the positive
bounds required by Theorem~\ref{vc:thm:compactness}.
Theorem~\ref{vc:thm:closedness}, using
\eqref{cone:eq:auxiliary-gap} at $b$, yields a smooth invariant
cscK metric in the exact class $2\pi A_b$. Its energy is zero, so
$b\in\mathcal A$. Hence $\mathcal A$ is relatively closed and
equals $J$. This step applies static compactness to terminal cscK
metrics. It does not pass to an infinite-time limit uniformly in
the parameter. In particular, an admissible metric at $12$ is a
conclusion of continuation, not an assumption needed to start it.

The local smoothness assertion follows from the scalar-curvature
implicit-function theorem. At $\omega_a$, the derivative in a
K\"ahler potential is the negative Lichnerowicz operator. It is
invertible on zero-mean functions since $H^0(X,T^{1,0}X)=0$.
Thus the background family $\omega_a+(b-a)\chi$ admits a smooth
local family of cscK corrections. Elliptic regularity gives
smoothness in all derivatives. Fixed-class uniqueness patches
these families and also gives their $G$-invariance. This is smooth
dependence of smooth cscK metrics on $X$. 

Finally choose a nonconstant smooth $G$-invariant real function
$f$ on $X$. For small $\varepsilon$, the form
$\omega_{a,\varepsilon}\coloneqq\omega_a+\varepsilon\ddc f$ is
positive and invariant. If $\mathcal L_a$ denotes the
Lichnerowicz operator, scalar-curvature linearization gives
\[
 \widehat{\Ca}(\omega_{a,\varepsilon})
 =\varepsilon^2
 \int_X(\mathcal L_af)^2\frac{\omega_a^2}{\int_X\omega_a^2}
 +O(\varepsilon^3).
\]
The coefficient is positive because the kernel consists of
constants. Therefore small nonzero $\varepsilon$ gives positive
energy below $\tau(a)$. The global-flow and convergence assertions
follow from Theorem~\ref{eqv:thm:real} and
Lemma~\ref{eqv:lem:real-gap}, as in the proof of
Theorem~\ref{real:thm:existence}(1). These flows are nonstationary.
\end{proof}

At $a=12$, the exact numerical data are
\[
 A_{12}^2=140,\qquad c\cdot A_{12}=32,\qquad
 \delta(12)=\frac{13}{105},\qquad
 \tau(12)=\frac{13}{3675}.
\]
After continuation has produced $\omega_{12}$, the perturbation
above can be chosen to satisfy
\[
 0<\widehat{\Ca}(\omega_{12,\varepsilon})
 <\frac{13}{7350}<\frac{13}{3675}.
\]
Since the volume is $280\pi^2$, this yields
$0<\Ca(\omega_{12,\varepsilon})<52\pi^2/105$
and the explicit Sobolev margin
\[
 \eta_0\coloneqq96\pi^2\delta(12)
 -6(2\pi A_{12})^2\widehat{\Ca}(\omega_{12,\varepsilon})
 >\frac{208\pi^2}{35}>0.
\]
To obtain the cscK metric in $2\pi A_{12}$, we apply
Proposition~\ref{cone:prop:gap}, through
Lemma~\ref{cone:lem:local-gap}, to the connected family
$\{A_a:3\leq a\leq12\}$. This verifies the local restricted gap
along the family, which lies in the controlled cone. The
anticanonical K\"ahler--Einstein metric at $a=3$ supplies the
starting metric. Theorem~\ref{real:thm:existence}(2) therefore
gives a smooth cscK metric in $2\pi A_a$ for every
$3\leq a\leq12$, in particular in $2\pi A_{12}$.
Small invariant perturbations of the
resulting metric then satisfy the initial-energy bounds above.

We do not know whether the class $A_{12}$ lies within the
parameter range of Arezzo--Pacard's gluing construction
\cite[Theorem~1.4]{APBlowupII}. Their theorem requires
sufficiently small exceptional coefficients but gives no
explicit threshold that can be compared with $1/12$.

\vskip10pt
Haozhao Li, Institute of Geometry and Physics, and School of Mathematical
Sciences, University of Science and Technology of China, No. 96 Jinzhai
Road, Hefei, Anhui Province, 230026, China.
Email: \href{mailto:hzli@ustc.edu.cn}{hzli@ustc.edu.cn}\par\medskip

Bing Wang, Institute of Geometry and Physics, and School of Mathematical
Sciences, University of Science and Technology of China, Hefei 230026,
China, and Hefei National Laboratory, Hefei 230088, China.
Email: \href{mailto:topspin@ustc.edu.cn}{topspin@ustc.edu.cn}\par\medskip

Kai Zheng, University of Chinese Academy of Sciences, Beijing, 100190,
China.\\
Email: \href{mailto:KaiZheng@amss.ac.cn}{KaiZheng@amss.ac.cn}

\end{document}